\pdfoutput=1
\documentclass{amsart}

\usepackage[letterpaper,margin=1in]{geometry}
\usepackage{amsaddr}
\usepackage[T1]{fontenc}
\usepackage{lmodern}
\usepackage{newpxtext}
\usepackage{newpxmath}
\DeclareEncodingSubset{TS1}{zpltlf}{0}
\DeclareEncodingSubset{TS1}{zpllf}{0}
\usepackage{microtype}
\usepackage{graphicx}
\usepackage[inline]{enumitem}
\usepackage{colortbl}
\usepackage{booktabs}
\usepackage{cite}
\usepackage{algorithm}
\usepackage{algpseudocode}
\usepackage{xcolor}
\usepackage{siunitx}
\definecolor{tuftsblue}{HTML}{3E8EDE}
\usepackage[colorlinks,allcolors=tuftsblue]{hyperref}
\usepackage{cleveref}

\crefformat{equation}{(#2#1#3)}
\Crefformat{equation}{(#2#1#3)}
\crefrangeformat{equation}{(#3#1#4)--(#5#2#6)}
\Crefrangeformat{equation}{(#3#1#4)--(#5#2#6)}
\crefmultiformat{equation}{(#2#1#3)}{ and~(#2#1#3)}{, (#2#1#3)}{, and~(#2#1#3)}
\Crefmultiformat{equation}{(#2#1#3)}{ and~(#2#1#3)}{, (#2#1#3)}{, and~(#2#1#3)}

\makeatletter
\renewcommand{\@captionfont}{\normalfont\small}
\def\paragraph{\@startsection{paragraph}{4}%
  \z@{.5\linespacing\@plus.7\linespacing}{-.5em}%
  {\normalfont\itshape}}
\makeatother

\theoremstyle{plain}
\newtheorem{theorem}{Theorem}
\newtheorem{lemma}[theorem]{Lemma}
\theoremstyle{remark}
\newtheorem{remark}{Remark}

\newcounter{run}
\newcommand{\runid}{\addtocounter{run}{1}\therun}

\newcommand{\igrad}{\ensuremath{\nabla}}
\newcommand{\idiv}{\ensuremath{\nabla\cdot}}
\newcommand{\ilap}{\rotatebox[origin=c]{180}{$\nabla$}}
\newcommand{\fs}[1]{\ensuremath{\mathcal{#1}}}
\newcommand{\scr}[1]{\ensuremath{\mathsf{#1}}}
\newcommand{\defeq}{\ensuremath{\mathrel{\mathop:}=}}
\newcommand{\eqdef}{\ensuremath{=\mathrel{\mathop:}}}

\newcommand{\tabadjust}{\centering\small\renewcommand{\arraystretch}{0.8}}
\definecolor{p1color}{RGB}{200,16,46}
\newcolumntype{R}{>{\columncolor{gray!20}}r}
\newcolumntype{E}{>{\columncolor{gray!20}[\tabcolsep][0pt]}r}
\newcolumntype{F}{>{\columncolor{white}[0pt][\tabcolsep]}r}

\newcommand{\fnum}[1]{\num[round-mode=places,scientific-notation=fixed,round-precision=2]{#1}}
\newcommand{\snum}[1]{\num[scientific-notation=true,round-precision=2,detect-weight=true,detect-family=true,round-mode=places,retain-explicit-plus=true,mode=text,output-exponent-marker=\text{e}]{#1}}

\newcommand{\hands}{\texttt{hands}}
\newcommand{\nirep}{\texttt{nirep}}
\newcommand{\rect}{\texttt{rect}}
\newcommand{\gmres}{{GMRES}}
\newcommand{\ngmres}{{NGMRES}}
\newcommand{\angmres}{{aNGMRES}}
\newcommand{\gangmres}{{GA-NGMRES}}
\newcommand{\rpgd}{{RPGD}}
\newcommand{\andacc}{{AA}}
\newcommand{\aaa}{{aAA}}
\newcommand{\gaaa}{{GA-AA}}
\newcommand{\nk}{{NK}}
\newcommand{\lbfgs}{{L-BFGS}}
\newcommand{\fp}{{FP}}
\newcommand{\gd}{{GD}}

\begin{document}

\title[NGMRES for PDE-Constrained Optimization Governed by Transport]{A Generalized Alternating NGMRES Method for PDE-Constrained Optimization Problems Governed by Transport Equations}

\author{Yunhui He}
\address{Department of Mathematics, University of Houston, Houston, TX, USA}
\email{yhe43@central.uh.edu}

\author{Andreas Mang}
\address{Department of Mathematics, Tufts University, Medford, MA, USA}
\email{andreas.mang@tufts.edu}


\begin{abstract}
We propose a generalized alternating nonlinear generalized minimal residual method (\gangmres) to accelerate first-order optimization schemes for PDE-constrained optimization problems governed by transport equations. Interpreting the update rule of a gradient descent scheme in the Hilbert space induced by the regularization model, i.e., a Sobolev gradient descent, as a fixed-point (\fp) iteration, \gangmres\ periodically alternates \ngmres\ updates with \fp\ steps. It requires only first-order derivatives. We consider advection, mass-preserving transport, and incompressible flows as PDE models; variational regularization of the control variable promotes diffeomorphic transport maps, which we verify numerically. With a line search safeguard, \gangmres\ converges globally. Since its \fp\ map need not be contractive, we extend a recent analysis of \ngmres\ for noncontractive iterations to the alternating scheme. We compare \gangmres\ to this Sobolev gradient descent, Anderson acceleration, a limited-memory BFGS method in the same Hilbert space, and Gauss--Newton--Krylov methods on real-world and synthetic data. \gangmres\ converges where the unaccelerated scheme does not and outperforms Anderson acceleration. On our problems, it does not, in general, outperform the (quasi-)Newton methods, but it closes much of the gap to them; which method is fastest depends on the preconditioner, the regularization parameter, the data, and the stopping tolerance. \gangmres\ is a simple and robust option when moderately accurate solutions suffice.
\end{abstract}

\maketitle

\section{Introduction}

In the present work, we propose a novel acceleration scheme for first-order optimization methods in the context of partial differential equation (PDE) constrained optimization problems governed by transport equations. The control variable of the considered formulations is a smooth, stationary vector field $v$. We use variational regularization models to stipulate adequate smoothness requirements on $v$, promoting that the computed transport map is a diffeomorphism. This requirement comes from the application. In biomedical imaging it is often reasonable to assume that maps between data do not introduce folding or singularities; applications that transfer data such as annotations or segmentations from one image to another require invertible maps. We refer to~\cite{modersitzki2004:numerical, sotiras2013:deformable, christensen2006:introduction, miller2015:hamiltonian, hering2023:learn2reg, chen2025:survey} for a general introduction and applications of diffeomorphic registration. We consider different PDE constraints, modeling \begin{enumerate*}[label=\it(\roman*)]\item value-preserving transport maps governed by an advection equation, \item mass-preserving transport maps governed by a continuity equation, and \item incompressible transport maps governed by Stokes-like systems\end{enumerate*}. The problems considered in this work fall into a class of inverse problems that are infinite-dimensional in principle, highly nonlinear, and ill-posed, leading to large-scale, ill-conditioned inversion operators. Our main contribution is the design and empirical evaluation of numerical schemes to accelerate the convergence of first-order optimization methods. We consider extensions of two schemes: \begin{enumerate*}[label=\it(\roman*)]\item a non-linear generalized minimal residual method (\ngmres)~\cite{washio1997:krylov} and \item Anderson acceleration (\andacc)~\cite{anderson1965:iterative}\end{enumerate*}. We will see that the proposed methodology remains robust across various hyperparameter choices, problem formulations, and data sets. It converges where the first-order scheme it accelerates stalls, and its runtime is nearly competitive with that of quasi-Newton and second-order methods; which method is fastest depends on the preconditioner, the regularization parameter, the data, and the stopping tolerance. We anticipate that the proposed methodology generalizes to numerous other applications.

\subsection{Outline of Method}\label{s:outline}

We consider the inverse problem of estimating a transport map from a data set $m_0 : \bar{\Omega} \to \mathbb{R}$ (the \emph{template}) to $m_1 : \bar{\Omega} \to \mathbb{R}$ (the \emph{reference}). These data are compactly supported on some domain $\Omega \subset \mathbb{R}^d$, $d \in \{1,2,3\}$, with closure $\bar{\Omega} \defeq \Omega \cup \partial\Omega$. We parameterize this transport map by a smooth velocity field $v \in \mathcal{V}_{\scr{ad}}$ (the control variable of our problem). Our applications are mostly in medical imaging; as such, the transported quantities are, in general, image intensities.

Formally, we are given two datasets $m_0 \in C^1(\bar{\Omega})$ and $m_1 \in C^1(\bar{\Omega})$ and seek a stationary velocity field $v \in \fs{V}_{\scr{ad}}$ that satisfies the PDE-constrained optimization problem~\cite{mang2015:inexact, hart2009:optimal}
\begin{subequations}\label{e:varopt}
\begin{align}\label{e:varopt:objective}
\underset{m \in \fs{M}_{\scr{ad}}, \, v \in \fs{V}_{\scr{ad}}}{\operatorname{minimize}}\quad
&
\operatorname{obj}(v) \defeq
\frac{1}{2}\int_{\Omega} (m(x,1) - m_1(x) )^2\, \mathrm{d}x
+ \frac{\alpha}{2}\int_{\Omega} \fs{L} v(x) \cdot v(x)\, \mathrm{d}x
\\
\begin{aligned}
\text{subject to}\,\,\, \\
\end{aligned}
\quad&
\begin{aligned}\label{e:varopt:transport}
\partial_t m(x,t) + \igrad m(x,t) \cdot v(x)  & = 0
&& \text{in}\,\,\Omega \times (0,1],
\\
m(x,0) &=  m_0(x)
&& \text{in}\,\,\Omega.
\end{aligned}
\end{align}
\end{subequations}

The state equation~\cref{e:varopt:transport} models the transport of the image intensities $m_0(x)$ subjected to $v$. The first term of the objective functional in~\cref{e:varopt:objective} is a squared $L^2$-distance that measures the proximity of $m$ at time $t=1$ (terminal state) and $m_1$. The second term is a regularization model that stipulates adequate smoothness requirements on $v$ to ensure that the computed transport map is a diffeomorphism. For simplicity and efficiency, we limit the exposition in the present work to a stationary velocity field $v(x)$; a solver and results for non-stationary $v$ can be found in our past work~\cite{mang2015:inexact, mang2017:lagrangian}. In the context of medical imaging, this problem is referred to as (diffeomorphic) image registration~\cite{modersitzki2004:numerical, sotiras2013:deformable, fischer2008:ill}.

The regularization model in~\cref{e:varopt:objective} is determined by the operator $\fs{L}$. We consider differential operators of Helmholtz type~\cite{beg2005:computing, mang2016:constrained},
\begin{equation}\label{e:regop}
\fs{L} \defeq \fs{A}^s,
\qquad
\fs{A} v \defeq -\beta_{\scr{lap}} \ilap v + \beta_{\scr{id}} v,
\qquad \beta_{\scr{lap}} > 0, \quad \beta_{\scr{id}} \ge 0, \quad s \in \mathbb{N},
\end{equation}

\noindent where $s$ determines the order of the regularization model. Unless stated otherwise, we set $\beta_{\scr{lap}}=1$, $\beta_{\scr{id}}=0$, and $s=2$, so that $\fs{L} = \ilap^2$ is the biharmonic operator and the regularization model is the $H^2$-seminorm $\frac{\alpha}{2}\|\ilap v\|^2_{L^2(\Omega)^d}$. For the incompressible model in~\Cref{s:stokes-flow} we use $s=3$ (tri-harmonic operator).
In~\cref{e:varopt:objective}, $\fs{V}_{\scr{ad}}$ and $\fs{M}_{\scr{ad}}$ denote the admissible sets for the control variable $v$ and the state variable $m$ (the transported image intensities of $m_0$), respectively. We require
\begin{equation}\label{e:admissible}
\fs{V}_{\scr{ad}} \defeq H^s_{\scr{per}}(\Omega)^d, \,\, s > d/2+1,
\qquad\text{and}\qquad
\fs{M}_{\scr{ad}} \defeq C^1(\bar{\Omega}\times[0,1]),
\end{equation}

\noindent where $\scr{per}$ denotes periodic boundary conditions, consistent with the spectral discretization in~\Cref{s:discretization}. For $s > d/2+1$ we have the embedding $\fs{V}_{\scr{ad}} \hookrightarrow C^{1,\nu}(\bar{\Omega})^d$, $0 < \nu < \min\{1,\, s-d/2-1\}$. Consequently, the characteristics of~\cref{e:varopt:transport} are uniquely defined, the associated deformation map is a diffeomorphism, and~\cref{e:varopt:transport} admits a unique classical solution $m \in \fs{M}_{\scr{ad}}$ for $m_0 \in C^1(\bar{\Omega})$. Theoretical results on the well-posedness of related optimal control formulations governed by transport equations can be found in~\cite{borzi2003:optimal, chen2011:image, chen2012:image, barbu2016:optimal, haubner2025:wellposedness}; the analysis in~\cite{chen2011:image} covers data of bounded variation, which includes our smooth data. The transport map is the solution operator of~\cref{e:varopt:transport}; the control $v$ is the velocity in~\cref{e:varopt:transport}, and the intensities of $m_0$ are transported along the characteristics of~\cref{e:varopt:transport}. We call the transport map diffeomorphic if the associated \emph{deformation map} $y$ (the endpoint of these characteristics at $t=1$) is a diffeomorphism. The deformation map does not appear in~\cref{e:varopt}; we compute it from $v$ a~posteriori for visualization and as a diagnostic.

\begin{remark}\label{r:modeling}
The constants in~\cref{e:regop} are dictated by the application. The $H^2$-seminorm in~\cref{e:varopt:objective} is a standard model in diffeomorphic image registration~\cite{mang2015:inexact, mang2016:constrained, beg2005:computing}. For $d=2$, the sufficient condition $s > d/2+1$ for a classical diffeomorphic flow holds for the model in~\Cref{s:stokes-flow} ($s=3$), but not for our default $H^2$ model ($s=2$), which lies at the endpoint of the embedding, nor for models based on an $H^1$-norm or $H^1$-seminorm ($s=1$); we assess diffeomorphic behavior numerically through the computed determinant of the deformation gradient (see \cite{mang2015:inexact, mang2016:constrained}). Since $\beta_{\scr{id}}=0$, $\fs{L}$ has a nontrivial kernel consisting of the constant fields; its numerical treatment is described in~\Cref{s:discretization}. Our framework supports $H^1$-, $H^2$-, and $H^3$-norms and seminorms; we use the $H^2$-seminorm as default for simplicity of presentation. The norms require tuning additional parameters.
\end{remark}

\subsection{Existing Work}

We consider PDE-constrained optimization problems to model inverse transport problems. Our PDE constraints include intensity and mass-preserving transport problems (advection and continuity equation), also accounting for incompressible flows (Stokes-like systems). The considered problem formulations have applications in medical imaging~\cite{lee2010:optimal, mang2015:inexact, hart2009:optimal, burger2018:variational, zapf2025:medical}, computer vision~\cite{borzi2003:optimal, chen2011:image, ruhnau2007:optical, barbu2016:optimal, chen2012:image, borzi2002:optimal, borzi2009:multigrid, haubner2025:wellposedness}, and optimal transport~\cite{benamou2000:computational, benzi2011:preconditioning, mang2017:lagrangian}.

Numerical methods for optimization problems governed by transport equations include first-order~\cite{lee2010:optimal, mang2015:inexact, hart2009:optimal, borzi2003:optimal, chen2011:image, ruhnau2007:optical, chen2012:image, borzi2002:optimal, borzi2009:multigrid} and second-order~\cite{mang2015:inexact, mang2017:lagrangian, mang2017:semilagrangian, mang2019:claire, mang2024:claire} optimization approaches. First-order methods are typically straightforward to implement, but often suffer from slow convergence. In contrast, second-order methods offer the potential for high accuracy and fast convergence in terms of the iteration count. But if implemented naively, they can become computationally prohibitive due to the computational costs associated with inverting the Hessian matrix. We note that our past work on effective numerical methods successfully addressed some of the underlying challenges~\cite{mang2015:inexact,  mang2017:semilagrangian, mang2017:lagrangian, mang2016:distributed, mang2016:constrained, brunn2020:multinode, brunn2021:claire, brunn2021:fast, mang2019:claire, himthani2022:claire}, culminating in high-performance code that allows us to solve the underlying inverse problem in under five seconds on a single graphics processing unit~\cite{mang2024:claire}.

We note that due to significant advancements in machine learning, many modern solvers for the inverse problems considered in our work use deep neural network architectures~\cite{joshi2023:r2net, chen2025:survey, liu2023:geometry, krebs2019:learning, balakrishnan2019:voxelmorph} or rely on automatic differentiation~\cite{hartman2023:elastic, franccois2021:metamorphic}. While these ideas have led to significant advancements, they yield high-throughput methodologies with near real-time capabilities, and often make the implementation straightforward, they also have significant drawbacks compared to adjoint-based optimization algorithms. These drawbacks include questionable generalizability to unseen data, no theoretical guarantees about the quality of the results during inference, massive offline costs for hyperparameter tuning and training, and a lack of interpretability.

In our past work~\cite{mang2015:inexact}, we provided numerical evidence that Newton--Krylov (\nk) methods outperform first-order methods for solving~\cref{e:varopt}. In the present work, we will revisit this point: First-order optimization approaches can be accelerated using specialized techniques to achieve faster convergence. Acceleration closes much, but not all, of the gap to \nk\ methods. A powerful acceleration method we study here is \ngmres. \ngmres\ was first proposed to accelerate nonlinear multigrid~\cite{washio1997:krylov}. In recent years, \ngmres\ has gained increased attention; see~\cite{oosterlee2000:krylov, sterck2012:nonlinear, sterck2021:asymptotic, oosterlee2003:multigrid, sterck2013:steepest, greif2025:ngmres,he25:ngmres,he2025:convergence}. Examples for \ngmres\ accelerated optimization algorithms can be found in~\cite{wang2021:asymptotic, riseth2019:objective, sterck2013:steepest}. The work in~\cite{greif2025:ngmres} presents a convergence analysis for \ngmres\ for linear systems, where the fixed-point (\fp) scheme is given by the Richardson iteration. The work in~\cite{he25:ngmres} considers nonlinear problems with a special type of \fp\ iteration. The recent analysis in~\cite{he2026:ngmresconv} covers general nonlinear \fp\ iterations, contractive and noncontractive, and shows that, up to higher order terms, \ngmres\ improves the linear convergence rate of the underlying \fp\ iteration by the gain of its least-squares problem. Analyses of alternating schemes, in which accelerated steps alternate with plain \fp\ steps, are restricted to linear problems~\cite{he2025:convergence, he2025:generalized}. For optimization problems, \cite{sterck2013:steepest} equips \ngmres\ with a line search safeguard on the accelerated iterate, so that the scheme inherits the global convergence of steepest descent with line search to a stationary point from any initial guess; we adopt this safeguard.

Despite its potential, the application and further development of \ngmres\ has yet to be thoroughly investigated. To the best of our knowledge, there is no application of \ngmres\ to PDE-constrained optimization problems governed by transport equations. The only work we are aware of in this context uses an \andacc\ scheme~\cite{willert2014:leveraging}. \andacc\ has been applied successfully in various other contexts~\cite{walker2011:anderson, anderson2019:comments,peng2018:anderson,an2017:anderson,wang2021:asymptotic} by virtue of its simplicity, ease of implementation, and fast convergence. We investigate both approaches and formulate novel  variants for the solution of the variational models considered in this work. We expect that our findings generalize to other problems.

\subsection{Contributions}\label{s:contributions}

We follow up on our prior work on designing effective numerical methods for optimization problems governed by transport equations~\cite{mang2015:inexact, mang2024:claire, mang2017:semilagrangian, mang2017:lagrangian, mang2016:distributed, mang2016:constrained, brunn2020:multinode, brunn2021:claire, brunn2021:fast, mang2019:claire, mang2023:operator, zhang2021:diffeomorphic}. Our main contributions in this work are:
\begin{enumerate*}[label=\it(\roman*)]
\item We propose a generalized alternating \ngmres\ (\gangmres) method to accelerate the convergence of first-order algorithms for solving large-scale, nonlinear PDE-constrained optimization problems governed by transport equations.
\item We provide a comprehensive empirical analysis of the proposed scheme. Specifically, we study the influence of hyperparameters on performance, investigate convergence as a function of mesh size and vanishing regularization, and evaluate different problem formulations (including advection, mass-preserving transport, and incompressible flows) using data sets of varying complexity. Unlike many existing approaches that infer transport maps from data, our formulations employ variational regularization of the control variable to stipulate adequate smoothness, promoting diffeomorphic transport maps. We report results for both synthetic and real-world examples.
\item We benchmark our \gangmres\ scheme against state-of-the-art first-order, quasi-Newton (\lbfgs), and second-order optimization methods and other acceleration schemes (\andacc-variants).
\item In \Cref{s:convergence}, we state which convergence guarantees apply to \gangmres. With the safeguard of~\cite{sterck2013:steepest}, it converges globally (\Cref{t:global}). A Gauss--Newton estimate based on \Cref{l:contraction} indicates that its \fp\ map is not contractive for five of the six regularization parameters we consider, so that the classical theory is unlikely to apply. We extend the residual bound of~\cite{he2026:ngmresconv} to the alternating and the safeguarded scheme with iteration-dependent step sizes (\Cref{t:onestep}) and prove local r-linear convergence under a condition on the computed linear factors (\Cref{t:local}).
\end{enumerate*}

In summary, our results demonstrate that the proposed method accelerates the convergence of first-order optimization algorithms. The gain is largest for small regularization parameters, for which the unaccelerated scheme does not converge within the iteration budget. The proposed \gangmres\ algorithm outperforms the well-known \andacc\ on the two image registration problems for which we compare them.  Its runtime is competitive with that of \nk\ methods; it is faster than the best \nk\ variant for the larger regularization parameters and slower for the smaller (practical) ones. A limited-memory BFGS method in the Hilbert space induced by the regularization operator is faster than \gangmres\ in every experiment in which we ran it and in which it reached the tolerance. On our brain imaging data, \gangmres\ overtakes the fastest \nk\ variant at a tighter tolerance, but not \lbfgs; below a tolerance of $10^{-2}$ only the \nk\ method converges reliably.

\subsection{Limitations}\label{s:limitations}

We only provide results for a Matlab prototype implementation. Our implementation is limited to the two-dimensional case (i.e., $d=2$). Extending our work to $d=3$ requires more work. The main obstacle is memory. \gangmres\ stores $2(w+1)$ vector fields for a window of size $w$, and \lbfgs\ two vector fields per secant pair. On a $256^3$ grid, a single velocity field occupies about $0.2$\,GB in single precision, so that the history alone requires about $8.5$\,GB for \gangmres\ with $w = 20$ and about $20$\,GB for \lbfgs\ with $50$ secant pairs, a substantial fraction of the memory of a single GPU. Our matrix-free \nk\ method stores no history of iterates and has been implemented in three dimensions on CPUs and GPUs, including multi-GPU settings for problems that exceed the memory of a single device~\cite{mang2019:claire, brunn2020:multinode, brunn2021:fast, mang2024:claire}.

Our results indicate that the proposed method is more sensitive to vanishing regularization parameters than \nk\ methods. Addressing this sensitivity requires additional work. At the default stopping tolerance, \lbfgs\ is faster than \gangmres\ in every experiment in which we ran it, and \nk\ with the zero-velocity preconditioner is faster for small regularization parameters.

Our scheme includes the safeguard of the \ngmres\ optimization algorithm of~\cite{sterck2013:steepest}; the descent test it requires comes at no additional cost, since the gradient at the \fp\ iterate is already needed for the least-squares problem. With it, \gangmres\ converges globally (\Cref{t:global}). Without it, we prove local convergence only under a condition on the computed linear factors that we cannot guarantee a~priori, and global convergence remains open (see the end of \Cref{s:convergence}).

\subsection{Outline of the Paper}

We present the methodology in \Cref{s:methods}. This includes a recapitulation of our strategies to solve and discretize the variational problem as well as the new acceleration schemes considered in this work. We report results in \Cref{s:results} and conclude with \Cref{s:conclusions}.

\section{Methods}\label{s:methods}

In this section, we present the methodology. We revisit the optimality conditions of the problem formulation in \Cref{s:optimality}. We discuss the numerical discretization in \Cref{s:discretization}. We introduce the approaches to solve the variational problem \cref{e:varopt} in \Cref{s:line-search}. In particular, we consider a regularization preconditioned gradient descent (\rpgd) scheme (see \Cref{s:rpgd}; first baseline method), propose several variants of acceleration schemes (see \Cref{s:acceleration}), and recapitulate our \nk\ method (see \Cref{s:nk-method}; second baseline method).

\subsection{Optimality Conditions}\label{s:optimality}

For simplicity of presentation, we limit the exposition for the derivation of the optimality conditions to the formulation in~\cref{e:varopt}. We note that we consider two other formulations in the results section.

To solve the variational problem in \cref{e:varopt}, we use the method of Lagrange multipliers and derive the optimality conditions in the continuum (optimize-then-discretize approach; a discretize-then-optimize approach for related formulations can be found in \cite{mang2017:lagrangian}). We introduce the Lagrange multiplier $\lambda : \Omega \times [0,1] \to \mathbb{R}$ and form the Lagrangian
\begin{equation}\label{e:lagrangian}
\begin{aligned}
\ell(v,m,\lambda) \defeq &
\frac{1}{2}\int_{\Omega} (m(x,1) - m_1(x) )^2\, \mathrm{d} x
+ \frac{\alpha}{2}\int_{\Omega} \fs{L} v(x) \cdot v(x)\, \mathrm{d}x + \int_{\Omega} \big(m(x,0) - m_0(x)\big)\lambda(x,0)\, \mathrm{d}x\\
& + \int_0^1 \int_{\Omega} \big(\partial_t m(x,t) + \igrad m(x,t) \cdot v(x)\big) \lambda(x,t)\, \mathrm{d}x \,\mathrm{d}t.
\end{aligned}
\end{equation}

Taking variation of $\ell$ with respect to $v$ yields the reduced gradient
\begin{equation}\label{e:reduced-gradient}
g(v) \defeq \alpha\fs{L} v(x) + \int_0^1 \lambda(x,t) \igrad m(x,t) \,\mathrm{d}t \quad \text{in}\,\,\Omega,
\end{equation}

\noindent where $m$ and $\lambda$ are found by solving the state and adjoint equations, respectively. The state equation is given by~\cref{e:varopt:transport}. Formally, it is obtained by computing first variations of $\ell$ with respect to $\lambda$. The adjoint equation is obtained by computing variations with respect to $m$; we obtain the final value problem
\begin{equation}\label{e:adjoint}
\begin{aligned}
-\partial_t \lambda(x,t) - \idiv (\lambda(x,t) v(x))  & = 0
&& \text{in}\,\,\Omega \times [0,1),
\\
\lambda(x,1) &=  -(m(x,1) - m_1(x))
&& \text{in}\,\,\Omega.
\end{aligned}
\end{equation}

\noindent We solve for $\lambda$ by integrating \cref{e:adjoint} backward in time. Consequently, the evaluation of~\cref{e:reduced-gradient} necessitates the solution of two PDEs: For a given trial velocity $v$, we have to solve \cref{e:varopt:transport} forward in time to obtain $m$ for all $t \in [0,1]$. Then, given $m$ at time $t = 1$, we solve~\cref{e:adjoint} backward in time to obtain $\lambda$ for all $t \in [0,1]$. Having found $m$ and $\lambda$ given a trial $v$, we can evaluate \cref{e:reduced-gradient}. Notice that every evaluation of the distance functional in~\cref{e:varopt:objective} also requires the solution of~\cref{e:varopt:transport}.

The multiplier $\lambda$ lies in the dual of the space in which we measure the residual of the constraint~\cref{e:varopt:transport}; for $L^2(\Omega\times(0,1))$, this dual can be identified with $L^2(\Omega\times(0,1))$, and the duality pairing takes the integral form used in~\cref{e:lagrangian}. Under the assumptions stated in~\Cref{s:outline} ($v \in \fs{V}_{\scr{ad}}$ with $s > d/2+1$, so that $v \in C^{1,\nu}(\bar{\Omega})^d$, and $m_0, m_1 \in C^1(\bar{\Omega})$), the state equation~\cref{e:varopt:transport} has a unique classical solution, given by the method of characteristics. The adjoint equation~\cref{e:adjoint} is in conservative form; its solution is obtained along the same characteristics but involves the Jacobian determinant of the flow of $v$ and hence $\idiv v$, which is, in general, only H\"older continuous. Hence,~\cref{e:adjoint} has a unique weak solution (equivalently, a solution in the sense of characteristics) that is continuous on $\bar{\Omega}\times[0,1]$ but, in general, not continuously differentiable. This suffices for all terms in~\cref{e:lagrangian} and~\cref{e:reduced-gradient} to be well defined. For the same reason, the solution $\lambda$ of~\cref{e:adjoint} is, in general, not an element of $\fs{M}_{\scr{ad}}$. These statements require $s > d/2+1$; for our default model ($s=2$, $d=2$), which does not satisfy this condition (see \Cref{r:modeling}), they are not guaranteed. The reduced gradient~\cref{e:reduced-gradient} is an element of $\fs{V}_{\scr{ad}}^\ast$; its second term is continuous on $\bar{\Omega}$, whereas the regularization operator maps $\fs{V}_{\scr{ad}}$ to its dual, so that $\alpha\fs{L} v$ is understood in the distributional sense. The inverse of the regularization operator, with its kernel treated as in~\Cref{s:discretization}, maps the reduced gradient back to $\fs{V}_{\scr{ad}}$; it is the preconditioner of the scheme in~\Cref{s:rpgd}.

\subsection{Numerical Discretization}\label{s:discretization}

We subdivide the time interval $[0,1]$ into $n_t \in \mathbb{N}$ cells of size $h_t = 1/n_t$. Integrals are discretized using a trapezoidal rule. We subdivide $\Omega = [0,\omega_1] \times \cdots \times [0,\omega_d] = [0,2\pi]^d \subset\mathbb{R}^d$ into $n_x = (n_1,\ldots, n_d) \in \mathbb{N}^d$ cells of size $h_i = 2 \pi / n_i$ along each spatial direction $x_i$, $i=1,\ldots,d$. We discretize spatial derivatives using a pseudo-spectral method with a Fourier basis~\cite{mang2015:inexact, mang2016:distributed}. That is, we approximate functions $u$ on $\Omega$ as
\[
\textstyle
u(x)
=\sum_{\xi \in \mathbb{Z}^d} \hat{u}_{\xi} \exp\left( i \sum_{j=1}^d 2\pi \xi_j x_j/\omega_j\right)
= \sum_{\xi \in \mathbb{Z}^d} \hat{u}_{\xi} \exp( i \langle \xi, x \rangle)
\]

\noindent with $\omega_j = 2\pi$, $x = (x_1,\ldots,x_d)\in\mathbb{R}^d$, $\xi = (\xi_1,\ldots, \xi_d) \in \mathbb{Z}^d$, $-n_j/2 + 1 \le \xi_j \le n_j/2$. The mapping between the spectral coefficients $\hat{u}_{\xi}$ and $u$ is done using Fast Fourier Transforms (FFTs). The associated regular grid locations are $x_l = 2\pi l \oslash n_x$ with $l = (l_1,\ldots,l_d) \in \mathbb{N}^d$, $0 \le l_i \le n_i - 1$, $i=1,\ldots,d$; $\oslash$ denotes the Hadamard division. Consequently, we can effectively (and for smooth data with high accuracy) apply and invert differential operators (at the cost of two FFTs and a diagonal scaling). By virtue of our model choices, some of the high order differential operators $\mathcal{L}$ have a nontrivial kernel. We modify their symbol on the kernel to make them invertible; that is, we set the spectral entries of $\mathcal{L}$ that are zero to one before computing the inverse. This modification is applied only when we invert $\mathcal{L}$, i.e., in the preconditioner in~\cref{e:searchdir}; the variational problem~\cref{e:varopt}, its optimality conditions, and the reduced gradient are unchanged. It is necessary because we set $\beta_{\scr{id}}=0$ in~\cref{e:regop}, so that the symbol of $\mathcal{L}$ vanishes for $\xi=0$, whereas~\cref{e:searchdir} requires $P^{(k)} \succ 0$. Since $P^{(k)} = \alpha\mathcal{L}$, setting these entries to one preconditions the mean component of the gradient by $\alpha I$, i.e., the mean of $v$ is updated by a gradient step scaled by $1/\alpha$. Using the pseudo-inverse instead (setting the corresponding entries of the inverse to zero) would freeze the mean of $v$ at its initial value, and the iterates would converge to a point at which only the non-constant part of the reduced gradient vanishes; this changes the problem. For any invertible preconditioner, the fixed points of~\cref{e:rpgd} are characterized by $g(v) = 0$, so our choice leaves the computed solution unchanged and affects only the metric and thus the path of the iterates.

We use a semi-Lagrangian scheme for numerical time integration. This time-integrator is unconditionally stable; it uses explicit, second-order Runge--Kutta methods. More details about our space-time discretization can be found in~\cite{mang2017:semilagrangian, mang2016:distributed, mang2015:inexact, mang2016:constrained}.

\subsection{Line Search Methods}\label{s:line-search}

Computing a minimizer for \cref{e:varopt} requires us to solve the nonlinear system
\begin{equation}\label{e:reduced-space-opt}
g(v) = 0
\end{equation}

\noindent for $v$. Here, $g$ denotes the reduced gradient in~\cref{e:reduced-gradient}. We use iterative numerical methods to solve \cref{e:reduced-space-opt}. In particular, we consider line search methods of the general form
\begin{equation}\label{e:linesearch}
v^{(k+1)} = v^{(k)} + \rho^{(k)} s^{(k)}, \quad k = 0,1,2,\ldots,
\end{equation}

\noindent where $\rho^{(k)} > 0$ denotes the line search parameter, $v^{(k)} \in \mathbb{R}^{dn}$, $n = \prod_{i=1}^d n_i$, is the discretized control variable $v$ in lexicographical ordering, $s^{(k)} \in \mathbb{R}^{dn}$ denotes the search direction, and $k \in \mathbb{N}_0$ denotes the iteration number. We globalize this scheme using a backtracking line search subject to the Armijo condition. That is, we accept the step size $\rho > 0$ if
\begin{equation}\label{e:armijo}
\operatorname{obj}(v^{(k)} + \rho s^{(k)}) < \operatorname{obj}(v^{(k)}) + \rho c \langle g(v^{(k)}), s^{(k)} \rangle,
\end{equation}

\noindent where $c > 0$ is set to $10^{-4}$ and $\langle\,\cdot\,,\,\cdot\,\rangle$ denotes the $L^2$ pairing, discretized as the grid sum weighted by $h^d=\prod_{i=1}^d h_i$. At each iteration $k$, we initialize the search with $\rho=1$ and backtrack by multiplying $\rho$ by a factor of $1/2$ until~\cref{e:armijo} holds. For first-order methods we accelerate this scheme by keeping $\rho^{(k)}$ in memory, i.e., we initialize the backtracking with the scale we found at iteration $k-1$. This is based on the empirical observation that $\rho^{(k)}$ on average does not change significantly across iterations $k$. This allows us to significantly reduce the number of objective function evaluations required during backtracking. To ensure that our estimate for $\rho^{(k)}$ is not overly pessimistic, we increase the stored search parameter by a factor of 2 for the next evaluation if~\cref{e:armijo} holds for the first backtracking step. We use this strategy for the baseline method described in \Cref{s:rpgd} and for the \fp\ map that underlies the proposed algorithms introduced in \Cref{s:acceleration}. The line search globalizes the \fp\ step but not the accelerated iterate obtained by the extrapolation in line~\ref{l:ngmres-update} of \Cref{a:NGMRESw}, which is accepted without a descent test; hence, \Cref{a:NGMRESw} and \Cref{a:AAw} are not globalized in the sense of~\cref{e:armijo}. We therefore also consider the safeguard of the \ngmres\ optimization algorithm of~\cite{sterck2013:steepest}, in which the accelerated iterate is accepted only if the direction from the \fp\ iterate to it is a descent direction, a line search is performed along that direction, and the scheme falls back to the \fp\ iterate otherwise. \Cref{a:gangmres-sg} states the resulting procedure, which we use for the results in~\Cref{s:results}.

In general, the search direction in \cref{e:linesearch} is given by
\begin{equation}\label{e:searchdir}
s^{(k)} = - (P^{(k)})^{-1} g(v^{(k)}),
\end{equation}

\noindent where $g(v^{(k)})$ is the discretized analogue of the reduced gradient in~\cref{e:reduced-gradient} at iteration $k$ and $P^{(k)} \succ 0$ is an $nd \times nd$ matrix introduced to improve the convergence. The choice of $P^{(k)}$ determines the optimization approach.

\subsubsection{(Regularization) Preconditioned Gradient Descent (\rpgd)}\label{s:rpgd}

For $P^{(k)} = I_{dn} = \operatorname{diag}(1,\ldots,1) \in \mathbb{R}^{dn,dn}$ the scheme in \cref{e:linesearch} corresponds to a first-order gradient descent (\gd) algorithm. We do not consider this scheme in the present work. An alternative strategy is to use the discretized regularization operator $L \in \mathbb{R}^{dn,dn}$ for $P^{(k)}$, scaled by the regularization parameter, i.e., $P^{(k)} = \alpha L$. We refer to this approach as regularization preconditioned \gd\ (\rpgd) method. The iterative scheme becomes
\begin{equation}\label{e:rpgd}
v^{(k+1)} = v^{(k)} - \rho^{(k)} (\alpha L)^{-1} g(v^{(k)}), \quad k = 0,1,2,\ldots
\end{equation}

The reduced gradient~\cref{e:reduced-gradient} is an element of $\fs{V}_{\scr{ad}}^\ast$ (see~\Cref{s:optimality}), and $\alpha L$, with $L$ the discretization of the operator $\fs{L}$ in~\cref{e:regop}, is the Riesz map associated with the inner product $\langle v, w \rangle_{\fs{V}_{\scr{ad}}} \defeq \langle \alpha\fs{L} v, w \rangle_{L^2(\Omega)^d}$ induced by the regularization model in~\cref{e:varopt:objective}. Where no confusion can arise, we use $\fs{L}$ for both the operator and its discretization. For $\beta_{\scr{id}}=0$, this is a semi-inner product, and the statement holds on the subspace of fields with vanishing mean, on which it is an inner product equivalent to $\|\cdot\|_{H^s}$. In the spectral basis, $\fs{A}$ has the symbol $\beta_{\scr{lap}}|\xi|^2 + \beta_{\scr{id}}$, so $L$ is diagonal with entries $(\beta_{\scr{lap}}|\xi|^2+\beta_{\scr{id}})^s$; for our choice of constants, these reduce to $|\xi|^{2s}$. Hence,~\cref{e:rpgd} is gradient descent on $\fs{V}_{\scr{ad}}$ with the gradient in the metric of the control space, i.e., a Sobolev gradient~\cite{neuberger2010:sobolev}, as search direction, and the updates lie in $\fs{V}_{\scr{ad}}$ by construction. For the treatment of the kernel, see~\Cref{s:discretization}. This approach is well-established; it often exhibits a faster convergence rate than standard \gd. Since we use a spectral discretization, applying the inverse of $L$ only involves two FFTs and a diagonal scaling; it has vanishing costs.

\subsubsection{Proposed Acceleration Schemes}\label{s:acceleration}

We propose new variants of the \ngmres\ method to solve~\cref{e:reduced-space-opt} for $v$. To do so, we view the scheme in~\cref{e:rpgd} as a \fp\ iteration $v^{(k+1)} = q_k(v^{(k)})$. More precisely,
\begin{equation}\label{e:fixedpoint}
v^{(k+1)} = v^{(k)} - \rho^{(k)} (\alpha L)^{-1} g(v^{(k)}) \eqdef q_k(v^{(k)}).
\end{equation}

With this, we define the $k$th residual as $r_k(v^{(k)})= v^{(k)} - q_k(v^{(k)}) = \rho^{(k)}(\alpha L)^{-1}g(v^{(k)})$, reflecting that $r_k(v^{(k)}) \to 0$ as $v^{(k)}$ approaches a stationary point of~\cref{e:varopt}. Conversely, $r_k(v^{(k)}) \to 0$ implies $g(v^{(k)}) \to 0$, as desired (see~\cref{e:reduced-space-opt}), provided that the step sizes $\rho^{(k)}$ remain bounded away from zero.

The step size $\rho^{(k)}$ in~\cref{e:fixedpoint} is determined by the line search of~\Cref{s:line-search} and is, in general, not constant across iterations. Thus,~\cref{e:fixedpoint} defines a sequence of maps $q_k(v) \defeq v - \rho^{(k)}(\alpha L)^{-1}g(v)$ and residuals $r_k(v) \defeq v - q_k(v) = \rho^{(k)}(\alpha L)^{-1}g(v)$, indexed by the iteration $k$ at which they are applied. These maps differ only in the scalar $\rho^{(k)}$; the preconditioned gradient $(\alpha L)^{-1}g(v)$ does not depend on $k$. For any $\rho^{(k)} > 0$, $q_k(v) = v$ if and only if $g(v) = 0$, so the fixed points of~\cref{e:fixedpoint} are the stationary points of~\cref{e:varopt}, irrespective of the step size. We keep the line search because a constant step size would have to be selected a~priori; too large a value causes divergence, and too small a value slows convergence. Since $\rho^{(k)}$ does not change significantly across iterations (see~\Cref{s:line-search}), the sequence of maps is approximately stationary in our experiments.

\begin{algorithm}
\caption{Windowed \ngmres\ with depth $w$: \ngmres($w$) for solving $g(v)=0$.}\label{a:NGMRESw}
\small
\begin{algorithmic}[1]
\State {\bf input:} initial guess $v^{(0)} = 0$, integers $w \geq 0$, $n_{\scr{iter}} > 0$, tolerance $\epsilon_{\scr{rel}} > 0$
\State {\bf initialize:} $k \gets 0$, stop $\gets$ 0
\While {$\neg$ stop}
    \State $w^{(k)} \gets \min\{k,w\}$
    \State ${\{\beta_i^{(k)}\}} \gets \operatorname{arg\,min}_{\{\beta_i\}}\big\|g(q_k(v^{(k)}))+\sum_{i=0}^{w^{(k)}} \beta_i\left(g(q_k(v^{(k)}))-g(v^{(k-i)})\right)\big\|^2$ \label{l:ngmres-lsq}
    \State $v^{(k+1)} \gets q_k(v^{(k)}) + \sum_{i=0}^{w^{(k)}} \beta^{(k)}_i\left(q_k(v^{(k)})- v^{(k-i)}\right)$\label{l:ngmres-update}
    \State stop $\gets \|g(v^{(k+1)})\|_{\infty} \leq \epsilon_{\scr{rel}}\|g(v^{(0)})\|_{\infty}\,\, \lor \,\, k \ge n_{\scr{iter}}$\label{l:ngmres-stop}
    \State $k \gets k + 1$
\EndWhile
\State {\bf output:} $v^{(k)}$
\end{algorithmic}
\end{algorithm}

In practice, the \fp\ iteration \cref{e:fixedpoint} might converge slowly or even diverge. We seek an acceleration method to speed up the convergence of \cref{e:fixedpoint}.

A candidate method to accomplish this is the \ngmres\ algorithm presented in \Cref{a:NGMRESw}. The performance of this algorithm is controlled by the hyperparameter $w \in \mathbb{N}$, the depth or window size. In real applications, the choice of the depth $w$ is problem-dependent. In practice, $w$ is typically small to avoid the high computational cost of solving the least-squares problem in line~\ref{l:ngmres-lsq} of \Cref{a:NGMRESw} and to mitigate the risk of rank deficiency. Here, $\|\cdot\|$ denotes the $\ell^2$ norm. We remark that one may consider a different norm for the least-squares problem in line~\ref{l:ngmres-lsq} of \Cref{a:NGMRESw}. We use the relative $\ell^{\infty}$-norm of the gradient $g(v^{(k)})$ at iteration $k$ with a tolerance of $\epsilon_{\scr{rel}}$ as a stopping criterion (see line~\ref{l:ngmres-stop} in \Cref{a:NGMRESw}).

A second candidate to accelerate the convergence of the iterative scheme~\cref{e:fixedpoint} is \andacc; see \Cref{a:AAw}.

\begin{algorithm}
\caption{Windowed \andacc\ with depth $w$: \andacc($w$) for solving $g(v)=0$.} \label{a:AAw}
\small
\begin{algorithmic}[1]
\State {\bf input:} initial guess $v^{(0)} = 0$, integers $w \geq 0$, $n_{\scr{iter}} > 0 $, tolerance $\epsilon_{\scr{rel}} > 0$
\State {\bf initialize:} $k \gets 0$, stop $\gets$ 0
\While {$\neg$ stop}
    \State $w^{(k)} \gets \min\{k,w\}$ and $r_k(v^{(k)}) \gets v^{(k)} - q_k(v^{(k)})$
    \State $\{\mu_i^{(k)}\} \gets \operatorname{arg\,min}_{\{\mu_i\}}\big\|r_k(v^{(k)}) + \sum_{i=1}^{w^{(k)}} \mu_i \left(r_k(v^{(k)})-r_{k-i}(v^{(k-i)})\right)\big\|^2$ \label{l:aa-lsq}
    \State $v^{(k+1)} \gets q_k(v^{(k)}) + \sum_{i=1}^{w^{(k)}} \mu^{(k)}_i\left(q_k(v^{(k)})- q_{k-i}(v^{(k-i)})\right)$\label{l:aa-update}
    \State stop $\gets \|g(v^{(k+1)})\|_{\infty} \leq \epsilon_{\scr{rel}}\|g(v^{(0)})\|_{\infty}\,\, \lor \,\, k \ge n_{\scr{iter}}$
    \State $k \gets k + 1$
\EndWhile
\State {\bf output:} $v^{(k)}$
\end{algorithmic}
\end{algorithm}

At each iteration both algorithms require us to solve a small least-squares problem (see line~\ref{l:ngmres-lsq} in \Cref{a:NGMRESw} and line \ref{l:aa-lsq} in \Cref{a:AAw}, respectively). The difference between these two methods lies in lines \ref{l:aa-lsq} and \ref{l:aa-update}. When \ngmres\ with untruncated depth, i.e., $w=\infty$ or, equivalently, $w^{(k)}=\min\{k,w\}=k$ for every iteration, is applied to solve linear systems using a Richardson iteration with a constant step size, so that $q_k = q$ for all $k$, it has been shown that \ngmres\ generates the same iterates as the classical generalized minimal residual method (\gmres) provided that the norms of the residuals of \gmres\ monotonically decrease~\cite{greif2025:ngmres}. However, in this situation, the iterates $v^{(k+1)}$ generated by \andacc\ can be recovered from the \gmres\ iterates, i.e., $v^{(k+1)}=q(\hat{v}^{(k)})$, where $\hat{v}^{(k)}$ is the iterate of \gmres~\cite{potra2013:characterization,walker2011:anderson}. \ngmres\ and \andacc\ approaches can be treated as a multisecant method~\cite{yang2022:anderson, he2025:convergence}. When applied to nonlinear problems, the behavior of these two methods becomes complex and remains insufficiently understood. To the best of our knowledge, no direct, in-depth quantitative comparison of the performance of \andacc\ and \ngmres\ has been reported in the past. One of the primary objectives of this work is to evaluate and compare the effectiveness of these two acceleration techniques across a range of problems. We solve both least-squares problems by a QR factorization, so that the two schemes differ only in the residual they minimize. For \andacc, we also implemented the normal equations of the least-squares problem, regularized by adding $10^{-6}$ times the identity to the Gram matrix (a Tikhonov shift); this choice makes no difference without the safeguard and a substantial one with it, in favor of QR, which is why we use QR for both methods. Other approaches for solving these least-squares problems are discussed in, e.g.,~\cite{lockhart2022performance,loffeld2016considerations}.

In \Cref{a:NGMRESw} and \Cref{a:AAw}, each iteration requires solving a least-squares problem. This can become computationally expensive, especially for large scale problems, and may result in an ill-conditioned system, especially when the approximations become sufficiently accurate. In this case, the vectors $g(q_k(v^{(k)}))-g(v^{(k-i)})$, $i=0,\ldots,w^{(k)}$ (or $r_k(v^{(k)})-r_{k-i}(v^{(k-i)})$, $i=1,\ldots,w^{(k)}$, for \andacc), tend to exhibit a near-linear dependence. The index ranges differ; at the same depth $w$, the least-squares problem of \ngmres($w$) has $w+1$ unknowns and that of \andacc($w$) has $w$, since \ngmres\ forms its differences with respect to $g(q_k(v^{(k)}))$, which is not one of the stored vectors, whereas \andacc\ forms them with respect to the most recent stored vector $r_k(v^{(k)})$. Both schemes store $w^{(k)}+1$ vectors and thus have the same memory requirement. These vectors differ in their dependence on the step size: \ngmres\ stores $v^{(k-i)}$ and $g(v^{(k-i)})$, which do not depend on it, whereas \andacc\ stores $q_{k-i}(v^{(k-i)})$ and $r_{k-i}(v^{(k-i)}) = \rho^{(k-i)}(\alpha L)^{-1}g(v^{(k-i)})$, so that its least-squares problem combines residuals scaled by different step sizes. The convergence analyses of \andacc\ that we are aware of assume a fixed map (e.g.,~\cite{pollock2021:anderson}). To reduce the number of least-squares problems we have to solve and to mitigate their ill-conditioning, inspired by the generalized alternating Anderson (\andacc--\fp)~\cite{he2025:generalized} and the alternating \ngmres~\cite{he2025:convergence}, where \ngmres\ is applied at periodic intervals within the \fp\ iteration, we propose a generalized alternating \ngmres\ method. The parameters that control the proposed algorithms are the depth $w \in \mathbb{N}$ and the step counts $\sigma\in\mathbb{N}$ and $\tau \in \mathbb{N}_0$, respectively. The proposed method alternates between $\sigma$ steps of \ngmres($w$) and $\tau$ steps of \fp\ iterations, repeating this pattern throughout the iteration process. We denote the proposed method as \gangmres($w;\sigma,\tau$) or simply \gangmres. The proposed numerical scheme is presented in \Cref{a:gangmres-sg}. The same idea can be extended to \andacc, and we denote the corresponding variant as \gaaa($w;\sigma,\tau$); in our implementation, the initial \fp\ step $v^{(1)} = q_0(v^{(0)})$ of \gaaa\ is not counted in the alternation pattern. We note that the work of~\cite{he2025:generalized} employs a reverse ordering in the mixed scheme--specifically, a combination of \fp\ iterations followed by \andacc, denoted by \aaa($w$)[$\sigma$]--\fp[$\tau$]. This scheme can be extended to \ngmres; we denote this approach by \angmres($w$)[$\sigma$]--\fp[$\tau$].

\begin{algorithm}
\caption{Safeguarded generalized alternating \ngmres: \gangmres($w;\sigma,\tau$). This is the procedure used for the results in \Cref{s:results}. Steps~I--III follow the \ngmres\ optimization algorithm of~\cite{sterck2013:steepest}. Omitting Step~III, i.e., accepting $\hat{v}$ unconditionally, gives the unsafeguarded scheme; \Cref{s:results} reports both variants. The numbers of trial steps of the line searches in Step~I and Step~III are denoted by $n_{\scr{ls}}^{(k)} \ge 1$ and $n_{\scr{sg}}^{(k)} \ge 1$, respectively.}
\label{a:gangmres-sg}
\small
\begin{algorithmic}[1]
\State {\bf input:} initial guess $v^{(0)} = 0$, integers $w \in \mathbb{N}$, $\sigma \in \mathbb{N}$, $\tau \in \mathbb{N}_0$, $n_{\scr{iter}} > 0$, tolerances $\epsilon_{\scr{rel}} > 0$, $c > 0$
\State {\bf initialize:} $k \gets 0$, $k_{\scr{res}} \gets 0$, stop $\gets 0$
\While {$\neg$ stop}
    \Statex \hspace{\algorithmicindent}$\triangleright$ \emph{Step I: preconditioning (one \fp\ step; see \cref{e:fixedpoint})}
    \State evaluate $\operatorname{obj}(v^{(k)})$ and $g(v^{(k)})$
    \State stop $\gets \|g(v^{(k)})\|_{\infty} \leq \epsilon_{\scr{rel}}\|g(v^{(0)})\|_{\infty}\,\, \lor \,\, k \ge n_{\scr{iter}}$\label{l:sg-stop}
    \If{stop} \textbf{break} \EndIf
    \State $s^{(k)} \gets -(\alpha L)^{-1} g(v^{(k)})$
    \State $\rho^{(k)} \gets$ step size satisfying \cref{e:armijo} along $s^{(k)}$
    \State $u^{(k)} \gets v^{(k)} + \rho^{(k)} s^{(k)} \eqdef q_k(v^{(k)})$
    \State evaluate $\operatorname{obj}(u^{(k)})$ and $g(u^{(k)})$
    \Statex \hspace{\algorithmicindent}$\triangleright$ \emph{Step II: acceleration}
    \If{$\operatorname{mod}(k,\sigma+\tau) \ge \sigma$}\label{l:sg-switch}
        \State $v^{(k+1)} \gets u^{(k)}$
    \Else
        \State $w^{(k)} \gets \min\{w,\, k - k_{\scr{res}}\}$
        \State ${\{\beta_i^{(k)}\}} \gets \operatorname{arg\,min}_{\{\beta_i\}}\big\|g(u^{(k)})+\sum_{i=0}^{w^{(k)}} \beta_i\left(g(u^{(k)})-g(v^{(k-i)})\right)\big\|^2$\label{l:sg-lsq}
        \State $\hat{v} \gets u^{(k)} + \sum_{i=0}^{w^{(k)}} \beta^{(k)}_i\left(u^{(k)}- v^{(k-i)}\right)$\label{l:sg-cand}
        \Statex \hspace{\algorithmicindent}\hspace{\algorithmicindent}$\triangleright$ \emph{Step III: safeguard}
        \State $\hat{s} \gets \hat{v} - u^{(k)}$
        \If{$\langle g(u^{(k)}), \hat{s}\rangle < 0$}
            \State $\rho_{\scr{sg}}^{(k)} \gets$ step size satisfying \cref{e:armijo} at $u^{(k)}$ along $\hat{s}$
        \Else
            \State $\rho_{\scr{sg}}^{(k)} \gets 0$
        \EndIf
        \If{$\rho_{\scr{sg}}^{(k)} > 0$}
            \State $v^{(k+1)} \gets u^{(k)} + \rho_{\scr{sg}}^{(k)} \hat{s}$
        \Else
            \State $v^{(k+1)} \gets u^{(k)}$, \quad $k_{\scr{res}} \gets k+1$
        \EndIf
    \EndIf
    \State $k \gets k + 1$
\EndWhile
\State {\bf output:} $v^{(k)}$
\end{algorithmic}
\end{algorithm}

\Cref{a:gangmres-sg} makes the cost per iteration explicit. Since we cache the objective and the reduced gradient, a point already evaluated costs no further PDE solve; evaluating the objective costs one PDE solve, and evaluating it together with the reduced gradient costs two. Step~I costs one PDE solve per trial step of the line search plus two for $\operatorname{obj}(u^{(k)})$ and $g(u^{(k)})$; the values at $v^{(k)}$ are cached if $v^{(k)}=u^{(k-1)}$ and cost two PDE solves otherwise. Step~II and the descent test of Step~III require no PDE solve, since $g(u^{(k)})$ is available; if the test succeeds, the line search of Step~III costs one PDE solve per trial step. Charging the evaluation at $v^{(k+1)}$ to the iteration that produces it, a \fp\ step or a fallback to $u^{(k)}$ costs $2+n_{\scr{ls}}^{(k)}$ PDE solves, as does an \rpgd\ step (plus $n_{\scr{sg}}^{(k)}$ if the line search of Step~III was run and failed), and an iteration that accepts the accelerated iterate costs $4+n_{\scr{ls}}^{(k)}+n_{\scr{sg}}^{(k)}$ ($4+n_{\scr{ls}}^{(k)}$ without the safeguard). The PDE counts and runtimes in \Cref{s:results} follow this accounting. \andacc\ does not evaluate the reduced gradient at $q_k(v^{(k)})$, but the safeguard requires it in every \andacc\ step; if the accelerated iterate is rejected, then $v^{(k+1)}=q_k(v^{(k)})$, and the evaluation is reused. The safeguard therefore adds $2+n_{\scr{sg}}^{(k)}$ PDE solves to every \andacc\ step whose accelerated iterate is accepted, $n_{\scr{sg}}^{(k)}$ PDE solves to each step in which the accelerated iterate is rejected following a failed line search in Step~III, and none otherwise.

Of the two orderings, we report results only for \gangmres($w;\sigma,\tau$), which starts each period with \ngmres\ steps. Without the safeguard, it yields on average better results than \angmres($w$)[$\sigma$]--\fp[$\tau$] on our test problems; with the safeguard of \Cref{a:gangmres-sg}, the two orderings perform equally well (see \Cref{s:results:safeguard}), and we report \gangmres\ for definiteness. Starting with \ngmres\ steps is also motivated by the observation that the initial guess can significantly influence \andacc/\ngmres\ performance~\cite{de2024:anderson, he25:ngmres}.

Without Step~III, we observe the following special cases of \gangmres($w;\sigma,\tau$):
\begin{enumerate*}[label=\it(\roman*)]
\item For $\tau = 0$, the method reduces to \ngmres($w$) for any given $\sigma$ (and, formally, to \ngmres($\infty$) for $w = \infty$).
\item For $\sigma=1$ and $\tau \neq 0$, \gangmres($w;\sigma,\tau$) coincides with the alternating \ngmres\ proposed in~\cite{he2025:convergence}.
\end{enumerate*}
With Step~III and $\tau = 0$, \Cref{a:gangmres-sg} is the \ngmres\ optimization algorithm of~\cite{sterck2013:steepest}, including the restart of the window ($w^{(k)}=0$) whenever the accelerated iterate is rejected~\cite{sterck2013:steepest}.

We investigate the performance of \gangmres\ applied to a first-order optimization method for PDE-constrained optimization problems governed by transport equations as a function of the hyperparameters $w, \sigma, \tau$.

\subsubsection{Convergence}\label{s:convergence}

We state the convergence guarantees for \Cref{a:gangmres-sg}; all of them hold for the iteration-dependent maps $q_k$. The proofs are given in \Cref{s:appendix:proofs}, and we verify the results numerically in \Cref{s:appendix:theory}. We equip the discretized control space $\fs{V}$ with the inner product $\langle v,w\rangle_{\fs{V}} \defeq \langle \alpha\tilde{\fs{L}}v, w\rangle_{L^2}$ of the regularization model, where $\tilde{\fs{L}}$ is $\fs{L}$ with its kernel treated as in~\Cref{s:discretization}, so that $(\alpha L)^{-1}$ in~\cref{e:fixedpoint} and \Cref{a:gangmres-sg} stands for $(\alpha\tilde{\fs{L}})^{-1}$, and its dual $\fs{V}^\ast$, which contains the reduced gradient, with the dual norm
\begin{equation}\label{e:dualnorm}
\|g\|_{\fs{V}^\ast} = \sqrt{\langle g, (\alpha\tilde{\fs{L}})^{-1} g\rangle_{L^2}}
                    = \|(\alpha\tilde{\fs{L}})^{-1/2} g\|_{L^2}.
\end{equation}

\noindent In this metric, $-(\alpha\tilde{\fs{L}})^{-1}g$ is the steepest descent direction, and the maps in~\cref{e:fixedpoint} satisfy
\begin{equation}\label{e:qk-norm}
\|q_k(v)-v\|_{\fs{V}} = \rho^{(k)}\|g(v)\|_{\fs{V}^\ast},
\qquad
\langle g(v), q_k(v)-v\rangle_{L^2} = -\rho^{(k)}\|g(v)\|_{\fs{V}^\ast}^2 .
\end{equation}

\noindent The least-squares problem in line~\ref{l:sg-lsq} of \Cref{a:gangmres-sg} can be posed either in the norm defined in~\cref{e:dualnorm}, as recommended in~\cite{he2026:ngmresconv}, or in the Euclidean norm, as done in \Cref{s:results} of this work. \Cref{t:onestep} below covers both choices. \Cref{t:local} holds as stated for the former and also holds for the latter, with $\|\cdot\|_2$ replacing $\|\cdot\|_{\fs{V}}$ and $\|\cdot\|_{\fs{V}^\ast}$ and with constants that depend on the mesh through the equivalence constants of these norms, i.e., through the spectrum of $\alpha\tilde{\fs{L}}$; in three dimensions, the dual norm may be preferable~\cite{he2026:ngmrespicard}. We compare the two variants in \Cref{s:appendix:theory}.

\paragraph{Global Convergence with the Safeguard}\label{s:conv-global}

\begin{theorem}\label{t:global}
Let $c < 1$ in~\cref{e:armijo}, let $\operatorname{obj}$ be bounded below, let $g$ be the gradient of $\operatorname{obj}$ with respect to $\langle\cdot,\cdot\rangle_{L^2}$, and let $g$ be Lipschitz continuous from $\fs{V}$ to $\fs{V}^\ast$ with constant $L_g$ on a convex set containing the iterates and all trial points of the line search of Step~I. Then, for any $w$, $\sigma$, and $\tau$, the step sizes of Step~I satisfy $\rho^{(k)} \ge \min\{1,(1-c)/L_g\}$, and \Cref{a:gangmres-sg} with $\epsilon_{\scr{rel}} = 0$ and $n_{\scr{iter}} = \infty$ either stops at a zero of $g$ or generates iterates that satisfy
\[
\lim_{k\to\infty}\|g(v^{(k)})\|_{\fs{V}^\ast} = 0 .
\]
In particular, the stopping criterion in line~\ref{l:sg-stop} is met after finitely many iterations for every $\epsilon_{\scr{rel}} > 0$.
\end{theorem}

\noindent The theorem makes no assumption on the maps $q_k$, the least-squares coefficients, or the alternation: the safeguard ensures that each iterate decreases the objective at least as much as the \fp\ step, and the line search, which starts from the previous step size, keeps $\rho^{(k)}$ bounded away from zero, which yields $\lim$ instead of $\liminf$ (see~\cite{sterck2013:steepest} for a similar argument under the Wolfe conditions). The theorem concerns the reduced gradient, not the iterates: if the level set $\{v : \operatorname{obj}(v) \le \operatorname{obj}(v^{(0)})\}$ is bounded, every accumulation point of the iterates is a zero of $g$, but convergence of the iterates to a single point requires further assumptions. Since we discretize the optimality conditions in the continuum, the computed $g$ is the gradient of the discrete objective only up to discretization errors. Suppose that the gradient $\nabla\operatorname{obj}$ of the discrete objective with respect to $\langle\cdot,\cdot\rangle_{L^2}$ satisfies the assumptions of \Cref{t:global} on $g$, that $\|\nabla\operatorname{obj}(v^{(k)})-g(v^{(k)})\|_{\fs{V}^\ast} \le \epsilon_{\scr{g}}\|g(v^{(k)})\|_{\fs{V}^\ast}$ for all $k$ with a uniform $\epsilon_{\scr{g}} < 1-c$, and that, as in our implementation, the line search of Step~III falls back to $u^{(k)}$ after finitely many trials. Then the proof carries over with the step-size bound $\rho^{(k)} \ge \min\{1,(1-c-\epsilon_{\scr{g}})/L_g\}$ and yields $\|g(v^{(k)})\|_{\fs{V}^\ast} \to 0$ and thus $\|\nabla\operatorname{obj}(v^{(k)})\|_{\fs{V}^\ast} \to 0$. The consistency of the discretization does not guarantee this relative-error condition near a stationary point, where $\|g\|_{\fs{V}^\ast}$ tends to zero but the discretization error need not.

\paragraph{Contractivity of the Fixed-Point Maps}\label{s:conv-contraction}

Without Step~III, the argument in the proof of \Cref{t:global} does not apply, since the accepted iterate need not decrease the objective.  The classical theory of \ngmres\ and \andacc\ instead assumes that the \fp\ map is a contraction near the solution. We now give a necessary and sufficient condition for this property.

\begin{lemma}\label[lemma]{l:contraction}
Let $g(v^\ast) = 0$, let $g$ be continuously differentiable near $v^\ast$, and let the reduced space Hessian $\fs{H} \defeq g'(v^\ast)$ be self-adjoint and positive definite with respect to $\langle\cdot,\cdot\rangle_{L^2}$. Then $\fs{B} \defeq (\alpha\tilde{\fs{L}})^{-1}\fs{H}$ is self-adjoint and positive definite with respect to $\langle\cdot,\cdot\rangle_{\fs{V}}$. Let $\lambda_{\min}(\fs{B})$ and $\lambda_{\max}(\fs{B})$ denote the smallest and the largest eigenvalue of $\fs{B}$. For $\rho > 0$, the map $q(v) \defeq v - \rho(\alpha\tilde{\fs{L}})^{-1}g(v)$ satisfies the following:
\begin{enumerate*}[label=(\roman*)]
\item if $\rho\lambda_{\max}(\fs{B}) < 2$, then $q$ is a contraction in $\|\cdot\|_{\fs{V}}$ on a neighborhood of $v^\ast$;
\item if $\rho\lambda_{\max}(\fs{B}) \ge 2$, then $q$ is not a contraction on any neighborhood of $v^\ast$, in any norm.
\end{enumerate*}
\end{lemma}

\noindent Applied to $q_k$, \Cref{l:contraction} yields the condition
\begin{equation}\label{e:contraction}
\rho^{(k)}\lambda_{\max}(\fs{B}) < 2,
\end{equation}

\noindent which a contraction-based theory requires to hold uniformly in $k$. At the computed solution, the Gauss--Newton approximation of $g'(v^\ast)$ violates this condition for five of the six regularization parameters considered, with $\rho^{(k)}\lambda_{\max}(\fs{B})$ between $3.1$ and $12.4$ (see \Cref{s:appendix:theory}). Since \Cref{l:contraction} concerns the exact Jacobian $g'(v^\ast)$, and the second-order terms that the Gauss--Newton approximation neglects can change its eigenvalues, this indicates, but does not prove, that the maps $q_k$ are not contractive. The classical theory is thus unlikely to explain the convergence of \gangmres\ that we observe.

\paragraph{Residual Bound and Acceleration}\label{s:conv-bound}
For the analysis, we write every iterate of \Cref{a:gangmres-sg} (\fp, accelerated, or safeguarded) as an affine combination of the \fp\ iterate and the iterates in the window $\fs{W}_k \defeq \{k-w^{(k)},\ldots,k\}$, where $w^{(k)} \defeq 0$ for \fp\ steps,
\begin{equation}\label{e:affine}
v^{(k+1)} = \zeta_{k+1}^{(k)}\,q_k(v^{(k)}) + \sum_{j\in\fs{W}_k}\zeta_j^{(k)}v^{(j)},
\qquad
\zeta_{k+1}^{(k)} + \sum_{j\in\fs{W}_k}\zeta_j^{(k)} = 1 .
\end{equation}

\noindent A \fp\ step, including a fallback in Step~III, has $\zeta_{k+1}^{(k)} = 1$; the accelerated candidate $\hat{v}$ in line~\ref{l:sg-cand}, which is $v^{(k+1)}$ without Step~III, has $\zeta_{k+1}^{(k)} = 1 + \sum_i\beta_i^{(k)}$ and $\zeta_{k-i}^{(k)} = -\beta_i^{(k)}$; and Step~III, whose line search starts from the full step, so that $\rho_{\scr{sg}}^{(k)} \in (0,1]$, gives $v^{(k+1)} = (1-\rho_{\scr{sg}}^{(k)})u^{(k)} + \rho_{\scr{sg}}^{(k)}\hat{v}$, i.e., $\zeta_{k+1}^{(k)} = 1 + \rho_{\scr{sg}}^{(k)}\sum_i\beta_i^{(k)}$ and $\zeta_{k-i}^{(k)} = -\rho_{\scr{sg}}^{(k)}\beta_i^{(k)}$. Following~\cite{he2026:ngmresconv}, we combine the residuals with the same coefficients and define
\begin{equation}\label{e:lsresidual}
\hat{g}^{(k+1)} \defeq \zeta_{k+1}^{(k)}g(q_k(v^{(k)})) + \sum_{j\in\fs{W}_k}\zeta_j^{(k)}g(v^{(j)}),
\end{equation}

\noindent which, for $\hat{v}$, is the residual of the least-squares problem in line~\ref{l:sg-lsq} at its minimizer. We further define the linear factors (whenever the denominators are nonzero)
\begin{equation}\label{e:theta-gamma}
\theta^{(k+1)} \defeq \frac{\|\hat{g}^{(k+1)}\|_{\fs{V}^\ast}}{\|g(v^{(k)})\|_{\fs{V}^\ast}},
\qquad
\gamma^{(k+1)} \defeq \frac{\|\hat{g}^{(k+1)}\|_{\fs{V}^\ast}}{\|g(q_k(v^{(k)}))\|_{\fs{V}^\ast}} .
\end{equation}

\begin{theorem}\label{t:onestep}
Let $g$ be Fr\'echet differentiable with
\begin{equation}\label{e:lipschitz}
\|g'(x)-g'(y)\|_{\fs{V}\to\fs{V}^\ast} \le \sigma_1\|x-y\|_{\fs{V}}
\end{equation}
on a convex set that contains $v^{(j)}$, $j\in\fs{W}_k$, $q_k(v^{(k)})$, and $v^{(k+1)}$.  Then, the iterates generated by \Cref{a:gangmres-sg}, with or without Step~III and with the least-squares problem in line~\ref{l:sg-lsq} posed in the norm defined in~\cref{e:dualnorm} or in the Euclidean norm, satisfy
\begin{subequations}\label{e:ngmres-bound}
\begin{align}
\|g(v^{(k+1)})\|_{\fs{V}^\ast} &\le \theta^{(k+1)}\,\|g(v^{(k)})\|_{\fs{V}^\ast} + \psi^{(k)}, \label{e:ngmres-bound-theta}\\
\|g(v^{(k+1)})\|_{\fs{V}^\ast} &\le \gamma^{(k+1)}\,\|g(q_k(v^{(k)}))\|_{\fs{V}^\ast} + \psi^{(k)}, \label{e:ngmres-bound-gamma}
\end{align}
\end{subequations}
with the higher-order term
\begin{equation}\label{e:hot}
\begin{aligned}
\psi^{(k)} \defeq \frac{\sigma_1}{2}\Big(&|\zeta_{k+1}^{(k)}|\sum_{j\in\fs{W}_k}|\zeta_j^{(k)}|\,\|q_k(v^{(k)})-v^{(j)}\|_{\fs{V}}^2
+ \sum_{i,j\in\fs{W}_k,\,i>j}|\zeta_i^{(k)}||\zeta_j^{(k)}|\,\|v^{(i)}-v^{(j)}\|_{\fs{V}}^2\Big).
\end{aligned}
\end{equation}
If the least-squares problem is posed in~\cref{e:dualnorm}, the norm in which~\cref{e:theta-gamma} measures the residuals, then $\gamma^{(k+1)} \le 1$, with equality for \fp\ steps, and
\begin{equation}\label{e:theta-fp}
\theta^{(k+1)} \le
\begin{cases}
1, & \text{if } v^{(k+1)} = \hat{v},\\[2pt]
\theta_{\scr{fp}}^{(k+1)} \defeq \|g(q_k(v^{(k)}))\|_{\fs{V}^\ast}/\|g(v^{(k)})\|_{\fs{V}^\ast}, & \text{otherwise},
\end{cases}
\end{equation}
with equality in the second case for \fp\ steps.
\end{theorem}
\noindent  \Cref{t:onestep} extends the estimate in the proofs of Theorems~3.2, 3.8, and~3.11 of~\cite{he2026:ngmresconv} to \Cref{a:gangmres-sg}; its proof uses only~\cref{e:affine} and~\cref{e:lipschitz} and therefore covers all steps of \Cref{a:gangmres-sg} and varying step sizes. Inequality~\cref{e:ngmres-bound-gamma} quantifies the acceleration: up to $\psi^{(k)}$, the step reduces the residual of the \fp\ iterate $u^{(k)}$ by the factor $\gamma^{(k+1)}$, which, for $\hat{v}$, is
\begin{equation}\label{e:gamma-angle}
\gamma^{(k+1)} = \sin \angle\Big( g(q_k(v^{(k)})),\;
   \operatorname{span}\big\{\,g(q_k(v^{(k)})) - g(v^{(k-i)})\,\big\}_{i=0}^{w^{(k)}} \Big),
\end{equation}

\noindent with the angle measured in the inner product of $\fs{V}^\ast$. Inequality~\cref{e:ngmres-bound-theta} shows that, to first order, an accelerated step never increases the residual, whereas a \fp\ step with a noncontractive map can ($\theta^{(k+1)}_{\scr{fp}} > 1$), so that, unlike for \ngmres($w$)~\cite{he2026:ngmresconv}, $\theta^{(k+1)} \le 1$ does not hold for every step. Over a period, the linear factor is the product of the $\sigma+\tau$ factors $\theta^{(k+1)}$. For a general noncontractive map, \cite{he2026:ngmresconv} bounds $\psi^{(k)}$ in terms of residuals through all earlier iterates, with constants that grow like $\bar\zeta^{2(k-j)}$ in the distance to iterate $j$, where $\bar\zeta$ bounds the coefficients (Theorem~3.11 of~\cite{he2026:ngmresconv}). Near a nondegenerate solution, only the iterates in the window enter.

\paragraph{Local Convergence}\label{s:conv-local}

\begin{theorem}\label{t:local}
Let $g(v^\ast) = 0$, let $g'(v^\ast)$ be invertible, and let~\cref{e:lipschitz} hold on the ball $\fs{N} \defeq \{v : \|v-v^\ast\|_{\fs{V}} \le r_0\}$ for some $r_0 > 0$. Pose the least-squares problem in the norm defined in~\cref{e:dualnorm}, let $\omega \defeq \sigma+\tau$, and consider iterates of \Cref{a:gangmres-sg} with or without Step~III, that satisfy
\begin{enumerate*}[label=(A\arabic*)]
\item\label{as:rho} $\rho^{(k)} \le \bar\rho$ for all $k$;
\item\label{as:zeta} $|\zeta_{k+1}^{(k)}| + \sum_{j\in\fs{W}_k}|\zeta_j^{(k)}| \le \bar\zeta$ for all $k$;
\item\label{as:period} $\prod_{k=n\omega}^{(n+1)\omega-1}\theta^{(k+1)} \le \bar\theta$ for all $n \ge 0$, for some $\bar\theta < 1$.
\end{enumerate*}
Then, for every $\hat\theta \in (\bar\theta, 1)$, there exist $\delta > 0$ and $C > 0$, which depend on $\hat\theta$, $\bar\theta$, $\bar\rho$, $\bar\zeta$, $\omega$, $w$, $\sigma_1$, $r_0$, $\|g'(v^\ast)\|$, and $\|g'(v^\ast)^{-1}\|$ but not on the iterates, such that $\|v^{(0)}-v^\ast\|_{\fs{V}} \le \delta$ implies
\[
\|g(v^{(k)})\|_{\fs{V}^\ast} \le C\,\hat\theta^{k/\omega}\,\|g(v^{(0)})\|_{\fs{V}^\ast}
\quad\text{and}\quad
\|v^{(k)}-v^\ast\|_{\fs{V}} \le c_E\,C\,\hat\theta^{k/\omega}\,\|g(v^{(0)})\|_{\fs{V}^\ast}
\]
for all $k$, where $c_E \defeq 2\|g'(v^\ast)^{-1}\|_{\fs{V}^\ast\to\fs{V}}$.
\end{theorem}

\noindent \Cref{t:local} does not require contractive maps; a period may contain \fp\ steps that increase the residual, provided that the accelerated steps compensate. The key ingredient is the local error bound $\|v-v^\ast\|_{\fs{V}} \le c_E\|g(v)\|_{\fs{V}^\ast}$ (\Cref{l:errorbound}), which bounds $\psi^{(k)}$ by the residuals in the window (\Cref{l:hot}) and thus takes the role of contractivity. By~\cref{e:qk-norm}, \ref{as:rho} is Assumption~2.1(2) of~\cite{he2026:ngmresconv} for the maps $q_k$, with $\kappa_0 = \bar\rho$. If the maps $q_k$ are contractions in $\|\cdot\|_{\fs{V}}$ near $v^\ast$ with a common constant $\kappa < 1$, then $\theta^{(k+1)} = \gamma^{(k+1)}\theta^{(k+1)}_{\scr{fp}} \le \theta^{(k+1)}_{\scr{fp}} \le \kappa$ near $v^\ast$, since $\|g(v)\|_{\fs{V}^\ast} = \|v-q_k(v)\|_{\fs{V}}/\rho^{(k)}$; the proof of \Cref{t:local} then goes through without~\ref{as:period} and yields any rate above $\kappa$, as~\cite{he25:ngmres} proves for \ngmres($w$) with a fixed contractive map.

\paragraph{Open Questions}\label{s:open-questions}

\Cref{t:local} is conditional, since \ref{as:rho}--\ref{as:period} concern the computed sequence. Assumption~\ref{as:rho} and Assumption~\ref{as:zeta} can be enforced by capping the step size of the line search (which our implementation does not do) and by restarting or filtering the least-squares problem~\cite{he2026:ngmresconv}, respectively, whereas Assumption~\ref{as:period} cannot be enforced in this way, since the gain of the least-squares problem can vanish and the \fp\ steps of a noncontractive map increase the residual. In our experiments, the average factor per period, i.e., the geometric mean of the factors $\theta^{(k+1)}$ of a run raised to the power $\omega$, lies between $0.13$ and $0.61$, but individual periods have products above one for $p=(5,1)$ and $p=(4,2)$ (up to $2.58$ and $1.06$; see \Cref{t:theory-cycle}), so that these runs converge although they do not satisfy~\ref{as:period}. Moreover, the radius $\delta$ given by the proof shrinks like $\bar\zeta^{-2}$, and, for the coefficients we measure (median of $\sum_j|\zeta_j^{(k)}|$ up to about $170$), the higher-order terms are not negligible at our tolerance, so that $\theta^{(k+1)}$ does not predict the observed rate (see \Cref{s:appendix:theory}). Global convergence without Step~III remains open.

\subsubsection{Quasi-Newton Method in the Metric of the Regularization Model}\label{s:lbfgs}

As a baseline, we consider the limited-memory BFGS (\lbfgs) method associated with the inner product $\langle v,w\rangle_{\fs{V}} = \langle \alpha\fs{L}v, w\rangle_{L^2}$, i.e., the quasi-Newton counterpart of~\cref{e:rpgd} (see \Cref{s:rpgd}). Like \rpgd, it uses only first-order derivative information; like \nk, it builds a model of the curvature. We use the standard two-loop recursion~\cite[Alg.~7.4]{nocedal2006:numerical} applied to the raw reduced gradient $g \in \fs{V}^\ast$, with secant pairs $\delta v^{(j)} = v^{(j+1)} - v^{(j)}$ and $\delta g^{(j)} = g(v^{(j+1)}) - g(v^{(j)})$ and the scaled inverse Riesz map $\chi^{(k)}(\alpha\fs{L})^{-1}$ as the initial inverse Hessian, where $\chi^{(k)} = \langle \delta g^{(k-1)}, \delta v^{(k-1)}\rangle / \|\delta g^{(k-1)}\|^2_{\fs{V}^\ast}$ is the analogue of the customary scaling of the identity. Since the gradient differences are kept in $\fs{V}^\ast$, the curvature pairing $\langle \delta g^{(j)}, \delta v^{(j)}\rangle$ is the duality pairing, which is positive whenever the reduced space Hessian is positive definite along the step; we skip pairs that violate this condition~\cite{nocedal2006:numerical}, so $\chi^{(k)}$ is well defined and positive. With an empty memory, the scheme returns the \rpgd\ search direction. We use the Armijo line search of \Cref{s:line-search}, initialized with the unit step at every iteration, as is customary for quasi-Newton methods. The only parameter is the memory $n_{\scr{mem}}$, which we sweep and report in brackets in the tables. The scheme stores $2n_{\scr{mem}}$ fields, comparable to the $2(w+1)$ fields of \gangmres($w$); one iteration costs two PDE solves plus one per trial step of the line search, and the recursion itself requires no PDE solves.

\subsubsection{Newton--Krylov Method}\label{s:nk-method}

The update rule in \cref{e:linesearch} becomes a Newton method if we select $P^{(k)} = H^{(k)}$ in~\cref{e:searchdir}, where $H^{(k)} \in \mathbb{R}^{dn,dn}$ denotes the reduced space Hessian~\cite{mang2015:inexact}. In all experiments, we use the Gauss--Newton approximation of $H^{(k)}$, which neglects the second-order terms of the data mismatch and is positive semidefinite by construction. At every (outer) iteration, we have to solve a large-scale, ill-conditioned linear system $H^{(k)} s^{(k)} = -g(v^{(k)})$ to find the search direction $s^{(k)}$~\cite{mang2015:inexact}. This poses significant computational challenges. To amortize the underlying computational costs and make this approach computationally tractable, we have developed an effective numerical framework~\cite{mang2017:semilagrangian, mang2017:lagrangian, brunn2020:multinode, mang2019:claire, mang2024:claire}. In particular, we have designed a matrix-free, inexact Newton--Krylov (\nk) method for numerical optimization. To achieve optimal performance, we have proposed several strategies to precondition the reduced space Hessian. We consider three variants for preconditioning in this work: a spectral preconditioner (we refer to this approach by \texttt{ireg})~\cite{mang2015:inexact}, a two-level preconditioner (denoted by \texttt{2lrpcsym})~\cite{mang2017:semilagrangian}, and a zero-velocity approximation (denoted by \texttt{h0rpc})~\cite{brunn2020:multinode}. Our numerical approach has been described and evaluated in detail in our past work~\cite{mang2015:inexact, mang2017:semilagrangian, brunn2020:multinode, mang2019:claire, mang2024:claire}.
We report the work of the preconditioner separately in \Cref{s:results}. The spectral preconditioner \texttt{ireg} uses the inverse of the regularization operator, $\mathcal{H}_{\scr{reg}}^{-1} = (\alpha\mathcal{L})^{-1}$, a diagonal scaling in the spectral domain that requires no PDE solves. The zero-velocity variant \texttt{h0rpc} evaluates the Gauss--Newton Hessian at $v = 0$, $\mathcal{H}_{0}\tilde{v} = \alpha\mathcal{L}\tilde{v} + (\igrad m_0 \otimes \igrad m_0)\tilde{v}$, which is a constant operator; applying it requires no PDE solves either, although we invert it with an inner Krylov solve, since our framework is matrix free. The two-level variant \texttt{2lrpcsym} applies the Hessian itself on a coarser grid and therefore solves PDEs at every application. Hence, the columns \#pdes$_{\scr{pc}}$ and pdes$_{\scr{pc}}$ of \Cref{t:h2-nk+rpgd_hands} and \Cref{t:h2-nk+rpgd} are populated for \texttt{2lrpcsym} only, and for this variant the number of PDE solves of the solver alone is a misleading measure of the total work. We note that we demonstrated in our original work~\cite{mang2015:inexact} that our \nk\ algorithm outperforms the \rpgd\ approach described in \Cref{s:rpgd} in terms of convergence and runtime.
\section{Numerical Results}\label{s:results}

We study the performance of the proposed scheme. We include an empirical convergence analysis (see \Cref{s:results:conv}), a study of the performance of our approaches as a function of the regularization parameter $\alpha$ (see \Cref{s:results:regularization}), and experiments that explore mesh convergence (see \Cref{s:results:mesh}). As part of the convergence analysis, we compare \gangmres\ against the baselines (\rpgd, \lbfgs, and \nk) and against \gaaa, examine the effect of the safeguard and of the stopping tolerance, and check the convergence theory of \Cref{s:convergence} numerically (see \Cref{s:results:theory}; details in \Cref{s:appendix:theory}). We conduct experiments for incompressible flows in \Cref{s:stokes-flow} and mass-preserving flows in \Cref{s:optimal-transport}. Several of these studies sweep many hyperparameter settings. For brevity, we do not report every run in detail. We include the individual runs for representative cases; for the others, we summarize the findings in the text or in summary tables, e.g., by the number of runs that reach the tolerance, by median counts, or by the fastest configuration, and we state in each case which runs we omit. In each table, we highlight in red the run(s) that converged the quickest. We describe the hardware and software we use in \Cref{s:appendix:hardware}.

\subsection{Parameter Setting}\label{s:parameters}

To evaluate the performance of the proposed acceleration schemes, we explore a range of hyperparameter choices. We consider $w \in \{1, 5, 10, 15, 20, 25, 50\}$ with $p=(\sigma,\tau)$; we explore various choices for $p$, which we report in the respective tables. Note that when $p=(1, 0)$, \gangmres($w;\sigma,\tau$) reduces to \ngmres($w$).  When $w=2(n_{\scr{iter}})$ and $p=(1,0)$, \gangmres($w;\sigma,\tau$) is \ngmres($\infty$). We also consider \gangmres($\infty;\sigma,\tau$) and \gangmres($w;\sigma,w+1-\sigma$). For the latter, in the linear case, the iterate at step $(w+1)j$ for $j=1,2,\ldots$, is the same as the iterate of restarted \gmres\ (i.e., \gmres($w+1$)) at the same step. The same parameter settings are applied to the \gaaa($w;\sigma,\tau$) method. We consider the \rpgd\ scheme, an \lbfgs\ method, and different variants of the \nk\ scheme as a baseline for comparison. We report exemplary results in \Cref{t:h2-nk+rpgd_hands}.

Runtimes are the median of five repetitions; counters and convergence metrics are deterministic. We use the relative $\ell^\infty$-norm of the reduced gradient as a stopping criterion. The tolerance is $\epsilon_{\scr{rel}} = 5 \times 10^{-2}$. The maximum number of (outer) iterations $n_{\scr{iter}} \in \mathbb{N}$ is set to 200. We add an $\ast$ to the total runtime if a method does not reach the tolerance in $n_{\scr{iter}}$ iterations.

We consider several different datasets. For non-smooth data, we consider biomedical imaging data (the \hands\ dataset~\cite{modersitzki2004:numerical} with native resolution $128\times 128$ and two images taken from the \nirep\ dataset~\cite{christensen2006:introduction} with native resolution $300 \times 300$; the ids are \texttt{na01} and \texttt{na06}) and one synthetic dataset that can be generated in arbitrary resolution (the \rect\ dataset).  We apply a Gaussian smoothing to this data before executing the solver. The standard deviation is set to $\eta{} h_i$, where $h_i = 2\pi/n_i$ is the mesh size and $\eta{} \in \mathbb{N}_0$. We set $\eta{} = 1$ if not stated otherwise. We also consider two synthetic datasets that are smooth.

\begin{table}
\caption{Convergence results for \rpgd, \lbfgs\ and \nk\ for the \hands\ data. The images are of size $128\times128$ (native resolution). The regularization parameter is set to $\alpha = 10^{-3}$. We report the number of (outer) iterations (\#iter), the number of PDE solves performed by the solver (\#pdes), the number of PDE solves performed inside the preconditioner (\#pdes$_{\scr{pc}}$), the number of Hessian matvecs (\#mvs), the relative change of the mismatch (dist), and the relative reduction of the $\ell^\infty$-norm of the gradient (grad). We also report various execution times (accumulative; in seconds). From left to right, we report the time for the evaluation of the PDEs of the solver (pdes; fraction of the total runtime in brackets), the time for the PDE solves inside the preconditioner (pdes$_{\scr{pc}}$; fraction of the total runtime in brackets), the evaluation of the Hessian matvec (mvs; fraction of the total runtime in brackets), and the time to solution (total runtime; tts; runtimes with $\ast$ indicate that the algorithm did not converge before the maximum number of iterations was reached). \#pdes and pdes count the same operations, as do \#pdes$_{\scr{pc}}$ and pdes$_{\scr{pc}}$. The PDE solves inside the Hessian matvecs are included in \#pdes and pdes, so the columns pdes and mvs overlap and must not be added. An entry --- marks a quantity that does not apply: \rpgd\ forms no Hessian matvecs, and the preconditioners \texttt{ireg} and \texttt{h0rpc} involve no PDE solves. The maximum number of iterations is set to 200. Runtimes are the median of five repetitions.}\label{t:h2-nk+rpgd_hands}
\tabadjust
\begin{tabular}{FlrrrrrrRRRE}\toprule
\multicolumn{8}{c}{} & \multicolumn{4}{c}{\cellcolor{gray!20}time (in seconds)} \\
run & method & \#iter & \#pdes & \#pdes$_{\scr{pc}}$ & \#mvs & dist & grad & pdes & pdes$_{\scr{pc}}$ & mvs & tts \\
\midrule
\runid & \rpgd                    & 200 &  732 &  --- &  --- & \snum{8.085751e-02} & \snum{1.278378e-01} & \fnum{4.613946e+00} (\fnum{8.271728e-01}) & ---                                 & ---                                 & $\ast$\fnum{5.577972e+00}\\
\runid & \nk\ (\texttt{ireg})     &   7 &  246 &  --- &  113 & \snum{6.769309e-02} & \snum{3.712912e-02} & \fnum{1.854048e+00} (\fnum{8.286574e-01}) & ---                                 & \fnum{2.014710e+00} (\fnum{9.004642e-01}) & \fnum{2.237413e+00}\\
\runid & \nk\ (\texttt{2lrpcsym}) &   7 &   56 &  250 &   18 & \snum{6.764156e-02} & \snum{4.025264e-02} & \fnum{3.933609e-01} (\fnum{2.870768e-01}) & \fnum{5.339867e-01} (\fnum{3.897062e-01}) & \fnum{3.320612e-01} (\fnum{2.423400e-01}) & \fnum{1.370229e+00}\\
\rowcolor{p1color!20}
\runid & \nk\ (\texttt{h0rpc})    &   8 &  107 &  --- &   42 & \snum{6.848275e-02} & \snum{3.557295e-02} & \fnum{7.755525e-01} (\fnum{6.010085e-01}) & ---                                 & \fnum{7.445407e-01} (\fnum{5.769762e-01}) & \fnum{1.290418e+00}\\
\runid & L-BFGS (5)               &  47 &  154 &  --- &  --- & \snum{6.729986e-02} & \snum{4.832224e-02} & \fnum{1.233855e+00} (\fnum{8.127567e-01}) & ---                                 & ---                                 & \fnum{1.518112e+00}\\
\runid & L-BFGS (10)              &  45 &  146 &  --- &  --- & \snum{6.663732e-02} & \snum{4.842747e-02} & \fnum{1.166819e+00} (\fnum{8.013624e-01}) & ---                                 & ---                                 & \fnum{1.456044e+00}\\
\runid & L-BFGS (20)              &  32 &  105 &  --- &  --- & \snum{6.679515e-02} & \snum{4.980672e-02} & \fnum{9.365879e-01} (\fnum{8.065412e-01}) & ---                                 & ---                                 & \fnum{1.161240e+00}\\
\rowcolor{p1color!20}
\runid & L-BFGS (50)              &  28 &   92 &  --- &  --- & \snum{6.654545e-02} & \snum{4.718971e-02} & \fnum{8.481548e-01} (\fnum{8.165824e-01}) & ---                                 & ---                                 & \fnum{1.038664e+00}\\
\bottomrule
\end{tabular}
\end{table}

\subsection{Convergence and Performance Analysis}\label{s:results:conv}

We compare the convergence of the proposed scheme to several efficient methods for numerical optimization. We report results for \gangmres, \gaaa, \rpgd, \lbfgs, and the \nk\ method. For the latter, we consider three variants of preconditioners for the reduced space Hessian: a spectral preconditioner (\texttt{ireg}), a two-level preconditioner (\texttt{2lrpcsym}), and a zero-velocity approximation (\texttt{h0rpc}). We consider two datasets: \hands\ (resolution: $128\times128$) and \nirep\ (resolution: $300\times300$). We set the regularization parameter $\alpha$ to $10^{-3}$; we found that this choice yields a good agreement between the transported intensities of $m_0(x)$ and $m_1(x)$.

\subsubsection{Comparison against the Baselines}\label{s:results:baselines}

\begin{figure}
\centering
\includegraphics[width=0.6\textwidth]{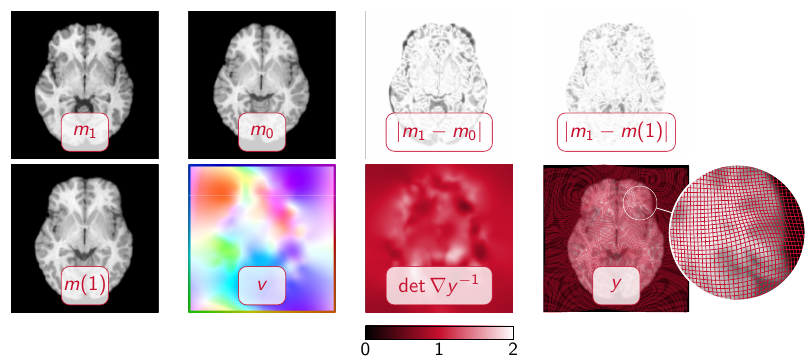}
\caption{We show exemplary results for the baseline model ($H^2$ regularization; compressible velocity). The results correspond to the run with $w=20$, $p=(5,1)$ of the scheme with the safeguard (see \Cref{t:nirep-ngmres-sg}). Top row (from left to right): $(i)$ the reference image $m_1$; $(ii)$ the template image $m_0$ (image to be transported), $(iii)$ the residual differences between $m_0$ and $m_1$ (white: small difference; black: large difference); and $(iv)$ the residual differences between the terminal state $m$ at $t=1$ and $m_1$ after solving for the optimal $v$. Bottom row (from left to right): $(i)$ final state $m$ at $t=1$; $(ii)$ optimal control variable $v$ (color indicates orientation); $(iii)$ determinant of the inverse deformation gradient $\det\igrad y^{-1}$ (the values are all positive at the grid points, numerical evidence of local orientation preservation and of no detected folding; the colormap has fixed limits $[0,2]$ (black: $0$; red: $1$, i.e., volume preserving; white: $2$) and clips values outside this range; since the sampled determinant is strictly positive in all results reported in this work, clipping can only occur at the upper end); and $(iv)$ computed deformation map $y$.}
\label{f:newton-vs-ngmres_nirep}
\end{figure}

\begin{figure}
\centering
\includegraphics[width=0.75\textwidth]{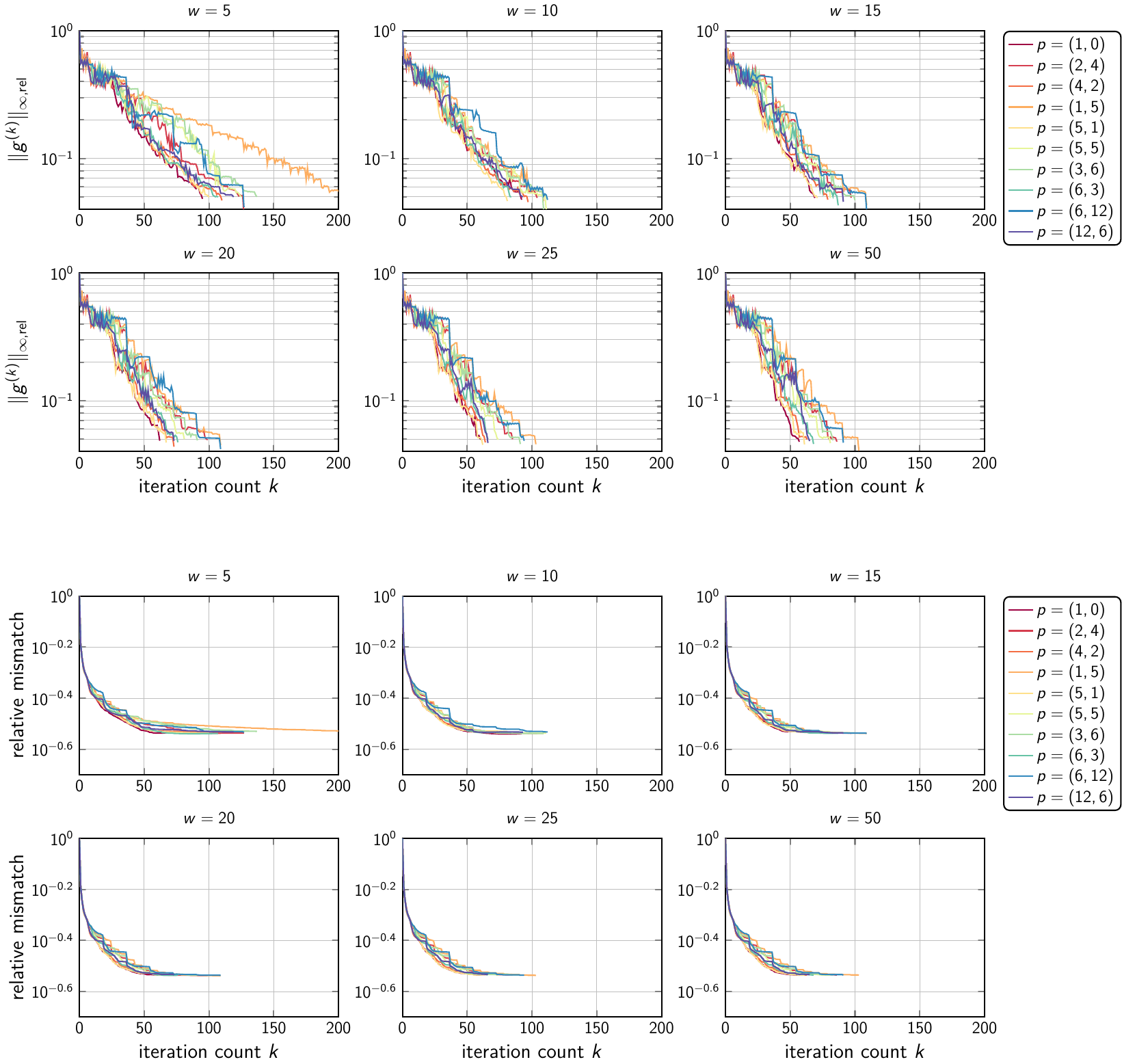}
\caption{Convergence plots for \gangmres\ with the safeguard of the accelerated iterate. We show the relative norm of the gradient $g^{(k)}$ (top block) and the relative mismatch (bottom block) as a function of the iteration count $k$ for the hyperparameters $w$ and $p =(\sigma, \tau)$. The plots show the $60$ runs with $w \in \{5,10,15,20,25,50\}$ and the ten sequences $p$ that are common to these windows; the runs are those of \Cref{t:nirep-ngmres-sg} and \Cref{t:nirep-ngmres-sg_cont}. The safeguard has its most visible effect on the non-alternating sequence $p = (1,0)$ at $w = 25$. Without it, the norm of the gradient stagnates at about $7\times 10^{-2}$ from roughly iteration $75$ until the iteration cap; with it, the run converges in $61$ iterations and no longer separates from the alternating sequences. Without the safeguard, the run with $w = 20$, $p = (5,1)$ also increases the norm of the gradient temporarily by a factor of about $2.5$ before converging in $161$ iterations; with the safeguard, no run shown exceeds its running minimum by more than a factor of $1.5$. In both variants, the histories are not monotone, and the run with $w = 5$, $p = (1,5)$ does not reach the tolerance within the iteration budget.}
\label{f:nirep-300x300-na06-t0-na01-ngmres-conv-sg}
\end{figure}

We report qualitative results in \Cref{f:newton-vs-ngmres_nirep}. To quantify the performance of the considered numerical schemes, we report the number of iterations, the number of PDE solves, the number of Hessian matvecs, the relative change of the data mismatch, and the relative change of the $\ell^\infty$-norm of the reduced gradient. We also report various execution times (in seconds), including the time spent on solving PDEs, Hessian matrix products, the evaluation of $q_k$, the evaluation of the reduced gradient $g$, the time spent on solving the least squares system that appears in the \gaaa\ and \gangmres\ schemes, and the time-to-solution (total runtime). The results for the baseline methods (\rpgd, \lbfgs, and \nk) for the \nirep\ dataset are included in \Cref{t:h2-nk+rpgd}. The results for \gangmres\ are in \Cref{t:nirep-ngmres-sg} and \Cref{t:nirep-ngmres-sg_cont}, with convergence plots in \Cref{f:nirep-300x300-na06-t0-na01-ngmres-conv-sg}. We report the sweeps for the harder \nirep\ problem; for \hands, we quote the corresponding numbers in the text but, for brevity, do not include the individual runs. Unless stated otherwise, all results use the scheme with the safeguard of the accelerated iterate in \Cref{a:gangmres-sg}. The global convergence result of \Cref{s:convergence} (\Cref{t:global}) requires the safeguard, whereas the residual bound and the local convergence result stated there (\Cref{t:onestep,t:local}) apply to both variants; \Cref{s:results:safeguard} summarizes how the two variants differ in practice.

\subsubsection{Effect of the Safeguard}\label{s:results:safeguard}

The safeguard has little effect on \gangmres, so we report only the variant with the safeguard here. On \hands, every run converges in both variants and the median iteration count is $38$ for both; the line search of Step~III costs about $17\%$ more PDE solves in the median ($218$ against $186$). On \nirep, $81$ of the $86$ runs converge within the given budget of 200 iterations without the safeguard and $82$ with it. The median iteration count drops from $100$ to $93$, and the median number of PDE solves is unchanged at $471$ (see \Cref{t:aa-summary}). The safeguard matters in the non-alternating case. For the untruncated window ($w = 400$, $p = (1,0)$), the run without the safeguard stagnates; after $200$ iterations, it has reduced the norm of the gradient only by a factor of about two. With the safeguard, it converges in $57$ iterations (see \Cref{t:nirep-ngmres-sg_cont}). The descent test reuses the reduced gradient that \ngmres\ evaluates anyway, so only the line search adds cost. We therefore enable the safeguard throughout. For brevity, we do not include a table for the individual runs without it.

Almost no \gaaa\ run reaches the tolerance $\epsilon_{\scr{rel}} = 5 \times 10^{-2}$ within $n_{\scr{iter}} = 200$ iterations (see \Cref{s:results:aa}). We provide a more detailed evaluation of \gaaa\ below.

Lastly, we considered two variants for the alternating sequence for our \gangmres\ scheme in line \ref{l:sg-switch} of \Cref{a:gangmres-sg} (variant 1: $\operatorname{mod}(k,\sigma+\tau) \ge \sigma$, i.e., \gangmres($w;\sigma,\tau$); variant 2: $\operatorname{mod}(k,\sigma+\tau)<\tau$, i.e., \angmres($w$)[$\sigma$]--\fp[$\tau$]). All results reported in this work are for variant 1. We chose variant 1 because it is the faster one without the safeguard. With the safeguard of \Cref{a:gangmres-sg}, the two variants both reach the tolerance for $82$ of the $86$ settings and their median iteration counts agree to within $3\%$ ($93$ against $90$). Variant 2 needs about $28\%$ more PDE solves in the median ($601$ against $471$) and is about $18\%$ slower ($20.3$ against $17.3$\,seconds). For brevity, we do not report the individual runs.

\begin{table}
\caption{Convergence results for \rpgd, \lbfgs\ and our \nk\ method for the \nirep\ data. The images are of size $300\times 300$ (native resolution). Setup, reported quantities and layout are as in \Cref{t:h2-nk+rpgd_hands}. Runtimes are the median of five repetitions.}\label{t:h2-nk+rpgd}
\tabadjust
\addtolength{\tabcolsep}{-0.12em}
\begin{tabular}{FlrrrrrrRRRE}\toprule
\multicolumn{8}{c}{} & \multicolumn{4}{c}{\cellcolor{gray!20}time (in seconds)} \\
run & method & \#iter & \#pdes & \#pdes$_{\scr{pc}}$ & \#mvs & dist & grad & pdes & pdes$_{\scr{pc}}$ & mvs & tts \\
\midrule
\runid & \rpgd                    & 200 &  711 &  --- &  --- & \snum{3.517696e-01} & \snum{3.738752e-01} & \fnum{2.075934e+01} (\fnum{8.811661e-01}) & ---                                 & ---                                 & $\ast$\fnum{2.355894e+01}\\
\runid & \nk\ (\texttt{ireg})     &  19 & 1652 &  --- &  798 & \snum{2.886617e-01} & \snum{4.244886e-02} & \fnum{5.518747e+01} (\fnum{8.827094e-01}) & ---                                 & \fnum{5.936521e+01} (\fnum{9.495312e-01}) & \fnum{6.252054e+01}\\
\runid & \nk\ (\texttt{2lrpcsym}) &  19 &  164 & 1510 &   54 & \snum{2.883969e-01} & \snum{4.638026e-02} & \fnum{5.050468e+00} (\fnum{2.321933e-01}) & \fnum{9.690895e+00} (\fnum{4.455352e-01}) & \fnum{4.162867e+00} (\fnum{1.913862e-01}) & \fnum{2.175113e+01}\\
\rowcolor{p1color!20}
\runid & \nk\ (\texttt{h0rpc})    &  15 &  216 &  --- &   86 & \snum{2.922136e-01} & \snum{4.909737e-02} & \fnum{6.444696e+00} (\fnum{6.866354e-01}) & ---                                 & \fnum{6.002716e+00} (\fnum{6.395456e-01}) & \fnum{9.385908e+00}\\
\runid & L-BFGS (5)               & 112 &  348 &  --- &  --- & \snum{2.893164e-01} & \snum{4.912668e-02} & \fnum{1.136580e+01} (\fnum{8.497008e-01}) & ---                                 & ---                                 & \fnum{1.337624e+01}\\
\runid & L-BFGS (10)              & 108 &  336 &  --- &  --- & \snum{2.900116e-01} & \snum{4.907234e-02} & \fnum{1.134975e+01} (\fnum{8.378165e-01}) & ---                                 & ---                                 & \fnum{1.354682e+01}\\
\runid & L-BFGS (20)              &  92 &  291 &  --- &  --- & \snum{2.899805e-01} & \snum{4.860600e-02} & \fnum{1.035247e+01} (\fnum{8.246365e-01}) & ---                                 & ---                                 & \fnum{1.255398e+01}\\
\rowcolor{p1color!20}
\runid & L-BFGS (50)              &  62 &  196 &  --- &  --- & \snum{2.903500e-01} & \snum{4.625843e-02} & \fnum{7.942120e+00} (\fnum{8.266105e-01}) & ---                                 & ---                                 & \fnum{9.608055e+00}\\
\bottomrule
\end{tabular}
\end{table}

\begin{table}
\caption{Convergence results for \gangmres\ with the safeguard of the accelerated iterate (Step~III of \Cref{a:gangmres-sg}) for the \nirep\ data ($300\times300$, native resolution; $\alpha = 10^{-3}$). We report the hyperparameters $w$ and $p = (\sigma,\tau)$, the number of (outer) iterations (\#iter), the number of PDE solves (\#pdes), the relative change of the mismatch (dist), and the relative reduction of the $\ell^\infty$-norm of the gradient (grad). Execution times (accumulative; in seconds) are given for the PDEs (pdes; fraction of the total runtime in brackets), the evaluation of $q_k$, the reduced gradient $g$ at $q_k(v^{(k)})$, the line search of the safeguard (\scr{sg}), the least-squares problem (ls), and the time to solution (tts; $\ast$: tolerance not reached within the iteration budget). Runtimes are the median of five repetitions.}\label{t:nirep-ngmres-sg}
\tabadjust\footnotesize
\begin{tabular}{FrcrrrrRRRRRE}\toprule
\multicolumn{7}{c}{} & \multicolumn{6}{c}{\cellcolor{gray!20}time (in seconds)} \\
run & $w$ & $(\sigma,\tau)$ & \#iter & \#pdes & dist & grad & pdes & $q_k$ & $g$ & sg & ls & tts \\
\midrule
\runid &   1 &   (1,0) & 125 & 854 & \snum{2.933908e-01} & \snum{4.976251e-02} & \fnum{2.365979e+01} (\fnum{8.739769e-01}) & \fnum{1.542602e+01} & \fnum{7.772354e+00} & \fnum{3.267573e+00} & \fnum{8.879504e-02} & \fnum{2.707141e+01}\\
\runid &   5 &         & 95 & 659 & \snum{2.911747e-01} & \snum{4.822407e-02} & \fnum{1.920526e+01} (\fnum{8.640497e-01}) & \fnum{1.282876e+01} & \fnum{5.883538e+00} & \fnum{2.628730e+00} & \fnum{2.530375e-01} & \fnum{2.222703e+01}\\
\runid &  10 &         & 92 & 635 & \snum{2.890279e-01} & \snum{4.741620e-02} & \fnum{1.834579e+01} (\fnum{8.435533e-01}) & \fnum{1.212137e+01} & \fnum{5.704521e+00} & \fnum{2.517997e+00} & \fnum{4.392070e-01} & \fnum{2.174822e+01}\\
\runid &  15 &         & 69 & 485 & \snum{2.932831e-01} & \snum{4.898185e-02} & \fnum{1.476515e+01} (\fnum{8.372058e-01}) & \fnum{1.001054e+01} & \fnum{4.336032e+00} & \fnum{2.013177e+00} & \fnum{4.485292e-01} & \fnum{1.763622e+01}\\
\runid &  20 &         & 62 & 436 & \snum{2.920990e-01} & \snum{4.853183e-02} & \fnum{1.369623e+01} (\fnum{8.146365e-01}) & \fnum{9.401752e+00} & \fnum{3.881426e+00} & \fnum{1.796055e+00} & \fnum{5.359846e-01} & \fnum{1.681269e+01}\\
\runid &  25 &         & 61 & 427 & \snum{2.924669e-01} & \snum{4.963884e-02} & \fnum{1.333944e+01} (\fnum{8.115840e-01}) & \fnum{9.157038e+00} & \fnum{3.791494e+00} & \fnum{1.718272e+00} & \fnum{6.299719e-01} & \fnum{1.643630e+01}\\
\runid &  50 &         & 57 & 392 & \snum{2.913510e-01} & \snum{4.768022e-02} & \fnum{1.228070e+01} (\fnum{7.758737e-01}) & \fnum{8.695890e+00} & \fnum{3.452024e+00} & \fnum{1.419281e+00} & \fnum{8.858065e-01} & \fnum{1.582823e+01}\\
\midrule
\runid &   1 &   (5,1) & 161 & 1016 & \snum{2.949504e-01} & \snum{4.998749e-02} & \fnum{2.707159e+01} (\fnum{8.698142e-01}) & \fnum{1.675214e+01} & \fnum{1.001689e+01} & \fnum{3.639790e+00} & \fnum{9.442267e-02} & \fnum{3.112342e+01}\\
\runid &   5 &         & 100 & 636 & \snum{2.904148e-01} & \snum{4.972881e-02} & \fnum{1.859396e+01} (\fnum{8.630659e-01}) & \fnum{1.209257e+01} & \fnum{6.211759e+00} & \fnum{2.251521e+00} & \fnum{2.209030e-01} & \fnum{2.154408e+01}\\
\runid &  10 &         & 81 & 528 & \snum{2.949264e-01} & \snum{4.627579e-02} & \fnum{1.608876e+01} (\fnum{8.499798e-01}) & \fnum{1.051013e+01} & \fnum{5.258149e+00} & \fnum{2.095278e+00} & \fnum{3.397134e-01} & \fnum{1.892840e+01}\\
\runid &  15 &         & 71 & 457 & \snum{2.933501e-01} & \snum{4.892110e-02} & \fnum{1.445179e+01} (\fnum{8.397717e-01}) & \fnum{9.662917e+00} & \fnum{4.557771e+00} & \fnum{1.643246e+00} & \fnum{4.133544e-01} & \fnum{1.720918e+01}\\
\runid &  20 &         & 67 & 430 & \snum{2.920950e-01} & \snum{4.693546e-02} & \fnum{1.336729e+01} (\fnum{8.199016e-01}) & \fnum{8.998106e+00} & \fnum{4.279242e+00} & \fnum{1.531596e+00} & \fnum{4.971276e-01} & \fnum{1.630353e+01}\\
\runid &  25 &         & 62 & 392 & \snum{2.925284e-01} & \snum{4.512947e-02} & \fnum{1.282343e+01} (\fnum{8.149312e-01}) & \fnum{8.838300e+00} & \fnum{3.986484e+00} & \fnum{1.312692e+00} & \fnum{5.379729e-01} & \fnum{1.573559e+01}\\
\runid &  50 &         & 61 & 384 & \snum{2.908364e-01} & \snum{4.543488e-02} & \fnum{1.231291e+01} (\fnum{7.762904e-01}) & \fnum{8.495848e+00} & \fnum{3.845888e+00} & \fnum{1.204097e+00} & \fnum{8.514792e-01} & \fnum{1.586122e+01}\\
\midrule
\runid &   1 &   (1,5) & 200 & 836 & \snum{2.912587e-01} & \snum{5.735165e-02} & \fnum{2.291685e+01} (\fnum{8.674700e-01}) & \fnum{1.245375e+01} & \fnum{1.235433e+01} & \fnum{9.153700e-01} & \fnum{2.538446e-02} & $\ast$\fnum{2.641803e+01}\\
\runid &   5 &         & 200 & 839 & \snum{2.965101e-01} & \snum{5.511406e-02} & \fnum{2.278227e+01} (\fnum{8.578085e-01}) & \fnum{1.233378e+01} & \fnum{1.228806e+01} & \fnum{9.958197e-01} & \fnum{9.509217e-02} & $\ast$\fnum{2.655868e+01}\\
\runid &  10 &         & 109 & 468 & \snum{2.900441e-01} & \snum{4.614739e-02} & \fnum{1.394774e+01} (\fnum{8.539583e-01}) & \fnum{8.422758e+00} & \fnum{6.608088e+00} & \fnum{5.130915e-01} & \fnum{9.938963e-02} & \fnum{1.633304e+01}\\
\runid &  15 &         & 109 & 469 & \snum{2.901278e-01} & \snum{4.301158e-02} & \fnum{1.403504e+01} (\fnum{8.511750e-01}) & \fnum{8.300231e+00} & \fnum{6.693015e+00} & \fnum{5.431635e-01} & \fnum{1.382196e-01} & \fnum{1.648902e+01}\\
\runid &  20 &         & 109 & 468 & \snum{2.898950e-01} & \snum{4.475914e-02} & \fnum{1.391711e+01} (\fnum{8.397113e-01}) & \fnum{8.149158e+00} & \fnum{6.725588e+00} & \fnum{5.441629e-01} & \fnum{1.884818e-01} & \fnum{1.657368e+01}\\
\runid &  25 &         & 103 & 441 & \snum{2.903558e-01} & \snum{4.558015e-02} & \fnum{1.365260e+01} (\fnum{8.335384e-01}) & \fnum{8.174936e+00} & \fnum{6.395495e+00} & \fnum{5.261463e-01} & \fnum{2.315485e-01} & \fnum{1.637908e+01}\\
\runid &  50 &         & 103 & 443 & \snum{2.912175e-01} & \snum{4.053060e-02} & \fnum{1.336431e+01} (\fnum{7.911981e-01}) & \fnum{7.977668e+00} & \fnum{6.379331e+00} & \fnum{5.266325e-01} & \fnum{4.760682e-01} & \fnum{1.689123e+01}\\
\midrule
\runid &   1 &   (5,5) & 144 & 766 & \snum{2.942997e-01} & \snum{4.966536e-02} & \fnum{2.076932e+01} (\fnum{8.701286e-01}) & \fnum{1.245780e+01} & \fnum{8.829877e+00} & \fnum{2.011261e+00} & \fnum{5.223258e-02} & \fnum{2.386925e+01}\\
\runid &   5 &         & 125 & 687 & \snum{2.949669e-01} & \snum{4.991647e-02} & \fnum{1.904862e+01} (\fnum{8.618982e-01}) & \fnum{1.142065e+01} & \fnum{7.723527e+00} & \fnum{2.146104e+00} & \fnum{1.754148e-01} & \fnum{2.210077e+01}\\
\runid &  10 &         & 111 & 614 & \snum{2.906802e-01} & \snum{3.934893e-02} & \fnum{1.731377e+01} (\fnum{8.515419e-01}) & \fnum{1.054706e+01} & \fnum{6.690013e+00} & \fnum{1.967732e+00} & \fnum{2.685079e-01} & \fnum{2.033225e+01}\\
\runid &  15 &         & 81 & 427 & \snum{2.927835e-01} & \snum{4.890978e-02} & \fnum{1.283255e+01} (\fnum{8.409879e-01}) & \fnum{8.367773e+00} & \fnum{4.885138e+00} & \fnum{9.814049e-01} & \fnum{2.670292e-01} & \fnum{1.525889e+01}\\
\runid &  20 &         & 81 & 427 & \snum{2.925759e-01} & \snum{4.995566e-02} & \fnum{1.323988e+01} (\fnum{8.299312e-01}) & \fnum{8.491623e+00} & \fnum{5.033062e+00} & \fnum{1.014101e+00} & \fnum{3.645588e-01} & \fnum{1.595299e+01}\\
\runid &  25 &         & 74 & 393 & \snum{2.917148e-01} & \snum{4.952882e-02} & \fnum{1.219753e+01} (\fnum{8.183918e-01}) & \fnum{8.048150e+00} & \fnum{4.532487e+00} & \fnum{8.879725e-01} & \fnum{4.265740e-01} & \fnum{1.490427e+01}\\
\runid &  50 &         & 81 & 427 & \snum{2.912181e-01} & \snum{4.648715e-02} & \fnum{1.266779e+01} (\fnum{7.631904e-01}) & \fnum{8.250012e+00} & \fnum{4.883947e+00} & \fnum{9.191178e-01} & \fnum{8.595995e-01} & \fnum{1.659847e+01}\\
\midrule
\runid &   1 &   (4,2) & 159 & 907 & \snum{2.908248e-01} & \snum{4.588129e-02} & \fnum{2.439557e+01} (\fnum{8.691223e-01}) & \fnum{1.491530e+01} & \fnum{9.717491e+00} & \fnum{2.755128e+00} & \fnum{7.594192e-02} & \fnum{2.806920e+01}\\
\runid &   5 &         & 110 & 655 & \snum{2.926808e-01} & \snum{4.683560e-02} & \fnum{1.884530e+01} (\fnum{8.655992e-01}) & \fnum{1.174952e+01} & \fnum{6.845141e+00} & \fnum{2.317400e+00} & \fnum{1.904406e-01} & \fnum{2.177139e+01}\\
\runid &  10 &         & 97 & 569 & \snum{2.903014e-01} & \snum{4.545029e-02} & \fnum{1.722331e+01} (\fnum{8.465040e-01}) & \fnum{1.109950e+01} & \fnum{6.223354e+00} & \fnum{1.801964e+00} & \fnum{3.340892e-01} & \fnum{2.034640e+01}\\
\runid &  15 &         & 79 & 460 & \snum{2.931619e-01} & \snum{4.746359e-02} & \fnum{1.446046e+01} (\fnum{8.378644e-01}) & \fnum{9.653250e+00} & \fnum{5.066249e+00} & \fnum{1.331661e+00} & \fnum{3.762097e-01} & \fnum{1.725871e+01}\\
\runid &  20 &         & 73 & 427 & \snum{2.914395e-01} & \snum{4.357157e-02} & \fnum{1.333176e+01} (\fnum{8.244927e-01}) & \fnum{8.925061e+00} & \fnum{4.651686e+00} & \fnum{1.300690e+00} & \fnum{4.495470e-01} & \fnum{1.616965e+01}\\
\rowcolor{p1color!20}
\runid &  25 &         & 64 & 376 & \snum{2.925754e-01} & \snum{4.707619e-02} & \fnum{1.190423e+01} (\fnum{8.177442e-01}) & \fnum{8.151552e+00} & \fnum{3.989632e+00} & \fnum{1.018161e+00} & \fnum{4.766172e-01} & \fnum{1.455740e+01}\\
\runid &  50 &         & 63 & 365 & \snum{2.912775e-01} & \snum{4.965446e-02} & \fnum{1.189602e+01} (\fnum{7.777289e-01}) & \fnum{8.153250e+00} & \fnum{3.907593e+00} & \fnum{9.717626e-01} & \fnum{7.754630e-01} & \fnum{1.529585e+01}\\
\midrule
\runid &   1 &   (2,4) & 200 & 946 & \snum{2.996536e-01} & \snum{8.769818e-02} & \fnum{2.578800e+01} (\fnum{8.687387e-01}) & \fnum{1.452717e+01} & \fnum{1.251754e+01} & \fnum{1.913174e+00} & \fnum{5.074596e-02} & $\ast$\fnum{2.968441e+01}\\
\runid &   5 &         & 127 & 616 & \snum{2.908183e-01} & \snum{3.599098e-02} & \fnum{1.755813e+01} (\fnum{8.634184e-01}) & \fnum{1.032160e+01} & \fnum{7.882262e+00} & \fnum{1.375166e+00} & \fnum{1.174543e-01} & \fnum{2.033559e+01}\\
\runid &  10 &         & 104 & 493 & \snum{2.901161e-01} & \snum{4.906508e-02} & \fnum{1.511931e+01} (\fnum{8.594104e-01}) & \fnum{9.259736e+00} & \fnum{6.543274e+00} & \fnum{9.313479e-01} & \fnum{1.680530e-01} & \fnum{1.759265e+01}\\
\runid &  15 &         & 97 & 455 & \snum{2.902453e-01} & \snum{4.913023e-02} & \fnum{1.367081e+01} (\fnum{8.461352e-01}) & \fnum{8.458035e+00} & \fnum{6.037219e+00} & \fnum{7.253370e-01} & \fnum{2.246642e-01} & \fnum{1.615677e+01}\\
\runid &  20 &         & 99 & 462 & \snum{2.901503e-01} & \snum{4.942363e-02} & \fnum{1.397419e+01} (\fnum{8.335741e-01}) & \fnum{8.665429e+00} & \fnum{6.145636e+00} & \fnum{7.247287e-01} & \fnum{3.079779e-01} & \fnum{1.676418e+01}\\
\runid &  25 &         & 85 & 396 & \snum{2.914708e-01} & \snum{4.986111e-02} & \fnum{1.233106e+01} (\fnum{8.321381e-01}) & \fnum{7.913533e+00} & \fnum{5.140258e+00} & \fnum{5.816720e-01} & \fnum{2.941341e-01} & \fnum{1.481853e+01}\\
\runid &  50 &         & 86 & 403 & \snum{2.915598e-01} & \snum{4.776222e-02} & \fnum{1.247659e+01} (\fnum{7.905661e-01}) & \fnum{8.006180e+00} & \fnum{5.334299e+00} & \fnum{6.157438e-01} & \fnum{5.748667e-01} & \fnum{1.578185e+01}\\
\midrule
\runid &   1 &   (6,3) & 141 & 833 & \snum{2.941284e-01} & \snum{4.932081e-02} & \fnum{2.302882e+01} (\fnum{8.731577e-01}) & \fnum{1.410126e+01} & \fnum{8.827923e+00} & \fnum{2.857697e+00} & \fnum{6.968179e-02} & \fnum{2.637418e+01}\\
\runid &   5 &         & 107 & 643 & \snum{2.896009e-01} & \snum{4.981086e-02} & \fnum{1.825456e+01} (\fnum{8.639795e-01}) & \fnum{1.134960e+01} & \fnum{6.622677e+00} & \fnum{2.380386e+00} & \fnum{1.872590e-01} & \fnum{2.112847e+01}\\
\runid &  10 &         & 84 & 501 & \snum{2.924237e-01} & \snum{4.996899e-02} & \fnum{1.469641e+01} (\fnum{8.518131e-01}) & \fnum{9.595176e+00} & \fnum{5.086872e+00} & \fnum{1.619742e+00} & \fnum{2.816332e-01} & \fnum{1.725309e+01}\\
\runid &  15 &         & 87 & 510 & \snum{2.903300e-01} & \snum{4.284844e-02} & \fnum{1.507477e+01} (\fnum{8.346408e-01}) & \fnum{9.957332e+00} & \fnum{5.257759e+00} & \fnum{1.514032e+00} & \fnum{4.261838e-01} & \fnum{1.806138e+01}\\
\runid &  20 &         & 76 & 451 & \snum{2.912253e-01} & \snum{4.724895e-02} & \fnum{1.388063e+01} (\fnum{8.277065e-01}) & \fnum{9.097158e+00} & \fnum{4.725855e+00} & \fnum{1.426130e+00} & \fnum{4.675160e-01} & \fnum{1.676999e+01}\\
\runid &  25 &         & 66 & 387 & \snum{2.924404e-01} & \snum{4.985464e-02} & \fnum{1.222995e+01} (\fnum{8.195486e-01}) & \fnum{8.350457e+00} & \fnum{4.056170e+00} & \fnum{1.107364e+00} & \fnum{4.564316e-01} & \fnum{1.492278e+01}\\
\runid &  50 &         & 68 & 400 & \snum{2.914196e-01} & \snum{4.577080e-02} & \fnum{1.239765e+01} (\fnum{7.704254e-01}) & \fnum{8.345680e+00} & \fnum{4.177363e+00} & \fnum{1.154962e+00} & \fnum{8.941053e-01} & \fnum{1.609195e+01}\\
\midrule
\runid &   1 &   (3,6) & 173 & 825 & \snum{2.914690e-01} & \snum{4.983384e-02} & \fnum{2.206721e+01} (\fnum{8.662569e-01}) & \fnum{1.276821e+01} & \fnum{1.042727e+01} & \fnum{1.634478e+00} & \fnum{4.337117e-02} & \fnum{2.547421e+01}\\
\runid &   5 &         & 137 & 669 & \snum{2.951814e-01} & \snum{4.997240e-02} & \fnum{1.902162e+01} (\fnum{8.646407e-01}) & \fnum{1.105566e+01} & \fnum{8.513722e+00} & \fnum{1.614191e+00} & \fnum{1.285211e-01} & \fnum{2.199945e+01}\\
\runid &  10 &         & 109 & 523 & \snum{2.897710e-01} & \snum{3.859592e-02} & \fnum{1.538839e+01} (\fnum{8.552580e-01}) & \fnum{9.325799e+00} & \fnum{6.749909e+00} & \fnum{1.028367e+00} & \fnum{1.869253e-01} & \fnum{1.799269e+01}\\
\runid &  15 &         & 100 & 473 & \snum{2.907927e-01} & \snum{4.705007e-02} & \fnum{1.448984e+01} (\fnum{8.491753e-01}) & \fnum{8.988653e+00} & \fnum{6.150976e+00} & \fnum{7.869815e-01} & \fnum{2.378585e-01} & \fnum{1.706342e+01}\\
\runid &  20 &         & 91 & 433 & \snum{2.910682e-01} & \snum{4.886429e-02} & \fnum{1.323420e+01} (\fnum{8.370767e-01}) & \fnum{8.382093e+00} & \fnum{5.488586e+00} & \fnum{7.703972e-01} & \fnum{2.769654e-01} & \fnum{1.581002e+01}\\
\runid &  25 &         & 91 & 431 & \snum{2.910791e-01} & \snum{4.572199e-02} & \fnum{1.308713e+01} (\fnum{8.289284e-01}) & \fnum{8.096480e+00} & \fnum{5.624195e+00} & \fnum{7.399641e-01} & \fnum{3.387260e-01} & \fnum{1.578801e+01}\\
\runid &  50 &         & 83 & 394 & \snum{2.917427e-01} & \snum{4.892793e-02} & \fnum{1.243553e+01} (\fnum{7.946064e-01}) & \fnum{8.034151e+00} & \fnum{5.086489e+00} & \fnum{6.521984e-01} & \fnum{5.672448e-01} & \fnum{1.564993e+01}\\
\bottomrule
\end{tabular}
\end{table}

\begin{table}
\caption{Continuation of the results reported in \Cref{t:nirep-ngmres-sg}.}\label{t:nirep-ngmres-sg_cont}
\tabadjust\footnotesize
\begin{tabular}{rrcrrrrRRRRRE}\toprule
\multicolumn{7}{c}{} & \multicolumn{6}{c}{\cellcolor{gray!20}time (in seconds)} \\
run & $w$ & $(\sigma,\tau)$ & \#iter & \#pdes & dist & grad & pdes & $q_k$ & $g$ & sg & ls & tts \\
\midrule
\runid &   1 &  (12,6) & 135 & 782 & \snum{2.917743e-01} & \snum{3.511956e-02} & \fnum{2.141292e+01} (\fnum{8.706375e-01}) & \fnum{1.349480e+01} & \fnum{8.285612e+00} & \fnum{2.286110e+00} & \fnum{6.776825e-02} & \fnum{2.459454e+01}\\
\runid &   5 &         & 119 & 719 & \snum{2.935120e-01} & \snum{4.988622e-02} & \fnum{1.936411e+01} (\fnum{8.569078e-01}) & \fnum{1.210266e+01} & \fnum{7.206849e+00} & \fnum{2.543115e+00} & \fnum{2.093611e-01} & \fnum{2.259766e+01}\\
\runid &  10 &         & 93 & 555 & \snum{2.929246e-01} & \snum{4.984678e-02} & \fnum{1.609474e+01} (\fnum{8.501012e-01}) & \fnum{1.031016e+01} & \fnum{5.638796e+00} & \fnum{1.829047e+00} & \fnum{2.989045e-01} & \fnum{1.893273e+01}\\
\runid &  15 &         & 91 & 533 & \snum{2.903590e-01} & \snum{4.564380e-02} & \fnum{1.555368e+01} (\fnum{8.338988e-01}) & \fnum{9.988913e+00} & \fnum{5.605059e+00} & \fnum{1.692577e+00} & \fnum{4.160300e-01} & \fnum{1.865176e+01}\\
\runid &  20 &         & 73 & 427 & \snum{2.921274e-01} & \snum{4.841045e-02} & \fnum{1.308862e+01} (\fnum{8.292143e-01}) & \fnum{8.673190e+00} & \fnum{4.506525e+00} & \fnum{1.267075e+00} & \fnum{4.272439e-01} & \fnum{1.578436e+01}\\
\runid &  25 &         & 66 & 395 & \snum{2.932656e-01} & \snum{4.666886e-02} & \fnum{1.233627e+01} (\fnum{8.150462e-01}) & \fnum{8.325544e+00} & \fnum{4.061610e+00} & \fnum{1.219586e+00} & \fnum{5.026520e-01} & \fnum{1.513566e+01}\\
\runid &  50 &         & 65 & 389 & \snum{2.926611e-01} & \snum{4.826998e-02} & \fnum{1.208229e+01} (\fnum{7.702250e-01}) & \fnum{8.216186e+00} & \fnum{3.957467e+00} & \fnum{1.166390e+00} & \fnum{8.353030e-01} & \fnum{1.568670e+01}\\
\midrule
\runid &   1 &  (6,12) & 163 & 761 & \snum{2.922743e-01} & \snum{4.634374e-02} & \fnum{2.089592e+01} (\fnum{8.699309e-01}) & \fnum{1.206084e+01} & \fnum{9.975364e+00} & \fnum{1.413374e+00} & \fnum{3.731267e-02} & \fnum{2.402021e+01}\\
\runid &   5 &         & 127 & 617 & \snum{2.922651e-01} & \snum{4.149246e-02} & \fnum{1.767973e+01} (\fnum{8.671326e-01}) & \fnum{1.029678e+01} & \fnum{7.843298e+00} & \fnum{1.544062e+00} & \fnum{9.709458e-02} & \fnum{2.038873e+01}\\
\runid &  10 &         & 112 & 544 & \snum{2.935963e-01} & \snum{4.728102e-02} & \fnum{1.617506e+01} (\fnum{8.564913e-01}) & \fnum{9.968047e+00} & \fnum{6.891259e+00} & \fnum{1.211685e+00} & \fnum{1.952680e-01} & \fnum{1.888526e+01}\\
\runid &  15 &         & 109 & 519 & \snum{2.900984e-01} & \snum{3.885953e-02} & \fnum{1.503133e+01} (\fnum{8.391296e-01}) & \fnum{9.039931e+00} & \fnum{6.709490e+00} & \fnum{9.976849e-01} & \fnum{2.633494e-01} & \fnum{1.791300e+01}\\
\runid &  20 &         & 109 & 518 & \snum{2.905500e-01} & \snum{4.174695e-02} & \fnum{1.533938e+01} (\fnum{8.298750e-01}) & \fnum{9.277643e+00} & \fnum{6.594216e+00} & \fnum{9.818230e-01} & \fnum{3.502457e-01} & \fnum{1.848396e+01}\\
\runid &  25 &         & 94 & 452 & \snum{2.910408e-01} & \snum{4.847698e-02} & \fnum{1.371167e+01} (\fnum{8.254043e-01}) & \fnum{8.463884e+00} & \fnum{5.777474e+00} & \fnum{8.679311e-01} & \fnum{3.850825e-01} & \fnum{1.661206e+01}\\
\runid &  50 &         & 91 & 435 & \snum{2.914698e-01} & \snum{4.681612e-02} & \fnum{1.341556e+01} (\fnum{7.838685e-01}) & \fnum{8.333890e+00} & \fnum{5.655304e+00} & \fnum{8.113343e-01} & \fnum{6.702174e-01} & \fnum{1.711455e+01}\\
\midrule
\runid &   1 &   (1,1) & 200 & 1060 & \snum{3.001037e-01} & \snum{9.056446e-02} & \fnum{2.781685e+01} (\fnum{8.668855e-01}) & \fnum{1.585097e+01} & \fnum{1.224698e+01} & \fnum{3.214341e+00} & \fnum{7.448971e-02} & $\ast$\fnum{3.208826e+01}\\
\runid &  10 &   (7,4) & 89 & 512 & \snum{2.912357e-01} & \snum{4.973795e-02} & \fnum{1.503809e+01} (\fnum{8.500563e-01}) & \fnum{9.723009e+00} & \fnum{5.417252e+00} & \fnum{1.562600e+00} & \fnum{2.676246e-01} & \fnum{1.769070e+01}\\
\runid &  15 &  (10,6) & 83 & 474 & \snum{2.907449e-01} & \snum{4.916774e-02} & \fnum{1.431163e+01} (\fnum{8.386511e-01}) & \fnum{9.328480e+00} & \fnum{5.133292e+00} & \fnum{1.376403e+00} & \fnum{3.514964e-01} & \fnum{1.706506e+01}\\
\runid &  20 &  (13,8) & 76 & 446 & \snum{2.914697e-01} & \snum{4.949313e-02} & \fnum{1.343860e+01} (\fnum{8.233997e-01}) & \fnum{8.927230e+00} & \fnum{4.596952e+00} & \fnum{1.327582e+00} & \fnum{4.603358e-01} & \fnum{1.632087e+01}\\
\runid &  25 & (16,10) & 79 & 446 & \snum{2.915494e-01} & \snum{4.249136e-02} & \fnum{1.359339e+01} (\fnum{8.127147e-01}) & \fnum{8.912355e+00} & \fnum{4.869675e+00} & \fnum{1.266258e+00} & \fnum{5.011627e-01} & \fnum{1.672590e+01}\\
\runid &  50 & (39,12) & 65 & 409 & \snum{2.915501e-01} & \snum{4.510878e-02} & \fnum{1.258022e+01} (\fnum{7.714900e-01}) & \fnum{8.587219e+00} & \fnum{3.980288e+00} & \fnum{1.347329e+00} & \fnum{8.526292e-01} & \fnum{1.630640e+01}\\
\midrule
\runid &  10 &  (1,10) & 177 & 698 & \snum{2.940200e-01} & \snum{3.997183e-02} & \fnum{1.953486e+01} (\fnum{8.533721e-01}) & \fnum{1.060320e+01} & \fnum{1.089652e+01} & \fnum{3.718300e-01} & \fnum{8.621463e-02} & \fnum{2.289138e+01}\\
\runid &  15 &  (1,15) & 145 & 565 & \snum{2.911640e-01} & \snum{3.562045e-02} & \fnum{1.635719e+01} (\fnum{8.513603e-01}) & \fnum{9.092041e+00} & \fnum{8.958231e+00} & \fnum{2.196725e-01} & \fnum{7.398608e-02} & \fnum{1.921300e+01}\\
\runid &  20 &  (1,20) & 148 & 567 & \snum{2.927380e-01} & \snum{3.765537e-02} & \fnum{1.618749e+01} (\fnum{8.420407e-01}) & \fnum{9.001590e+00} & \fnum{9.080683e+00} & \fnum{1.747266e-01} & \fnum{8.001258e-02} & \fnum{1.922412e+01}\\
\runid &  25 &  (1,25) & 157 & 598 & \snum{2.921966e-01} & \snum{3.323215e-02} & \fnum{1.708533e+01} (\fnum{8.337856e-01}) & \fnum{9.143919e+00} & \fnum{9.749967e+00} & \fnum{1.545350e-01} & \fnum{8.816238e-02} & \fnum{2.049127e+01}\\
\runid &  50 &  (1,50) & 154 & 579 & \snum{2.930211e-01} & \snum{4.773140e-02} & \fnum{1.655645e+01} (\fnum{8.117332e-01}) & \fnum{8.926855e+00} & \fnum{9.482999e+00} & \fnum{8.738121e-02} & \fnum{1.131648e-01} & \fnum{2.039641e+01}\\
\midrule
\runid & 400 &   (1,0) & 57 & 392 & \snum{2.913499e-01} & \snum{4.767728e-02} & \fnum{1.262487e+01} (\fnum{7.799062e-01}) & \fnum{8.862761e+00} & \fnum{3.511291e+00} & \fnum{1.441079e+00} & \fnum{9.175515e-01} & \fnum{1.618768e+01}\\
\runid &     &   (1,1) & 73 & 383 & \snum{2.909504e-01} & \snum{4.170841e-02} & \fnum{1.197353e+01} (\fnum{7.725759e-01}) & \fnum{7.880816e+00} & \fnum{4.492121e+00} & \fnum{8.487408e-01} & \fnum{7.479200e-01} & \fnum{1.549819e+01}\\
\runid &     &   (2,2) & 73 & 385 & \snum{2.912664e-01} & \snum{4.697974e-02} & \fnum{1.208544e+01} (\fnum{7.680170e-01}) & \fnum{7.948949e+00} & \fnum{4.498433e+00} & \fnum{8.538021e-01} & \fnum{8.432205e-01} & \fnum{1.573590e+01}\\
\runid &     &   (5,5) & 81 & 427 & \snum{2.910007e-01} & \snum{4.100724e-02} & \fnum{1.294559e+01} (\fnum{7.502164e-01}) & \fnum{8.378520e+00} & \fnum{4.989498e+00} & \fnum{9.581141e-01} & \fnum{1.061301e+00} & \fnum{1.725581e+01}\\
\runid &     &   (8,8) & 82 & 434 & \snum{2.915344e-01} & \snum{4.629416e-02} & \fnum{1.335930e+01} (\fnum{7.552451e-01}) & \fnum{8.584770e+00} & \fnum{5.074966e+00} & \fnum{1.059866e+00} & \fnum{1.063199e+00} & \fnum{1.768870e+01}\\
\bottomrule
\end{tabular}
\end{table}

\begin{table}
\caption{Comparison against \gaaa. For each dataset, with and without the safeguard of \Cref{a:gangmres-sg}, we report how many of the $86$ runs of the hyperparameter sweep reach the tolerance $\epsilon_{\scr{rel}} = 5 \times 10^{-2}$ within $n_{\scr{iter}} = 200$ iterations, and the median number of (outer) iterations and of PDE solves over all runs of that sweep. Runs that do not reach the tolerance exhaust the iteration budget, so the medians for \gaaa\ reflect the cost of $200$ iterations, not of convergence. We did not run \gaaa\ without the safeguard on \nirep. The runs of \gangmres\ with the safeguard on \nirep\ are in \Cref{t:nirep-ngmres-sg,t:nirep-ngmres-sg_cont}; for brevity, the other individual runs are not included.}
\label{t:aa-summary}
\tabadjust
\begin{tabular}{llcrrr}\toprule
scheme & data & safeguard & converged & median \#iter & median \#pdes \\
\midrule
\gangmres\ & \hands\ & no  & $86/86$ &  $38$ &  $186$ \\
\gangmres\ & \hands\ & yes & $86/86$ &  $38$ &  $218$ \\
\gaaa\     & \hands\ & no  & $ 3/86$ & $200$ &  $741$ \\
\gaaa\     & \hands\ & yes & $10/86$ & $200$ &  $929$ \\
\midrule
\gangmres\ & \nirep\ & no  & $81/86$ & $100$ &  $471$ \\
\gangmres\ & \nirep\ & yes & $82/86$ &  $93$ &  $471$ \\
\gaaa\     & \nirep\ & yes & $ 0/86$ & $200$ &  $918$ \\
\bottomrule
\end{tabular}
\end{table}

\subsubsection{Comparison against Anderson Acceleration}\label{s:results:aa}

We compare \gangmres\ against \gaaa, the generalized alternating variant of \andacc, on the hyperparameter sweeps of \Cref{s:results:conv}; \Cref{t:aa-summary} summarizes the outcome. On \hands, \gangmres\ reaches the tolerance for all 86 hyperparameter settings with and without the safeguard, whereas \gaaa\ reaches it for three without and ten with the safeguard. On \nirep, \gaaa\ with the safeguard reaches it for none, at about twice the median number of PDE solves of \gangmres. The extra unknown in the least-squares problem of \gangmres\ at equal $w$ does not explain this, since \gaaa\ also fails with the untruncated window $w=400$. Both schemes solve the least-squares problem by a QR factorization. We also ran the \hands\ sweeps of \gaaa\ with the normal equations of the least-squares problem in place of QR, with $10^{-6}$ times the identity added to the Gram matrix for stability. Without the safeguard, the outcome does not change: three of the $86$ hyperparameter settings reach the tolerance either way, at a median of $743$ PDE solves against $741$ with QR. With the safeguard, the normal equations reach the tolerance for one setting instead of ten, at about the same median cost ($940$ against $929$ PDE solves).

Per iteration, \gaaa\ is cheaper, since its residual $r_k(v^{(k)})$ is free, whereas \ngmres\ minimizes $g(q_k(v^{(k)}))$ at the cost of one forward and one adjoint solve. Without the safeguard, \gaaa\ needs $3.7$ PDE solves per iteration against $4.9$ (ratio of the median counts), but $5.3$ times as many iterations and $4.0$ times as many PDE solves; it never needs fewer PDE solves than \gangmres. The safeguard reduces this per-iteration advantage, since its descent test requires the reduced gradient at $q_k(v^{(k)})$. It removes the erratic behavior of \gaaa\ (the relative norm of the gradient increases in a median of 43\% of the iterations without it and in 32\% with it, and the largest excursion above the running minimum drops from a factor of $11.6$ to $1.5$), but not its slow convergence. For brevity, we do not include the individual runs.

\subsubsection{Verification of the Convergence Theory}\label{s:results:theory}

\Cref{s:appendix:theory} checks the analysis of \Cref{s:convergence} numerically. For the Gauss--Newton approximation of the reduced space Hessian, the contraction condition fails for five of the six regularization parameters, which indicates that the \fp\ map is not a contraction and that the classical theory does not explain the observed convergence. On \ngmres\ steps, $\theta^{(k+1)}$ and $\gamma^{(k+1)}$ never exceed one, and the factors in~\cref{e:theta-fp} behave as predicted for the alternating scheme, but $\theta^{(k+1)}$ does not predict the observed rate. The window size acts through two opposing effects, a smaller angle in~\cref{e:gamma-angle} and larger coefficients $\sum_j|\zeta_j|$, so that a larger window does not necessarily reduce the iteration count.

\subsubsection{Comparison against the Quasi-Newton Method in the Same Metric}\label{s:results:qn}

For the \lbfgs\ scheme the only parameter we sweep is the memory $n_{\scr{mem}}$ (the number of stored secant pairs).
At the default tolerance $\epsilon_{\scr{rel}} = 5\times 10^{-2}$, \lbfgs\ is the fastest scheme on \hands\ ($1.04$ against $1.29$\,s for the best \nk\ variant). On \nirep\ at $\alpha = 10^{-3}$, \nk\ (\texttt{h0rpc}) is about two percent faster ($9.39$\,s against $9.61$\,s for \lbfgs\ with the best memory in \Cref{t:h2-nk+rpgd}, $n_{\scr{mem}} = 50$); with $n_{\scr{mem}} = 100$, which we include only in the regularization study in \Cref{t:nirep-nk-rpdg-regul-new}, \lbfgs\ needs $9.13$\,s. Comparing the best configuration of each, \lbfgs\ is faster than \gangmres\ by a factor of $1.63$ on \hands\ and $1.52$ on \nirep, and by $1.49$ to $1.80$ across the regularization parameters of \Cref{s:results:regularization}. It is also faster than the best of the \rpgd\ and \nk\ variants at five of these six regularization parameters; at the smallest, $\alpha = 5\times 10^{-4}$, \nk\ (\texttt{h0rpc}) is faster than every history size for \lbfgs\ we considered and supplies the reference of the speedup column of \Cref{t:nirep-ngmres-regul-new}. In the mesh study of \Cref{s:results:mesh}, which uses the tighter tolerance $\epsilon_{\scr{rel}} = 10^{-3}$, \lbfgs\ requires fewer PDE solves than \gangmres\ and is faster than both \gangmres\ and \nk\ in the six runs in which it reaches the tolerance; for $n_i = 64$ it does not. \gangmres\ is therefore not the fastest scheme for this problem. At tighter tolerances, it overtakes \nk\ (\texttt{h0rpc}) on \nirep, but not \lbfgs.

\subsubsection{Dependence on the Stopping Tolerance}\label{s:results:tolerance}

The results for the real-world data use $\epsilon_{\scr{rel}} = 5\times 10^{-2}$; those for the synthetic data in \Cref{s:results:mesh,s:stokes-flow,s:optimal-transport} use $10^{-3}$. The former tolerance is loose, and the ranking of the schemes depends on it. We therefore repeat the entire \nirep\ study at $10^{-2}$ and $5\times 10^{-3}$, including the four baselines, the four \lbfgs\ memories, and the $86$ settings of \gangmres\ in \Cref{t:nirep-ngmres-sg,t:nirep-ngmres-sg_cont}, each with the protocol of \Cref{s:results:conv}. We report the time to solution of the fastest configuration of each family that reaches the tolerance and how many configurations of each family do so in \Cref{t:tolerance}.

\begin{table}
\caption{Dependence on the stopping tolerance $\epsilon_{\scr{rel}}$ (\nirep; $\alpha = 10^{-3}$). Column conv.\ gives how many configurations of each method reach the tolerance: of the three preconditioners of \nk\ (\texttt{ireg}, \texttt{2lrpcsym}, \texttt{h0rpc}), of the four memory sizes of \lbfgs\ ($n_{\scr{mem}} \in \{5,10,20,50\}$), and of the $86$ hyperparameter settings of \gangmres. Column tts gives the time to solution in seconds of the fastest of these configurations and, in parentheses, that configuration: the preconditioner of \nk, the memory size $n_{\scr{mem}}$ of \lbfgs, and $w$; $\sigma,\tau$ of \gangmres. ---: no configuration reaches the tolerance within the iteration budget. \rpgd\ reaches none of the three tolerances within the iteration budget. Runtimes are the median of five repetitions.}\label{t:tolerance}
\centering
\tabadjust
\begin{tabular}{lrrrrrr}\toprule
& \multicolumn{2}{c}{\nk} & \multicolumn{2}{c}{\lbfgs} & \multicolumn{2}{c}{\gangmres} \\
\cmidrule(lr){2-3}\cmidrule(lr){4-5}\cmidrule(lr){6-7}
$\epsilon_{\scr{rel}}$ & tts & conv. & tts & conv. & tts & conv. \\\midrule
$5\times 10^{-2}$ & $ 9.39$ (\texttt{h0rpc}) & $3/3$ & $ 9.61$ ($50$) & $4/4$ & $14.56$ ($25$; $4,2$) & $82/86$ \\
$10^{-2}$         & $21.43$ (\texttt{h0rpc}) & $3/3$ & $15.76$ ($20$) & $3/4$ & $18.83$ ($50$; $4,2$) & $65/86$ \\
$5\times 10^{-3}$ & $31.56$ (\texttt{h0rpc}) & $3/3$ & $25.01$ ($50$) & $1/4$ & ---     & $ 0/86$ \\
\bottomrule
\end{tabular}
\end{table}

We draw three conclusions, only the first of which favors our scheme. First, tightening the tolerance increases the cost of \gangmres\ the least. Going from $5\times 10^{-2}$ to $10^{-2}$, the time to solution of the fastest configuration grows by a factor of $1.3$ for \gangmres, $1.6$ for \lbfgs, and $2.3$ for \nk. As a result, \gangmres\ is the slowest of the three at $5\times 10^{-2}$ but faster than \nk\ at $10^{-2}$; it does not overtake \lbfgs\ at any of the three tolerances. Second, \gangmres\ and \lbfgs\ are comparably sensitive to their hyperparameters: at $10^{-2}$, each reaches the tolerance for about three quarters of its configurations ($65$ of $86$ and three of four, respectively). Third, \gangmres\ does not reach tight tolerances on these data within the iteration budget. At $5\times 10^{-3}$, none of its $86$ configurations reaches the tolerance within $n_{\scr{iter}} = 200$ iterations, whereas \nk\ and \lbfgs\ still do. These runs do not stall; they are still descending when they hit the iteration limit (best: $6.68\times 10^{-3}$, median: $9.35\times 10^{-3}$). Below roughly $10^{-2}$, \nk\ is therefore the only one of the three schemes that converges reliably on \nirep. This limitation is specific to the non-smooth real-world data. On the smooth synthetic data of \Cref{s:results:mesh,s:stokes-flow,s:optimal-transport}, \gangmres\ reaches $\epsilon_{\scr{rel}} = 10^{-3}$ in all runs of the mesh convergence and Stokes studies and in $19$ of the $20$ runs of the mass-preserving transport study; the ranking does not reverse there, however, since \lbfgs\ remains faster in the mesh convergence study.

\subsubsection{Summary of Observations}\label{s:results:conv-observations}

The most important observations are:
\begin{enumerate*}[label=\it(\roman*)]
\item \gangmres\ converges where \rpgd\ does not. \rpgd\ does not reach the tolerance within $200$ iterations on either dataset, whereas \gangmres\ reaches it for all $86$ hyperparameter settings on \hands\ and for $82$ of them on \nirep.
\item \gangmres\ is not the fastest method. \nk\ with the zero-velocity preconditioner (\texttt{h0rpc}) and \lbfgs\ are faster than every configuration of \gangmres\ on both datasets. On \hands, the fastest runs take $1.29$\,s for \nk\ (\texttt{h0rpc}), $1.04$\,s for \lbfgs, and $1.69$\,s for \gangmres; on \nirep, they take $9.39$\,s, $9.61$\,s, and $14.56$\,s, respectively. \gangmres\ is faster than \nk\ with the spectral preconditioner (\texttt{ireg}) on both datasets. It is faster than \nk\ with the two-level preconditioner (\texttt{2lrpcsym}) on \nirep\ ($14.56$ against $21.75$\,s), but not on \hands\ ($1.69$ against $1.37$\,s).
\item \gangmres\ requires fewer PDE solves than \nk\ with \texttt{ireg} and \texttt{2lrpcsym}, if we count the PDE solves inside the preconditioner, but more than \nk\ with \texttt{h0rpc}. On \hands, \gangmres\ requires a median of $218$ PDE solves over the $86$ hyperparameter settings (for brevity, we do not include the table of these runs). \nk\ requires $246$ PDE solves with \texttt{ireg}, $306$ with \texttt{2lrpcsym} ($56$ in the solver and $250$ in the preconditioner), and $107$ with \texttt{h0rpc}. On \nirep, the median for \gangmres\ is $471$, against $1652$, $1674$, and $216$ for the three preconditioners. The advantage of \texttt{h0rpc} is specific to the image data. In the mesh convergence study, \gangmres\ requires fewer PDE solves than \nk\ with \texttt{h0rpc} at six of the seven mesh points (see \Cref{s:results:mesh}).
\item For $w \in \{10,15,20,25\}$, \gangmres\ is robust with respect to the sequence $p = (\sigma,\tau)$. All of these runs reach the tolerance, and their time to solution varies by a factor of $1.6$ on \nirep\ (from $14.56$ to $22.89$\,s) and of $1.7$ on \hands\ (from $1.69$ to $2.88$\,s). The four runs that do not reach the tolerance on \nirep\ all have $w \le 5$.
\item Larger windows do not necessarily help. We run ten sequences $p$ with both $w = 25$ and $w = 50$. On \nirep, $w = 50$ needs fewer iterations for six of them but is faster for only two; on \hands, it needs fewer iterations for none and is faster for three. The fastest run uses $w = 25$, $p = (4,2)$ on \nirep\ and $w = 15$, $p = (5,5)$ on \hands.
\item The least-squares problem is cheap compared to the PDE solves. With the safeguard, which restarts the window whenever it rejects the accelerated iterate, this holds even for the untruncated window: for $w = 400$, $p = (1,0)$ on \nirep, the least-squares problem takes $0.92$\,s and the PDE solves take $12.62$\,s. Without the safeguard, the two become comparable ($26.6$ against $30.0$\,s).
\item Sequences with $\sigma \ge \tau$ need fewer iterations than those with $\sigma < \tau$ (medians of $34$ against $44$ on \hands\ and of $81$ against $109$ on \nirep), and three of the four runs that do not reach the tolerance on \nirep\ have $\sigma < \tau$. The time to solution, however, differs little (medians of $17.21$ against $17.91$\,s on \nirep\ and of $1.93$ against $1.92$\,s on \hands).
\item The \gangmres\ scheme significantly outperforms the \gaaa\ scheme. With the safeguard, \gaaa\ reaches the tolerance for ten of the $86$ hyperparameter settings on \hands\ and for none on \nirep\ (see \Cref{s:results:aa}).
\end{enumerate*}

In summary, using first-order information only, \gangmres\ attains runtimes within the range spanned by the \nk\ variants we consider. At this tolerance, \lbfgs\ is the fastest method on \hands; on \nirep, it is about two percent slower than \nk\ with \texttt{h0rpc}.

\subsection{Regularization Parameter Sensitivity}\label{s:results:regularization}

We study computational performance for vanishing regularization parameters $\alpha \to 0$. We expect the performance to deteriorate as the regularization parameter becomes smaller (the problem becomes more ill-conditioned).

We consider \nk, \rpgd, \lbfgs\ (with memories $\{20, 50, 100\}$), and \gangmres. We test performance for the \nirep\ dataset (native resolution: $300 \times 300$). We select $\alpha$ in \[\{10^{-1}, 5 \times 10^{-2}, 10^{-2}, 5 \times 10^{-3}, 10^{-3}, 5 \times 10^{-4}\}.\] Based on the prior experiments, we limit $p = (\sigma,\tau)$ to $(5,1)$ and $(4,2)$ and select $w$ in $\{10,15,20,25\}$.

To quantify the performance of the considered numerical schemes, we report the number of iterations, the number of PDE solves, the number of Hessian matvecs, the relative change of the data mismatch, and the relative change of the $\ell^\infty$-norm of the reduced gradient. We also report various execution times (in seconds), including the time spent on solving PDEs, Hessian matrix products, the evaluation of $q_k$, the evaluation of the reduced gradient $g$, the time spent on solving the least squares system that appears in the \gangmres\ schemes, and the time-to-solution (total runtime).

The results for the baseline methods (\rpgd, \lbfgs, and \nk) are reported in \Cref{t:nirep-nk-rpdg-regul-new}. The results for \gangmres\ are reported in \Cref{t:nirep-ngmres-regul-new}. The speedup reported in \Cref{t:nirep-ngmres-regul-new} is based on the best performing baseline in \Cref{t:nirep-nk-rpdg-regul-new} for each choice of $\alpha$ (highlighted in red), i.e., the quickest \rpgd, \lbfgs\ or \nk\ run.

\begin{table}[!t]
\caption{Convergence results for \rpgd, \lbfgs\ and \nk\ for the \nirep\ data as a function of the regularization parameter $\alpha$. Setup, reported quantities and layout are as in \Cref{t:h2-nk+rpgd_hands}. Runtimes with $\ast$ indicate that the scheme did not reach the tolerance $\epsilon_{\scr{rel}}$ within the iteration budget. The highlighted row, the quickest run for each $\alpha$, is the reference for the speedup in \Cref{t:nirep-ngmres-regul-new}. Runtimes are the median of five repetitions.}
\label{t:nirep-nk-rpdg-regul-new}
\tabadjust\setcounter{run}{0}
\addtolength{\tabcolsep}{-0.12em}\begin{tabular}{FlrrrrrrRRRE}\toprule
\multicolumn{8}{c}{} & \multicolumn{4}{c}{\cellcolor{gray!20}time (in seconds)} \\
run & method & \#iter & \#pdes & \#pdes$_{\scr{pc}}$ & \#mvs & dist & grad & pdes & pdes$_{\scr{pc}}$ & mvs & tts \\
\midrule
\multicolumn{12}{l}{$\alpha = 10^{-1}$} \\
\midrule
\runid & \rpgd                    &  21 &   76 &  --- &  --- & \snum{4.655556e-01} & \snum{4.453441e-02} & \fnum{1.846811e+00} (\fnum{8.617794e-01}) & ---                                 & ---                                 & \fnum{2.143021e+00}\\
\runid & \nk\ (\texttt{ireg})     &  11 &  180 &  --- &   74 & \snum{4.726577e-01} & \snum{4.411519e-02} & \fnum{5.554159e+00} (\fnum{8.829455e-01}) & ---                                 & \fnum{5.271247e+00} (\fnum{8.379709e-01}) & \fnum{6.290489e+00}\\
\runid & \nk\ (\texttt{2lrpcsym}) &  12 &  101 &  168 &   33 & \snum{4.705102e-01} & \snum{4.002755e-02} & \fnum{3.040168e+00} (\fnum{5.786428e-01}) & \fnum{1.034326e+00} (\fnum{1.968659e-01}) & \fnum{2.477609e+00} (\fnum{4.715695e-01}) & \fnum{5.253963e+00}\\
\runid & \nk\ (\texttt{h0rpc})    &  14 &  177 &  --- &   68 & \snum{4.682336e-01} & \snum{3.478840e-02} & \fnum{5.403754e+00} (\fnum{7.627729e-01}) & ---                                 & \fnum{4.864872e+00} (\fnum{6.867064e-01}) & \fnum{7.084355e+00}\\
\rowcolor{p1color!20}
\runid & L-BFGS (20)              &  12 &   38 &  --- &  --- & \snum{4.645337e-01} & \snum{4.555599e-02} & \fnum{9.070968e-01} (\fnum{8.227931e-01}) & ---                                 & ---                                 & \fnum{1.102460e+00}\\
\runid & L-BFGS (50)              &  12 &   38 &  --- &  --- & \snum{4.645337e-01} & \snum{4.555599e-02} & \fnum{9.168660e-01} (\fnum{8.234903e-01}) & ---                                 & ---                                 & \fnum{1.113390e+00}\\
\runid & L-BFGS (100)             &  12 &   38 &  --- &  --- & \snum{4.645337e-01} & \snum{4.555599e-02} & \fnum{9.219381e-01} (\fnum{8.255621e-01}) & ---                                 & ---                                 & \fnum{1.116740e+00}\\
\midrule
\multicolumn{12}{l}{$\alpha = 5 \times 10^{-2}$} \\
\midrule
\runid & \rpgd                    &  32 &  118 &  --- &  --- & \snum{4.309854e-01} & \snum{4.681355e-02} & \fnum{2.868673e+00} (\fnum{8.649781e-01}) & ---                                 & ---                                 & \fnum{3.316469e+00}\\
\runid & \nk\ (\texttt{ireg})     &  12 &  241 &  --- &  103 & \snum{4.337903e-01} & \snum{3.892403e-02} & \fnum{7.452113e+00} (\fnum{8.819872e-01}) & ---                                 & \fnum{7.305092e+00} (\fnum{8.645867e-01}) & \fnum{8.449230e+00}\\
\runid & \nk\ (\texttt{2lrpcsym}) &  13 &  110 &  216 &   36 & \snum{4.320315e-01} & \snum{4.532482e-02} & \fnum{3.358583e+00} (\fnum{5.478778e-01}) & \fnum{1.348039e+00} (\fnum{2.199023e-01}) & \fnum{2.744461e+00} (\fnum{4.476974e-01}) & \fnum{6.130169e+00}\\
\runid & \nk\ (\texttt{h0rpc})    &  15 &  190 &  --- &   73 & \snum{4.308888e-01} & \snum{4.183959e-02} & \fnum{5.699517e+00} (\fnum{7.522500e-01}) & ---                                 & \fnum{5.137704e+00} (\fnum{6.780992e-01}) & \fnum{7.576626e+00}\\
\runid & L-BFGS (20)              &  14 &   45 &  --- &  --- & \snum{4.292795e-01} & \snum{4.380370e-02} & \fnum{1.104732e+00} (\fnum{8.225642e-01}) & ---                                 & ---                                 & \fnum{1.343034e+00}\\
\runid & L-BFGS (50)              &  14 &   45 &  --- &  --- & \snum{4.292795e-01} & \snum{4.380370e-02} & \fnum{1.112753e+00} (\fnum{8.248690e-01}) & ---                                 & ---                                 & \fnum{1.349005e+00}\\
\rowcolor{p1color!20}
\runid & L-BFGS (100)             &  14 &   45 &  --- &  --- & \snum{4.292795e-01} & \snum{4.380370e-02} & \fnum{1.097865e+00} (\fnum{8.206306e-01}) & ---                                 & ---                                 & \fnum{1.337831e+00}\\
\midrule
\multicolumn{12}{l}{$\alpha = 10^{-2}$} \\
\midrule
\runid & \rpgd                    & 200 &  722 &  --- &  --- & \snum{3.794229e-01} & \snum{5.083453e-02} & \fnum{1.784049e+01} (\fnum{8.636206e-01}) & ---                                 & ---                                 & $\ast$\fnum{2.065779e+01}\\
\runid & \nk\ (\texttt{ireg})     &  14 &  451 &  --- &  205 & \snum{3.799046e-01} & \snum{3.605962e-02} & \fnum{1.399018e+01} (\fnum{8.787256e-01}) & ---                                 & \fnum{1.442837e+01} (\fnum{9.062484e-01}) & \fnum{1.592098e+01}\\
\runid & \nk\ (\texttt{2lrpcsym}) &  15 &  128 &  424 &   42 & \snum{3.786731e-01} & \snum{3.664973e-02} & \fnum{3.928993e+00} (\fnum{4.365596e-01}) & \fnum{2.685261e+00} (\fnum{2.983656e-01}) & \fnum{3.218356e+00} (\fnum{3.575991e-01}) & \fnum{8.999900e+00}\\
\runid & \nk\ (\texttt{h0rpc})    &  14 &  181 &  --- &   70 & \snum{3.785235e-01} & \snum{3.004500e-02} & \fnum{5.468660e+00} (\fnum{7.412067e-01}) & ---                                 & \fnum{4.956791e+00} (\fnum{6.718293e-01}) & \fnum{7.378051e+00}\\
\runid & L-BFGS (20)              &  25 &   80 &  --- &  --- & \snum{3.773673e-01} & \snum{3.016614e-02} & \fnum{2.053133e+00} (\fnum{8.057488e-01}) & ---                                 & ---                                 & \fnum{2.548106e+00}\\
\runid & L-BFGS (50)              &  25 &   80 &  --- &  --- & \snum{3.773957e-01} & \snum{3.052797e-02} & \fnum{2.054819e+00} (\fnum{8.042196e-01}) & ---                                 & ---                                 & \fnum{2.555048e+00}\\
\rowcolor{p1color!20}
\runid & L-BFGS (100)             &  25 &   80 &  --- &  --- & \snum{3.773957e-01} & \snum{3.052797e-02} & \fnum{2.058961e+00} (\fnum{8.101282e-01}) & ---                                 & ---                                 & \fnum{2.541525e+00}\\
\midrule
\multicolumn{12}{l}{$\alpha = 5 \times 10^{-3}$} \\
\midrule
\runid & \rpgd                    & 200 &  725 &  --- &  --- & \snum{3.650153e-01} & \snum{1.549686e-01} & \fnum{1.856661e+01} (\fnum{8.683243e-01}) & ---                                 & ---                                 & $\ast$\fnum{2.138211e+01}\\
\runid & \nk\ (\texttt{ireg})     &  15 &  624 &  --- &  290 & \snum{3.601467e-01} & \snum{4.467131e-02} & \fnum{2.011880e+01} (\fnum{8.780755e-01}) & ---                                 & \fnum{2.111082e+01} (\fnum{9.213715e-01}) & \fnum{2.291238e+01}\\
\runid & \nk\ (\texttt{2lrpcsym}) &  16 &  137 &  578 &   45 & \snum{3.577160e-01} & \snum{3.896343e-02} & \fnum{4.236312e+00} (\fnum{3.854732e-01}) & \fnum{3.709192e+00} (\fnum{3.375092e-01}) & \fnum{3.473419e+00} (\fnum{3.160556e-01}) & \fnum{1.098990e+01}\\
\runid & \nk\ (\texttt{h0rpc})    &  15 &  202 &  --- &   79 & \snum{3.537369e-01} & \snum{4.480370e-02} & \fnum{6.163180e+00} (\fnum{7.254539e-01}) & ---                                 & \fnum{5.633804e+00} (\fnum{6.631423e-01}) & \fnum{8.495619e+00}\\
\runid & L-BFGS (20)              &  36 &  114 &  --- &  --- & \snum{3.550659e-01} & \snum{4.934590e-02} & \fnum{3.056734e+00} (\fnum{7.975369e-01}) & ---                                 & ---                                 & \fnum{3.832718e+00}\\
\runid & L-BFGS (50)              &  33 &  105 &  --- &  --- & \snum{3.550562e-01} & \snum{4.721313e-02} & \fnum{2.806060e+00} (\fnum{7.994864e-01}) & ---                                 & ---                                 & \fnum{3.509829e+00}\\
\rowcolor{p1color!20}
\runid & L-BFGS (100)             &  33 &  105 &  --- &  --- & \snum{3.550562e-01} & \snum{4.721313e-02} & \fnum{2.816377e+00} (\fnum{8.033216e-01}) & ---                                 & ---                                 & \fnum{3.505915e+00}\\
\midrule
\multicolumn{12}{l}{$\alpha = 10^{-3}$} \\
\midrule
\runid & \rpgd                    & 200 &  711 &  --- &  --- & \snum{3.517696e-01} & \snum{3.738752e-01} & \fnum{2.075934e+01} (\fnum{8.811661e-01}) & ---                                 & ---                                 & $\ast$\fnum{2.355894e+01}\\
\runid & \nk\ (\texttt{ireg})     &  19 & 1652 &  --- &  798 & \snum{2.886617e-01} & \snum{4.244886e-02} & \fnum{5.518747e+01} (\fnum{8.827094e-01}) & ---                                 & \fnum{5.936521e+01} (\fnum{9.495312e-01}) & \fnum{6.252054e+01}\\
\runid & \nk\ (\texttt{2lrpcsym}) &  19 &  164 & 1510 &   54 & \snum{2.883969e-01} & \snum{4.638026e-02} & \fnum{5.050468e+00} (\fnum{2.321933e-01}) & \fnum{9.690895e+00} (\fnum{4.455352e-01}) & \fnum{4.162867e+00} (\fnum{1.913862e-01}) & \fnum{2.175113e+01}\\
\runid & \nk\ (\texttt{h0rpc})    &  15 &  216 &  --- &   86 & \snum{2.922136e-01} & \snum{4.909737e-02} & \fnum{6.444696e+00} (\fnum{6.866354e-01}) & ---                                 & \fnum{6.002716e+00} (\fnum{6.395456e-01}) & \fnum{9.385908e+00}\\
\runid & L-BFGS (20)              &  92 &  291 &  --- &  --- & \snum{2.899805e-01} & \snum{4.860600e-02} & \fnum{1.035247e+01} (\fnum{8.246365e-01}) & ---                                 & ---                                 & \fnum{1.255398e+01}\\
\runid & L-BFGS (50)              &  62 &  196 &  --- &  --- & \snum{2.903500e-01} & \snum{4.625843e-02} & \fnum{7.942120e+00} (\fnum{8.266105e-01}) & ---                                 & ---                                 & \fnum{9.608055e+00}\\
\rowcolor{p1color!20}
\runid & L-BFGS (100)             &  59 &  187 &  --- &  --- & \snum{2.904002e-01} & \snum{4.881566e-02} & \fnum{7.628853e+00} (\fnum{8.352304e-01}) & ---                                 & ---                                 & \fnum{9.133831e+00}\\
\midrule
\multicolumn{12}{l}{$\alpha = 5 \times 10^{-4}$} \\
\midrule
\runid & \rpgd                    & 200 &  713 &  --- &  --- & \snum{3.489067e-01} & \snum{4.143061e-01} & \fnum{3.151784e+01} (\fnum{9.179996e-01}) & ---                                 & ---                                 & $\ast$\fnum{3.433318e+01}\\
\runid & \nk\ (\texttt{ireg})     &  17 & 1848 &  --- &  899 & \snum{2.533550e-01} & \snum{4.020409e-02} & \fnum{6.181754e+01} (\fnum{8.821935e-01}) & ---                                 & \fnum{6.688521e+01} (\fnum{9.545139e-01}) & \fnum{7.007254e+01}\\
\runid & \nk\ (\texttt{2lrpcsym}) &  21 &  182 & 2376 &   60 & \snum{2.530559e-01} & \snum{4.750696e-02} & \fnum{5.684566e+00} (\fnum{1.778020e-01}) & \fnum{1.560684e+01} (\fnum{4.881510e-01}) & \fnum{4.686841e+00} (\fnum{1.465951e-01}) & \fnum{3.197133e+01}\\
\rowcolor{p1color!20}
\runid & \nk\ (\texttt{h0rpc})    &  17 &  281 &  --- &  115 & \snum{2.589842e-01} & \snum{3.882408e-02} & \fnum{8.961385e+00} (\fnum{5.904125e-01}) & ---                                 & \fnum{8.512277e+00} (\fnum{5.608234e-01}) & \fnum{1.517818e+01}\\
\runid & L-BFGS (20)              & 117 &  370 &  --- &  --- & \snum{2.561287e-01} & \snum{4.833765e-02} & \fnum{2.310059e+01} (\fnum{8.898685e-01}) & ---                                 & ---                                 & \fnum{2.595956e+01}\\
\runid & L-BFGS (50)              &  85 &  268 &  --- &  --- & \snum{2.562962e-01} & \snum{4.966618e-02} & \fnum{1.955701e+01} (\fnum{8.830836e-01}) & ---                                 & ---                                 & \fnum{2.214627e+01}\\
\runid & L-BFGS (100)             &  67 &  211 &  --- &  --- & \snum{2.560217e-01} & \snum{4.427210e-02} & \fnum{1.836379e+01} (\fnum{9.121842e-01}) & ---                                 & ---                                 & \fnum{2.013166e+01}\\
\bottomrule
\end{tabular}
\end{table}

\begin{table}
\caption{Convergence results for \gangmres\ for the \nirep\ data as a function of the regularization parameter $\alpha$. Setup, reported quantities and layout are as in \Cref{t:nirep-ngmres-sg}, with an additional column for the speedup against the quickest of the \rpgd, \lbfgs\ and \nk\ runs in \Cref{t:nirep-nk-rpdg-regul-new} for the same $\alpha$ (the highlighted row of that table). This reference is \lbfgs\ for all but the smallest regularization parameter, where it is \nk\ (\texttt{h0rpc}). Runtimes are the median of five repetitions.}
\label{t:nirep-ngmres-regul-new}
\tabadjust\setcounter{run}{0}
\addtolength{\tabcolsep}{-0.12em}\begin{tabular}{FrcrrrrRRRRRRr}\toprule
\multicolumn{7}{c}{} & \multicolumn{6}{c}{\cellcolor{gray!20}time (in seconds)} & \\
run & $w$ & $(\sigma,\tau)$ & \#iter & \#pdes & dist & grad & pdes & $q_k$ & $g$ & sg & ls & tts & speedup\\
\midrule
\multicolumn{14}{l}{$\alpha = 10^{-1}$} \\
\midrule
\runid &  10 & (5,1) &  10 &   71 & \snum{4.643657e-01} & \snum{4.486330e-02} & \fnum{1.693226e+00} (\fnum{8.476363e-01}) & \fnum{9.593345e-01} & \fnum{6.516849e-01} & \fnum{2.512163e-01} & \fnum{1.632175e-02} & \fnum{1.997585e+00} & \fnum{5.518965e-01}\\
\runid &  15 & (5,1) &  10 &   71 & \snum{4.643657e-01} & \snum{4.486330e-02} & \fnum{1.693713e+00} (\fnum{8.461634e-01}) & \fnum{9.641013e-01} & \fnum{6.498919e-01} & \fnum{2.525553e-01} & \fnum{1.536762e-02} & \fnum{2.001639e+00} & \fnum{5.507788e-01}\\
\rowcolor{p1color!20}
\runid &  20 & (5,1) &  10 &   71 & \snum{4.643657e-01} & \snum{4.486330e-02} & \fnum{1.608639e+00} (\fnum{8.390593e-01}) & \fnum{9.202955e-01} & \fnum{6.214469e-01} & \fnum{2.398873e-01} & \fnum{1.551529e-02} & \fnum{1.917193e+00} & \fnum{5.750388e-01}\\
\runid &  25 & (5,1) &  10 &   71 & \snum{4.643657e-01} & \snum{4.486330e-02} & \fnum{1.662798e+00} (\fnum{8.438713e-01}) & \fnum{9.491905e-01} & \fnum{6.378085e-01} & \fnum{2.478845e-01} & \fnum{1.573112e-02} & \fnum{1.970441e+00} & \fnum{5.594994e-01}\\
\runid &  10 & (4,2) &  13 &   82 & \snum{4.647640e-01} & \snum{6.350361e-03} & \fnum{1.951560e+00} (\fnum{8.414551e-01}) & \fnum{1.069904e+00} & \fnum{8.398485e-01} & \fnum{2.489418e-01} & \fnum{2.064225e-02} & \fnum{2.319268e+00} & \fnum{4.753483e-01}\\
\runid &  15 & (4,2) &  13 &   82 & \snum{4.647640e-01} & \snum{6.350361e-03} & \fnum{1.952010e+00} (\fnum{8.409893e-01}) & \fnum{1.074080e+00} & \fnum{8.414271e-01} & \fnum{2.491543e-01} & \fnum{2.027217e-02} & \fnum{2.321087e+00} & \fnum{4.749758e-01}\\
\runid &  20 & (4,2) &  13 &   82 & \snum{4.647640e-01} & \snum{6.350361e-03} & \fnum{1.949569e+00} (\fnum{8.306095e-01}) & \fnum{1.083081e+00} & \fnum{8.509034e-01} & \fnum{2.519751e-01} & \fnum{1.994025e-02} & \fnum{2.347155e+00} & \fnum{4.697007e-01}\\
\runid &  25 & (4,2) &  13 &   82 & \snum{4.647640e-01} & \snum{6.350361e-03} & \fnum{1.958417e+00} (\fnum{8.432137e-01}) & \fnum{1.075298e+00} & \fnum{8.422963e-01} & \fnum{2.491870e-01} & \fnum{2.010833e-02} & \fnum{2.322563e+00} & \fnum{4.746739e-01}\\
\midrule
\multicolumn{14}{l}{$\alpha = 5 \times 10^{-2}$} \\
\midrule
\runid &  10 & (5,1) &  13 &   93 & \snum{4.295166e-01} & \snum{3.344970e-02} & \fnum{2.180280e+00} (\fnum{8.327427e-01}) & \fnum{1.206320e+00} & \fnum{8.299197e-01} & \fnum{4.008834e-01} & \fnum{2.598850e-02} & \fnum{2.618191e+00} & \fnum{5.109751e-01}\\
\runid &  15 & (5,1) &  13 &   93 & \snum{4.295166e-01} & \snum{3.344970e-02} & \fnum{2.236631e+00} (\fnum{8.464462e-01}) & \fnum{1.218377e+00} & \fnum{8.460517e-01} & \fnum{4.102628e-01} & \fnum{2.573071e-02} & \fnum{2.642378e+00} & \fnum{5.062980e-01}\\
\runid &  20 & (5,1) &  13 &   93 & \snum{4.295166e-01} & \snum{3.344970e-02} & \fnum{2.221849e+00} (\fnum{8.462646e-01}) & \fnum{1.209355e+00} & \fnum{8.390135e-01} & \fnum{4.065800e-01} & \fnum{2.605896e-02} & \fnum{2.625478e+00} & \fnum{5.095570e-01}\\
\runid &  25 & (5,1) &  13 &   93 & \snum{4.295166e-01} & \snum{3.344970e-02} & \fnum{2.181770e+00} (\fnum{8.346729e-01}) & \fnum{1.207899e+00} & \fnum{8.281690e-01} & \fnum{4.006016e-01} & \fnum{2.586283e-02} & \fnum{2.613922e+00} & \fnum{5.118097e-01}\\
\runid &  10 & (4,2) &  13 &   88 & \snum{4.296434e-01} & \snum{4.439796e-02} & \fnum{2.048626e+00} (\fnum{8.453483e-01}) & \fnum{1.064322e+00} & \fnum{8.286338e-01} & \fnum{3.772490e-01} & \fnum{1.995021e-02} & \fnum{2.423411e+00} & \fnum{5.520445e-01}\\
\runid &  15 & (4,2) &  13 &   88 & \snum{4.296434e-01} & \snum{4.439796e-02} & \fnum{2.058901e+00} (\fnum{8.466286e-01}) & \fnum{1.066944e+00} & \fnum{8.302640e-01} & \fnum{3.795071e-01} & \fnum{2.013712e-02} & \fnum{2.431883e+00} & \fnum{5.501214e-01}\\
\rowcolor{p1color!20}
\runid &  20 & (4,2) &  13 &   88 & \snum{4.296434e-01} & \snum{4.439796e-02} & \fnum{2.033798e+00} (\fnum{8.444108e-01}) & \fnum{1.056506e+00} & \fnum{8.204771e-01} & \fnum{3.757651e-01} & \fnum{1.995062e-02} & \fnum{2.408541e+00} & \fnum{5.554526e-01}\\
\runid &  25 & (4,2) &  13 &   88 & \snum{4.296434e-01} & \snum{4.439796e-02} & \fnum{2.080592e+00} (\fnum{8.433196e-01}) & \fnum{1.084176e+00} & \fnum{8.407224e-01} & \fnum{3.845686e-01} & \fnum{2.012238e-02} & \fnum{2.467145e+00} & \fnum{5.422586e-01}\\
\midrule
\multicolumn{14}{l}{$\alpha = 10^{-2}$} \\
\midrule
\runid &  10 & (5,1) &  22 &  148 & \snum{3.767607e-01} & \snum{4.544734e-02} & \fnum{3.627257e+00} (\fnum{8.253249e-01}) & \fnum{2.150183e+00} & \fnum{1.408215e+00} & \fnum{5.137512e-01} & \fnum{7.117858e-02} & \fnum{4.394945e+00} & \fnum{5.782837e-01}\\
\runid &  15 & (5,1) &  21 &  140 & \snum{3.769557e-01} & \snum{3.177611e-02} & \fnum{3.414970e+00} (\fnum{8.161039e-01}) & \fnum{2.045129e+00} & \fnum{1.335717e+00} & \fnum{4.494183e-01} & \fnum{8.206179e-02} & \fnum{4.184479e+00} & \fnum{6.073695e-01}\\
\runid &  20 & (5,1) &  21 &  139 & \snum{3.770261e-01} & \snum{3.396793e-02} & \fnum{3.489431e+00} (\fnum{8.234871e-01}) & \fnum{2.054435e+00} & \fnum{1.369448e+00} & \fnum{4.284575e-01} & \fnum{8.610300e-02} & \fnum{4.237384e+00} & \fnum{5.997863e-01}\\
\runid &  25 & (5,1) &  21 &  139 & \snum{3.770261e-01} & \snum{3.396793e-02} & \fnum{3.460889e+00} (\fnum{8.258354e-01}) & \fnum{2.080781e+00} & \fnum{1.356099e+00} & \fnum{4.310067e-01} & \fnum{8.273708e-02} & \fnum{4.190773e+00} & \fnum{6.064573e-01}\\
\runid &  10 & (4,2) &  25 &  155 & \snum{3.775257e-01} & \snum{3.845394e-02} & \fnum{3.714191e+00} (\fnum{8.294152e-01}) & \fnum{2.101777e+00} & \fnum{1.572637e+00} & \fnum{5.055155e-01} & \fnum{6.548217e-02} & \fnum{4.478084e+00} & \fnum{5.675474e-01}\\
\runid &  15 & (4,2) &  24 &  141 & \snum{3.770943e-01} & \snum{4.422024e-02} & \fnum{3.445102e+00} (\fnum{8.229435e-01}) & \fnum{2.008588e+00} & \fnum{1.524539e+00} & \fnum{3.361602e-01} & \fnum{7.332571e-02} & \fnum{4.186317e+00} & \fnum{6.071029e-01}\\
\rowcolor{p1color!20}
\runid &  20 & (4,2) &  23 &  138 & \snum{3.768023e-01} & \snum{4.145654e-02} & \fnum{3.423065e+00} (\fnum{8.185509e-01}) & \fnum{1.994695e+00} & \fnum{1.481053e+00} & \fnum{3.618905e-01} & \fnum{8.073667e-02} & \fnum{4.181859e+00} & \fnum{6.077500e-01}\\
\runid &  25 & (4,2) &  23 &  138 & \snum{3.768023e-01} & \snum{4.145654e-02} & \fnum{3.442539e+00} (\fnum{8.227620e-01}) & \fnum{1.999132e+00} & \fnum{1.486803e+00} & \fnum{3.654438e-01} & \fnum{7.998287e-02} & \fnum{4.184125e+00} & \fnum{6.074209e-01}\\
\midrule
\multicolumn{14}{l}{$\alpha = 5 \times 10^{-3}$} \\
\midrule
\runid &  10 & (5,1) &  35 &  231 & \snum{3.536670e-01} & \snum{4.940570e-02} & \fnum{5.771275e+00} (\fnum{8.219446e-01}) & \fnum{3.464616e+00} & \fnum{2.237759e+00} & \fnum{8.054217e-01} & \fnum{1.342766e-01} & \fnum{7.021488e+00} & \fnum{4.993123e-01}\\
\runid &  15 & (5,1) &  33 &  214 & \snum{3.555438e-01} & \snum{3.402592e-02} & \fnum{5.416618e+00} (\fnum{8.158993e-01}) & \fnum{3.244109e+00} & \fnum{2.119459e+00} & \fnum{7.228356e-01} & \fnum{1.648857e-01} & \fnum{6.638831e+00} & \fnum{5.280922e-01}\\
\runid &  20 & (5,1) &  31 &  197 & \snum{3.549846e-01} & \snum{2.994036e-02} & \fnum{4.953447e+00} (\fnum{7.999212e-01}) & \fnum{3.044783e+00} & \fnum{1.976180e+00} & \fnum{5.634957e-01} & \fnum{1.758155e-01} & \fnum{6.192419e+00} & \fnum{5.661625e-01}\\
\runid &  25 & (5,1) &  31 &  196 & \snum{3.548669e-01} & \snum{3.057696e-02} & \fnum{5.029675e+00} (\fnum{7.988419e-01}) & \fnum{3.073809e+00} & \fnum{2.009118e+00} & \fnum{5.725768e-01} & \fnum{1.970920e-01} & \fnum{6.296207e+00} & \fnum{5.568297e-01}\\
\runid &  10 & (4,2) &  37 &  222 & \snum{3.549017e-01} & \snum{4.519146e-02} & \fnum{5.641222e+00} (\fnum{8.317851e-01}) & \fnum{3.306601e+00} & \fnum{2.392857e+00} & \fnum{6.390362e-01} & \fnum{1.110020e-01} & \fnum{6.782066e+00} & \fnum{5.169391e-01}\\
\runid &  15 & (4,2) &  33 &  199 & \snum{3.543522e-01} & \snum{4.791372e-02} & \fnum{5.002427e+00} (\fnum{8.176807e-01}) & \fnum{2.990571e+00} & \fnum{2.110300e+00} & \fnum{5.205958e-01} & \fnum{1.313368e-01} & \fnum{6.117824e+00} & \fnum{5.730657e-01}\\
\rowcolor{p1color!20}
\runid &  20 & (4,2) &  32 &  190 & \snum{3.551167e-01} & \snum{4.315396e-02} & \fnum{4.874544e+00} (\fnum{8.142519e-01}) & \fnum{2.902814e+00} & \fnum{2.077497e+00} & \fnum{4.814482e-01} & \fnum{1.507598e-01} & \fnum{5.986531e+00} & \fnum{5.856338e-01}\\
\runid &  25 & (4,2) &  32 &  190 & \snum{3.550026e-01} & \snum{4.089928e-02} & \fnum{4.870507e+00} (\fnum{8.081707e-01}) & \fnum{2.905321e+00} & \fnum{2.075246e+00} & \fnum{4.813205e-01} & \fnum{1.678343e-01} & \fnum{6.026582e+00} & \fnum{5.817419e-01}\\
\midrule
\multicolumn{14}{l}{$\alpha = 10^{-3}$} \\
\midrule
\runid &  10 & (5,1) &  81 &  528 & \snum{2.949264e-01} & \snum{4.627579e-02} & \fnum{1.608876e+01} (\fnum{8.499798e-01}) & \fnum{1.051013e+01} & \fnum{5.258149e+00} & \fnum{2.095278e+00} & \fnum{3.397134e-01} & \fnum{1.892840e+01} & \fnum{4.825464e-01}\\
\runid &  15 & (5,1) &  71 &  457 & \snum{2.933501e-01} & \snum{4.892110e-02} & \fnum{1.445179e+01} (\fnum{8.397717e-01}) & \fnum{9.662917e+00} & \fnum{4.557771e+00} & \fnum{1.643246e+00} & \fnum{4.133544e-01} & \fnum{1.720918e+01} & \fnum{5.307533e-01}\\
\runid &  20 & (5,1) &  67 &  430 & \snum{2.920950e-01} & \snum{4.693546e-02} & \fnum{1.336729e+01} (\fnum{8.199016e-01}) & \fnum{8.998106e+00} & \fnum{4.279242e+00} & \fnum{1.531596e+00} & \fnum{4.971276e-01} & \fnum{1.630353e+01} & \fnum{5.602365e-01}\\
\runid &  25 & (5,1) &  62 &  392 & \snum{2.925284e-01} & \snum{4.512947e-02} & \fnum{1.282343e+01} (\fnum{8.149312e-01}) & \fnum{8.838300e+00} & \fnum{3.986484e+00} & \fnum{1.312692e+00} & \fnum{5.379729e-01} & \fnum{1.573559e+01} & \fnum{5.804567e-01}\\
\runid &  10 & (4,2) &  97 &  569 & \snum{2.903014e-01} & \snum{4.545029e-02} & \fnum{1.722331e+01} (\fnum{8.465040e-01}) & \fnum{1.109950e+01} & \fnum{6.223354e+00} & \fnum{1.801964e+00} & \fnum{3.340892e-01} & \fnum{2.034640e+01} & \fnum{4.489164e-01}\\
\runid &  15 & (4,2) &  79 &  460 & \snum{2.931619e-01} & \snum{4.746359e-02} & \fnum{1.446046e+01} (\fnum{8.378644e-01}) & \fnum{9.653250e+00} & \fnum{5.066249e+00} & \fnum{1.331661e+00} & \fnum{3.762097e-01} & \fnum{1.725871e+01} & \fnum{5.292303e-01}\\
\runid &  20 & (4,2) &  73 &  427 & \snum{2.914395e-01} & \snum{4.357157e-02} & \fnum{1.333176e+01} (\fnum{8.244927e-01}) & \fnum{8.925061e+00} & \fnum{4.651686e+00} & \fnum{1.300690e+00} & \fnum{4.495470e-01} & \fnum{1.616965e+01} & \fnum{5.648751e-01}\\
\rowcolor{p1color!20}
\runid &  25 & (4,2) &  64 &  376 & \snum{2.925754e-01} & \snum{4.707619e-02} & \fnum{1.190423e+01} (\fnum{8.177442e-01}) & \fnum{8.151552e+00} & \fnum{3.989632e+00} & \fnum{1.018161e+00} & \fnum{4.766172e-01} & \fnum{1.455740e+01} & \fnum{6.274357e-01}\\
\midrule
\multicolumn{14}{l}{$\alpha = 5 \times 10^{-4}$} \\
\midrule
\runid &  10 & (5,1) & 127 &  816 & \snum{2.606810e-01} & \snum{4.535013e-02} & \fnum{3.369442e+01} (\fnum{8.821497e-01}) & \fnum{2.513733e+01} & \fnum{8.214924e+00} & \fnum{3.189350e+00} & \fnum{5.559913e-01} & \fnum{3.819581e+01} & \fnum{3.973781e-01}\\
\runid &  15 & (5,1) & 103 &  655 & \snum{2.596182e-01} & \snum{4.518046e-02} & \fnum{2.906857e+01} (\fnum{8.792218e-01}) & \fnum{2.231214e+01} & \fnum{6.642446e+00} & \fnum{2.332727e+00} & \fnum{5.827310e-01} & \fnum{3.306170e+01} & \fnum{4.590864e-01}\\
\runid &  20 & (5,1) &  91 &  577 & \snum{2.584386e-01} & \snum{4.655789e-02} & \fnum{2.723064e+01} (\fnum{8.688595e-01}) & \fnum{2.175466e+01} & \fnum{5.910651e+00} & \fnum{2.031516e+00} & \fnum{7.530636e-01} & \fnum{3.134067e+01} & \fnum{4.842965e-01}\\
\runid &  25 & (5,1) &  85 &  538 & \snum{2.590094e-01} & \snum{4.880622e-02} & \fnum{2.595057e+01} (\fnum{8.591126e-01}) & \fnum{2.030910e+01} & \fnum{5.528726e+00} & \fnum{1.920308e+00} & \fnum{8.642098e-01} & \fnum{3.020625e+01} & \fnum{5.024847e-01}\\
\runid &  10 & (4,2) & 129 &  753 & \snum{2.566528e-01} & \snum{4.970074e-02} & \fnum{3.184972e+01} (\fnum{8.854381e-01}) & \fnum{2.346953e+01} & \fnum{8.447620e+00} & \fnum{2.481813e+00} & \fnum{4.590702e-01} & \fnum{3.597058e+01} & \fnum{4.219608e-01}\\
\runid &  15 & (4,2) & 112 &  651 & \snum{2.575194e-01} & \snum{4.999209e-02} & \fnum{2.911479e+01} (\fnum{8.786745e-01}) & \fnum{2.212455e+01} & \fnum{7.258248e+00} & \fnum{2.029512e+00} & \fnum{5.711870e-01} & \fnum{3.313490e+01} & \fnum{4.580722e-01}\\
\rowcolor{p1color!20}
\runid &  20 & (4,2) &  94 &  542 & \snum{2.600450e-01} & \snum{4.834528e-02} & \fnum{2.609980e+01} (\fnum{8.729450e-01}) & \fnum{2.032594e+01} & \fnum{6.150895e+00} & \fnum{1.569438e+00} & \fnum{6.407057e-01} & \fnum{2.989856e+01} & \fnum{5.076557e-01}\\
\runid &  25 & (4,2) &  92 &  529 & \snum{2.595147e-01} & \snum{4.682473e-02} & \fnum{2.658866e+01} (\fnum{8.717471e-01}) & \fnum{2.085320e+01} & \fnum{5.919334e+00} & \fnum{1.495788e+00} & \fnum{7.447846e-01} & \fnum{3.050043e+01} & \fnum{4.976381e-01}\\
\bottomrule
\end{tabular}
\end{table}

The most important observations are:
\begin{enumerate*}[label=\it(\roman*)]
\item At this tolerance, \gangmres\ is not the fastest scheme. Against the fastest baseline for the same $\alpha$ (\lbfgs\ for $\alpha \ge 10^{-3}$, \nk\ (\texttt{h0rpc}) for $\alpha = 5\times 10^{-4}$), the speedup in \Cref{t:nirep-ngmres-regul-new} ranges from $0.40\times$ to $0.63\times$ with a median of $0.55\times$; the best \gangmres\ configuration reaches between $0.51\times$ and $0.63\times$. Against the \rpgd\ and \nk\ variants alone, \gangmres\ is competitive but not uniformly faster; the speedup then ranges from $0.40\times$ to $1.76\times$ with a median of $1.16\times$.
\item For $\alpha \le 10^{-2}$, the acceleration of \gangmres\ compared to the \rpgd\ is significant; for these values of $\alpha$, \rpgd\ does not converge within 200 iterations (see runs {15}, {22}, {29} and {36} in \Cref{t:nirep-nk-rpdg-regul-new}).
\item As $\alpha$ tends to zero, the performance of the \gangmres\ deteriorates significantly. For $\alpha \le 10^{-3}$ the best performing \nk\ scheme in \Cref{t:nirep-nk-rpdg-regul-new} is faster than every configuration of \gangmres; for $\alpha = 5 \times 10^{-4}$ by a factor of $2.0$ ($15.18$ vs.\ $29.90$ seconds).
\end{enumerate*}

We note that in practical applications, we typically perform a bisection search for an optimal $\alpha$ subject to bounds on the determinant of the deformation gradient. Likewise, we have designed a parameter continuation scheme that delivers faster convergence if we have identified an optimal regularization parameter. Both start with large regularization parameters, where the advantage of \gangmres\ is small or absent. For $\alpha = 10^{-1}$, \rpgd\ ($2.14$\,s) is slower than the four configurations of \gangmres\ with $p=(5,1)$ ($1.92$ to $2.00$\,s) and faster than the four with $p=(4,2)$ ($2.32$ to $2.35$\,s), and \lbfgs\ is faster than \gangmres\ by a factor between $1.49$ and $1.80$ for all $\alpha$ we consider. Whether \gangmres\ pays off within such a continuation scheme remains open. More details about the search for an optimal $\alpha$ and the continuation approach can be found in~\cite{mang2015:inexact, mang2019:claire, mang2024:claire}.

\subsection{Mesh Convergence}\label{s:results:mesh}

We assess the performance of the proposed scheme as a function of the mesh size.

We select $n_i$ in $\{64,128,256,512\}$. The associated number of time steps $n_t$ for the time integrator are $\{4,8,16,32\}$. We use the \rect\ dataset. We perform two experiments. First, we fix the smoothing of the input data to $\eta{} = 1$. This implies that the edges of the images become sharper as we increase the resolution. In the second experiment, we increase $\eta$ as the resolution increases. We expect mesh independent convergence for \nk\ methods for the latter setup. We set the regularization parameter to $\alpha = 10^{-3}$.

We compare the performance of the proposed \gangmres\ method to the \nk\ scheme and \lbfgs. The results can be found in \Cref{t:ngmres-mesh-new}. We use a smaller tolerance of $\epsilon_{\scr{rel}} = 10^{-3}$ for these experiments.

\begin{table}[!tp]
\caption{Convergence results for \gangmres\ (with the safeguard of the accelerated iterate), \lbfgs\ and \nk\ for the \rect\ data as a function of the mesh size; $n_{\scr{iter}} = 200$ and $\epsilon_{\scr{rel}} = 10^{-3}$. Reported quantities and layout combine those of \Cref{t:h2-nk+rpgd_hands} and \Cref{t:nirep-ngmres-sg}. We additionally report the mesh size $n_i$, the number of time steps $n_t$ and the smoothing parameter $\eta$, and, as in \Cref{t:h2-nk+rpgd_hands}, the PDE solves and PDE solve time of the preconditioner separately (\#pdes$_{\scr{pc}}$ and pdes$_{\scr{pc}}$). Runtimes with $\ast$ indicate that the scheme did not reach the tolerance; all of these runs stop before $n_{\scr{iter}}$ iterations because the line search fails. --- marks a quantity that does not apply. No row is highlighted, since the runs are for different mesh sizes, i.e., different problems. Runtimes are the median of five repetitions.}
\label{t:ngmres-mesh-new}
\tabadjust\setcounter{run}{0}
\addtolength{\tabcolsep}{-0.18em}\begin{tabular}{FrrrrrrrrRRRRRRE}\toprule
\multicolumn{9}{c}{} & \multicolumn{7}{c}{\cellcolor{gray!20}time (in seconds)} \\
run & $n_i$ & $n_t$ & $\eta$ & \#iter & \#pdes & \#pdes$_{\scr{pc}}$ & dist & grad & pdes & pdes$_{\scr{pc}}$ & $q_k$ & $g$ & sg & ls & tts\\
\midrule
\multicolumn{16}{l}{\gangmres(25;4,2)} \\
\midrule
\runid & 64 & 4 & 1 & 19 & 118 & --- & \snum{1.887416e-02} & \snum{9.013160e-04} & \fnum{3.285882e-01} (\fnum{8.207232e-01}) & --- & \fnum{1.960216e-01} & \fnum{1.406165e-01} & \fnum{3.348375e-02} & \fnum{4.157458e-03} & \fnum{4.003642e-01} \\
\runid & 128 & 8 &  & 25 & 152 & --- & \snum{2.364025e-02} & \snum{7.484925e-04} & \fnum{1.497753e+00} (\fnum{8.308124e-01}) & --- & \fnum{8.436798e-01} & \fnum{6.830670e-01} & \fnum{1.767886e-01} & \fnum{1.856779e-02} & \fnum{1.802757e+00} \\
\runid & 256 & 16 &  & 56 & 335 & --- & \snum{2.742452e-02} & \snum{7.697989e-04} & \fnum{1.744198e+01} (\fnum{8.839634e-01}) & --- & \fnum{8.631609e+00} & \fnum{7.851490e+00} & \fnum{2.636110e+00} & \fnum{1.328006e-01} & \fnum{1.973157e+01} \\
\runid & 512 & 32 &  & 75 & 444 & --- & \snum{3.069079e-02} & \snum{9.516854e-04} & \fnum{2.179526e+02} (\fnum{9.374665e-01}) & --- & \fnum{1.035141e+02} & \fnum{9.415645e+01} & \fnum{3.115576e+01} & \fnum{6.759226e-01} & \fnum{2.324911e+02} \\
\midrule
\runid & 128 & 8 & 2 & 16 & 107 & --- & \snum{1.884013e-02} & \snum{4.141357e-04} & \fnum{1.015811e+00} (\fnum{8.291840e-01}) & --- & \fnum{6.152065e-01} & \fnum{4.266077e-01} & \fnum{1.120505e-01} & \fnum{1.204063e-02} & \fnum{1.225073e+00} \\
\runid & 256 & 16 & 4 & 22 & 139 & --- & \snum{1.889019e-02} & \snum{2.147334e-04} & \fnum{7.307169e+00} (\fnum{8.833501e-01}) & --- & \fnum{4.096099e+00} & \fnum{3.103730e+00} & \fnum{7.085811e-01} & \fnum{4.930596e-02} & \fnum{8.272109e+00} \\
\runid & 512 & 32 & 8 & 15 & 96 & --- & \snum{1.890551e-02} & \snum{9.956908e-04} & \fnum{4.615003e+01} (\fnum{9.393505e-01}) & --- & \fnum{2.432281e+01} & \fnum{1.885878e+01} & \fnum{4.391622e+00} & \fnum{7.270171e-02} & \fnum{4.912973e+01} \\
\midrule
\multicolumn{16}{l}{\nk\ (\texttt{2lrpcsym})} \\
\midrule
\runid & 64 & 4 & 1 & 7 & 61 & 320 & \snum{1.882281e-02} & \snum{6.125471e-04} & \fnum{2.050495e-01} (\fnum{2.563049e-01}) & \fnum{4.339197e-01} & --- & --- & --- & --- & \fnum{8.000217e-01} \\
\runid & 128 & 8 &  & 10 & 83 & 252 & \snum{2.366716e-02} & \snum{5.300392e-04} & \fnum{1.040747e+00} (\fnum{4.104748e-01}) & \fnum{9.447994e-01} & --- & --- & --- & --- & \fnum{2.535471e+00} \\
\runid & 256 & 16 &  & 18 & 155 & 518 & \snum{2.750966e-02} & \snum{6.839325e-04} & \fnum{1.117472e+01} (\fnum{4.956473e-01}) & \fnum{7.773636e+00} & --- & --- & --- & --- & \fnum{2.254571e+01} \\
\runid & 512 & 32 &  & 32 & 281 & 948 & \snum{3.090339e-02} & \snum{9.523001e-04} & \fnum{1.918597e+02} (\fnum{6.403362e-01}) & \fnum{8.007737e+01} & --- & --- & --- & --- & \fnum{2.996233e+02} \\
\midrule
\runid & 128 & 8 & 2 & 6 & 47 & 104 & \snum{1.884541e-02} & \snum{5.676859e-04} & \fnum{5.363911e-01} (\fnum{4.700391e-01}) & \fnum{3.603000e-01} & --- & --- & --- & --- & \fnum{1.141163e+00} \\
\runid & 256 & 16 & 4 & 6 & 47 & 96 & \snum{1.891021e-02} & \snum{8.196733e-04} & \fnum{2.859278e+00} (\fnum{5.805302e-01}) & \fnum{1.292804e+00} & --- & --- & --- & --- & \fnum{4.925288e+00} \\
\runid & 512 & 32 & 8 & 6 & 47 & 82 & \snum{1.897573e-02} & \snum{7.581393e-04} & \fnum{2.728078e+01} (\fnum{7.483297e-01}) & \fnum{6.075417e+00} & --- & --- & --- & --- & \fnum{3.645557e+01} \\
\midrule
\multicolumn{16}{l}{\nk\ (\texttt{h0rpc})} \\
\midrule
\runid & 64 & 4 & 1 & 10 & 175 & --- & \snum{1.880647e-02} & \snum{5.147347e-04} & \fnum{6.257411e-01} (\fnum{5.899815e-01}) & --- & --- & --- & --- & --- & \fnum{1.060611e+00} \\
\runid & 128 & 8 &  & 20 & 394 & --- & \snum{2.363949e-02} & \snum{7.448134e-04} & \fnum{5.552599e+00} (\fnum{6.092222e-01}) & --- & --- & --- & --- & --- & \fnum{9.114243e+00} \\
\runid & 256 & 16 &  & 17 & 266 & --- & \snum{5.067517e-02} & \snum{1.645584e-01} & \fnum{1.898162e+01} (\fnum{8.697758e-01}) & --- & --- & --- & --- & --- & $\ast$\fnum{2.182358e+01} \\
\runid & 512 & 32 &  & 33 & 551 & --- & \snum{4.173397e-02} & \snum{2.011142e-01} & \fnum{3.888744e+02} (\fnum{9.167273e-01}) & --- & --- & --- & --- & --- & $\ast$\fnum{4.241985e+02} \\
\midrule
\runid & 128 & 8 & 2 & 11 & 192 & --- & \snum{1.882525e-02} & \snum{2.138460e-04} & \fnum{2.609190e+00} (\fnum{6.223414e-01}) & --- & --- & --- & --- & --- & \fnum{4.192538e+00} \\
\runid & 256 & 16 & 4 & 11 & 160 & --- & \snum{1.888197e-02} & \snum{8.962567e-04} & \fnum{1.168698e+01} (\fnum{8.134988e-01}) & --- & --- & --- & --- & --- & \fnum{1.436631e+01} \\
\runid & 512 & 32 & 8 & 12 & 209 & --- & \snum{1.893638e-02} & \snum{2.101383e-04} & \fnum{1.478413e+02} (\fnum{8.901812e-01}) & --- & --- & --- & --- & --- & \fnum{1.660800e+02} \\
\midrule
\multicolumn{16}{l}{L-BFGS (20)} \\
\midrule
\runid & 64 & 4 & 1 & 23 & 105 & --- & \snum{1.879320e-02} & \snum{1.339882e-03} & \fnum{2.770660e-01} (\fnum{8.058984e-01}) & --- & --- & --- & --- & --- & $\ast$\fnum{3.437976e-01} \\
\runid & 128 & 8 &  & 30 & 99 & --- & \snum{2.363582e-02} & \snum{7.101534e-04} & \fnum{9.187062e-01} (\fnum{7.897041e-01}) & --- & --- & --- & --- & --- & \fnum{1.163355e+00} \\
\runid & 256 & 16 &  & 58 & 180 & --- & \snum{2.746257e-02} & \snum{8.411148e-04} & \fnum{9.229252e+00} (\fnum{8.464928e-01}) & --- & --- & --- & --- & --- & \fnum{1.090293e+01} \\
\runid & 512 & 32 &  & 54 & 169 & --- & \snum{3.074520e-02} & \snum{9.623463e-04} & \fnum{8.067210e+01} (\fnum{9.206269e-01}) & --- & --- & --- & --- & --- & \fnum{8.762735e+01} \\
\midrule
\runid & 128 & 8 & 2 & 22 & 69 & --- & \snum{1.878296e-02} & \snum{8.792546e-04} & \fnum{6.470735e-01} (\fnum{7.955314e-01}) & --- & --- & --- & --- & --- & \fnum{8.133853e-01} \\
\runid & 256 & 16 & 4 & 22 & 68 & --- & \snum{1.889138e-02} & \snum{9.307206e-04} & \fnum{3.470346e+00} (\fnum{8.645246e-01}) & --- & --- & --- & --- & --- & \fnum{4.014167e+00} \\
\runid & 512 & 32 & 8 & 22 & 69 & --- & \snum{1.892359e-02} & \snum{8.563473e-04} & \fnum{3.325369e+01} (\fnum{9.294307e-01}) & --- & --- & --- & --- & --- & \fnum{3.577856e+01} \\
\midrule
\multicolumn{16}{l}{L-BFGS (50)} \\
\midrule
\runid & 64 & 4 & 1 & 23 & 105 & --- & \snum{1.879320e-02} & \snum{1.339882e-03} & \fnum{2.729347e-01} (\fnum{8.019513e-01}) & --- & --- & --- & --- & --- & $\ast$\fnum{3.403382e-01} \\
\runid & 128 & 8 &  & 29 & 92 & --- & \snum{2.363362e-02} & \snum{6.115660e-04} & \fnum{8.697930e-01} (\fnum{7.939609e-01}) & --- & --- & --- & --- & --- & \fnum{1.095511e+00} \\
\runid & 256 & 16 &  & 50 & 156 & --- & \snum{2.745422e-02} & \snum{9.556476e-04} & \fnum{7.963214e+00} (\fnum{8.474792e-01}) & --- & --- & --- & --- & --- & \fnum{9.396353e+00} \\
\runid & 512 & 32 &  & 56 & 173 & --- & \snum{3.073424e-02} & \snum{9.579554e-04} & \fnum{8.639300e+01} (\fnum{9.203975e-01}) & --- & --- & --- & --- & --- & \fnum{9.386488e+01} \\
\midrule
\runid & 128 & 8 & 2 & 22 & 69 & --- & \snum{1.878296e-02} & \snum{8.792546e-04} & \fnum{6.449744e-01} (\fnum{8.013505e-01}) & --- & --- & --- & --- & --- & \fnum{8.048593e-01} \\
\runid & 256 & 16 & 4 & 22 & 68 & --- & \snum{1.889138e-02} & \snum{9.307206e-04} & \fnum{3.451483e+00} (\fnum{8.655212e-01}) & --- & --- & --- & --- & --- & \fnum{3.987751e+00} \\
\runid & 512 & 32 & 8 & 22 & 69 & --- & \snum{1.892359e-02} & \snum{8.563473e-04} & \fnum{3.341922e+01} (\fnum{9.302458e-01}) & --- & --- & --- & --- & --- & \fnum{3.592515e+01} \\
\bottomrule
\end{tabular}
\end{table}

The most important observations are:
\begin{enumerate*}[label=\it(\roman*)]
\item The \gangmres\ scheme remains competitive in terms of runtime and iteration count as the mesh size decreases. Compared to \nk\ with the two-level preconditioner, \gangmres\ is faster in the four runs with a fixed smoothing parameter (by factors between $1.1$ and $2.0$) and slower in the three runs in which the smoothing increases with the mesh (by up to a factor of $1.7$). Counting the PDE solves of the preconditioner, \gangmres\ requires fewer PDE solves than \nk\ in all seven runs.
\item The ordering of the two \nk\ variants is reversed with respect to the image data. \nk\ with the zero-velocity preconditioner, the fastest of the \nk\ variants on \hands\ and \nirep, does not reach the tolerance at the two finest meshes of the fixed-smoothing experiment, and \gangmres\ is faster than it at each of the five mesh points at which it converges. This preconditioner is the reduced space Hessian evaluated at $v = 0$. With fixed smoothing, the edges of the data sharpen as the mesh is refined, so this approximation deteriorates with resolution; when the smoothing grows with the mesh, the scheme converges at every resolution. Its advantage is therefore specific to the image registration problems.
\item \lbfgs\ requires fewer PDE solves than \gangmres\ in all runs and is faster than both \gangmres\ and \nk\ wherever it reaches the tolerance, i.e., except for $n_i = 64$.
\item If we increase the smoothness parameter $\eta$ as we refine the mesh, the iteration counts of \gangmres\ and \nk{} depend much less on the mesh than for a fixed $\eta$ (\gangmres: $16$, $22$, $15$ instead of $19$ to $75$; \nk: $6$, $6$, $6$ instead of $7$ to $32$).
\end{enumerate*}

\subsection{Incompressible Diffeomorphisms}\label{s:stokes-flow}

We assess the performance of the proposed \gangmres\ algorithm for transport-dominated PDE-constrained optimization problems governed by incompressible flows.

In our past work, we have extended the problem formulation in \cref{e:varopt} to include PDE constraints for the divergence of $v$~\cite{mang2015:inexact, mang2016:constrained}. Adding the constraint $\idiv v = 0$ renders the flow incompressible. Similar formulations have been considered in~\cite{chen2011:image, ruhnau2007:optical}. We switch from $H^2$ to $H^3$ regularity for $v$, i.e., we set $s=3$ in~\cref{e:regop}. The variational problem is given by
\begin{subequations}\label{e:varopt-stokes}
\begin{align}
\underset{m \in \fs{M}_{\scr{ad}}, \, v \in \fs{V}_{\scr{ad}}}{\operatorname{minimize}}\quad
&
\frac{1}{2}\int_{\Omega} (m(x,1) - m_1(x) )^2\, \mathrm{d}x
+ \frac{\alpha}{2}\int_{\Omega} \fs{L} v(x) \cdot v(x)\, \mathrm{d}x
\\
\begin{aligned}
\text{subject to}\,\,\, \\ \\ \\
\end{aligned}
\quad&
\begin{aligned}\label{e:varopt:stokes-constraints}
\partial_t m(x,t) + \igrad m(x,t) \cdot v(x)  & = 0
&& \text{in}\,\,\Omega \times (0,1],
\\
m(x,0) &=  m_0(x)
&& \text{in}\,\,\Omega,
\\
\idiv v(x) &= 0
&& \text{in}\,\,\Omega.
\end{aligned}
\end{align}
\end{subequations}

We introduce a Lagrange multiplier $p : \Omega \to \mathbb{R}$ (a pressure) for the incompressibility constraint in~\cref{e:varopt:stokes-constraints}. Taking variations of the associated Lagrangian with respect to $v$ yields the Stokes-like system
\begin{equation}\label{e:stokes-system}
\alpha \fs{L} v(x) + \int_0^1 \lambda(x,t) \igrad m(x,t)\,\mathrm{d}t + \igrad p(x) = 0,
\qquad
\idiv v(x) = 0,
\qquad \text{in}\,\,\Omega.
\end{equation}

\noindent Since $\fs{L}$ commutes with $\idiv$ on the torus, we can eliminate the pressure from \cref{e:stokes-system}: Applying $\idiv$ to the first equation yields $\ilap p = -\idiv \int_0^1 \lambda \igrad m \,\mathrm{d}t$, so that the reduced gradient becomes
\begin{equation}\label{e:reduced-gradient-stokes}
g(v) \defeq \alpha \fs{L} v(x) + \fs{P}\left(\int_0^1 \lambda(x,t) \igrad m(x,t)\,\mathrm{d}t\right),
\qquad
\fs{P} \defeq \fs{I} - \igrad \ilap^{-1} \idiv,
\end{equation}

\noindent where $\fs{P}$ denotes the Leray projection onto the space of divergence-free vector fields and $\fs{I}$ is the identity. With $\beta_{\scr{lap}}=1$ and $\beta_{\scr{id}}=0$ we have $\fs{L} = (-\ilap)^3$; that is, the first term in~\cref{e:reduced-gradient-stokes} is a tri-harmonic operator applied to $v$, with spectral entries $|\xi|^6$. Both $\fs{L}$ and $\fs{P}$ are diagonal in the spectral basis, so applying them (and inverting $\fs{L}$) amounts to FFTs and a diagonal scaling. The spatial mesh is of size $256 \times 256$. The number of time steps is set to $n_t = 16$. The tolerance for the optimizer is set to $\epsilon_{\scr{rel}} = 10^{-3}$. No pre-smoothing is applied to the data. The regularization parameter $\alpha$ is set to $\alpha = 10^{-4}$. The data set is generated synthetically using smooth sinusoidal functions.

We report convergence results in \Cref{t:cinfty-stokes-new}.

\begin{table}
\caption{Convergence results for \gangmres\ for an incompressible Stokes flow. Reported quantities and layout are as in \Cref{t:nirep-ngmres-sg}. We use $w\in\{1,2,5,10\}$ instead of $w\in\{10,15,20,25\}$, since on these smooth data the scheme converges in at most eleven iterations for $w = 10$, so a window of ten is truncated in at most one iteration. We use $n_{\scr{iter}} = 200$ and $\epsilon_{\scr{rel}} = 10^{-3}$. Runtimes are the median of five repetitions.}
\label{t:cinfty-stokes-new}
\tabadjust\setcounter{run}{0}
\begin{tabular}{FrcrrrrRRRRRE}\toprule
\multicolumn{7}{c}{} & \multicolumn{6}{c}{\cellcolor{gray!20}time (in seconds)} \\
run & $w$ & $(\sigma,\tau)$ & \#iter & \#pdes & dist & grad & pdes & $q_k$ & $g$ & sg & ls & tts \\
\midrule
\runid &   1 &   (5,5) & 101 & 471 & \snum{8.464231e-04} & \snum{8.492680e-04} & \fnum{1.887436e+01} (\fnum{7.779258e-01}) & \fnum{1.020365e+01} & \fnum{1.167342e+01} & \fnum{1.733049e+00} & \fnum{2.040929e-02} & \fnum{2.426241e+01}\\
\runid &   2 &         & 84 & 409 & \snum{8.536313e-04} & \snum{2.811919e-04} & \fnum{1.635873e+01} (\fnum{7.779072e-01}) & \fnum{9.054237e+00} & \fnum{9.768279e+00} & \fnum{1.645414e+00} & \fnum{2.628700e-02} & \fnum{2.102916e+01}\\
\runid &   5 &         & 11 & 75 & \snum{8.664755e-04} & \snum{9.384124e-04} & \fnum{2.937817e+00} (\fnum{7.856404e-01}) & \fnum{1.979415e+00} & \fnum{1.273433e+00} & \fnum{3.089412e-01} & \fnum{7.204250e-03} & \fnum{3.739391e+00}\\
\runid &  10 &         & 11 & 75 & \snum{8.465618e-04} & \snum{2.502912e-04} & \fnum{2.956062e+00} (\fnum{7.860811e-01}) & \fnum{1.989732e+00} & \fnum{1.272640e+00} & \fnum{2.883396e-01} & \fnum{9.395667e-03} & \fnum{3.760505e+00}\\
\midrule
\runid &   1 &   (4,2) & 79 & 388 & \snum{8.355768e-04} & \snum{9.729589e-04} & \fnum{1.554143e+01} (\fnum{7.796727e-01}) & \fnum{8.609946e+00} & \fnum{9.156115e+00} & \fnum{1.596634e+00} & \fnum{2.084371e-02} & \fnum{1.993327e+01}\\
\runid &   2 &         & 56 & 283 & \snum{8.453996e-04} & \snum{9.863745e-04} & \fnum{1.140193e+01} (\fnum{7.777281e-01}) & \fnum{6.535981e+00} & \fnum{6.470411e+00} & \fnum{1.305141e+00} & \fnum{2.353575e-02} & \fnum{1.466056e+01}\\
\runid &   5 &         & 7 & 57 & \snum{8.337199e-04} & \snum{4.708061e-04} & \fnum{2.207895e+00} (\fnum{7.917133e-01}) & \fnum{1.552004e+00} & \fnum{8.292686e-01} & \fnum{2.463834e-01} & \fnum{5.251458e-03} & \fnum{2.788756e+00}\\
\runid &  10 &         & 7 & 57 & \snum{8.371182e-04} & \snum{4.506673e-04} & \fnum{2.188615e+00} (\fnum{7.889353e-01}) & \fnum{1.542095e+00} & \fnum{8.108871e-01} & \fnum{2.465047e-01} & \fnum{5.488792e-03} & \fnum{2.774137e+00}\\
\midrule
\runid &   1 &   (2,4) & 84 & 407 & \snum{8.300867e-04} & \snum{9.927467e-04} & \fnum{1.678295e+01} (\fnum{7.767255e-01}) & \fnum{9.664856e+00} & \fnum{9.939282e+00} & \fnum{1.241464e+00} & \fnum{1.464600e-02} & \fnum{2.160732e+01}\\
\runid &   2 &         & 139 & 593 & \snum{8.505497e-04} & \snum{9.050935e-04} & \fnum{2.425669e+01} (\fnum{7.771359e-01}) & \fnum{1.258502e+01} & \fnum{1.626142e+01} & \fnum{1.340499e+00} & \fnum{3.075729e-02} & \fnum{3.121293e+01}\\
\runid &   5 &         & 7 & 48 & \snum{8.491663e-04} & \snum{4.865818e-04} & \fnum{1.887345e+00} (\fnum{7.921341e-01}) & \fnum{1.260421e+00} & \fnum{8.239904e-01} & \fnum{1.315369e-01} & \fnum{2.718542e-03} & \fnum{2.382608e+00}\\
\rowcolor{p1color!20}
\runid &  10 &         & 7 & 48 & \snum{8.562708e-04} & \snum{4.491472e-04} & \fnum{1.873235e+00} (\fnum{7.902080e-01}) & \fnum{1.232285e+00} & \fnum{8.419117e-01} & \fnum{1.228422e-01} & \fnum{3.354042e-03} & \fnum{2.370559e+00}\\
\midrule
\runid &   1 &   (6,3) & 105 & 499 & \snum{8.381605e-04} & \snum{9.986108e-04} & \fnum{2.074963e+01} (\fnum{7.827817e-01}) & \fnum{1.158311e+01} & \fnum{1.247699e+01} & \fnum{1.980408e+00} & \fnum{2.602100e-02} & \fnum{2.650756e+01}\\
\runid &   2 &         & 77 & 378 & \snum{8.491997e-04} & \snum{9.640430e-04} & \fnum{1.507088e+01} (\fnum{7.766358e-01}) & \fnum{8.468996e+00} & \fnum{8.888062e+00} & \fnum{1.605118e+00} & \fnum{2.605617e-02} & \fnum{1.940534e+01}\\
\runid &   5 &         & 6 & 57 & \snum{8.407736e-04} & \snum{4.003030e-04} & \fnum{2.178594e+00} (\fnum{7.898549e-01}) & \fnum{1.601756e+00} & \fnum{6.952371e-01} & \fnum{3.047663e-01} & \fnum{7.083583e-03} & \fnum{2.758220e+00}\\
\runid &  10 &         & 6 & 57 & \snum{8.407736e-04} & \snum{4.003030e-04} & \fnum{2.142730e+00} (\fnum{7.875506e-01}) & \fnum{1.584788e+00} & \fnum{6.937176e-01} & \fnum{2.873808e-01} & \fnum{6.709000e-03} & \fnum{2.720753e+00}\\
\midrule
\runid &   1 &   (3,6) & 64 & 308 & \snum{8.406279e-04} & \snum{9.971055e-04} & \fnum{1.205366e+01} (\fnum{7.738940e-01}) & \fnum{6.846389e+00} & \fnum{7.388132e+00} & \fnum{1.021938e+00} & \fnum{1.071992e-02} & \fnum{1.557534e+01}\\
\runid &   2 &         & 27 & 149 & \snum{8.309355e-04} & \snum{9.827102e-04} & \fnum{5.770249e+00} (\fnum{7.829898e-01}) & \fnum{3.304059e+00} & \fnum{3.119423e+00} & \fnum{7.325163e-01} & \fnum{7.273125e-03} & \fnum{7.369507e+00}\\
\runid &   5 &         & 56 & 266 & \snum{8.385117e-04} & \snum{8.949339e-04} & \fnum{1.044351e+01} (\fnum{7.700489e-01}) & \fnum{5.913451e+00} & \fnum{6.477389e+00} & \fnum{7.824541e-01} & \fnum{3.136175e-02} & \fnum{1.356214e+01}\\
\runid &  10 &         & 10 & 63 & \snum{8.584690e-04} & \snum{3.297686e-04} & \fnum{2.386578e+00} (\fnum{7.796416e-01}) & \fnum{1.565239e+00} & \fnum{1.157560e+00} & \fnum{1.644977e-01} & \fnum{5.783417e-03} & \fnum{3.061121e+00}\\
\bottomrule
\end{tabular}
\end{table}

The most important observation is that the \gangmres\ scheme remains effective for the reformulation of our problem to account for incompressible transport maps. The solver remains effective for a range of hyperparameter values, providing excellent agreement between the transported intensities $m$ at time $t=1$ and the reference image $m_1$. The determinant of the inverse deformation gradient $\det\igrad y^{-1}$ associated with the computed flow map is equal to 1 to high numerical accuracy.

\subsection{Mass-Preserving Transport}\label{s:optimal-transport}

We assess the performance of the \gangmres\ scheme for  transport-dominated PDE-constrained optimization problems governed by the continuity equation. This formulation is related to optimal mass transport.

We consider two probability distributions $\pi_0 : \Omega \to [0,1]$, $\pi_1 : \Omega \to [0,1]$ that integrate to 1. Our formulation is related to the classical Monge--Kantorovich problem. In the classical formulation the density functions are assumed to have equal masses; in our formulation, density functions are allowed to yield masses that are not exactly equal. In general, we seek a map $y : \mathbb{R}^d \to \mathbb{R}^d$ (the deformation map) such that the pushforward $y_\# \pi_0 \approx \pi_1$. Likewise to the other models considered in this work, we do not parameterize $y$ directly; we model the transport by a smooth velocity field $v$. The key difference to the problem formulation in \cref{e:varopt} is that the sought after transport map ought to be mass-preserving. To do so, we consider the continuity equation as a PDE constraint. Related formulations have been considered in~\cite{mang2017:lagrangian, benzi2011:preconditioning, colombo2011:control, benamou2000:computational, ambrosio2014:continuity}. A key difference to many formulations for optimal mass transport is that our variational regularization model ensures that the computed transport map is a diffeomorphism. The variational problem formulation is given by
\begin{subequations}\label{e:varopt-optimal-transport}
\begin{align}
\underset{\pi \in \fs{P}_{\scr{ad}}, \, v \in \fs{V}_{\scr{ad}}}{\operatorname{minimize}}\quad
&
\frac{1}{2}\int_{\Omega} (\pi(x,1) - \pi_1(x) )^2\, \mathrm{d}x
+ \frac{\alpha}{2}\int_{\Omega} \fs{L} v(x) \cdot v(x)\, \mathrm{d}x
\\
\begin{aligned}
\text{subject to}\,\,\, \\ \\
\end{aligned}
\quad&
\begin{aligned}\label{e:varopt:continuity}
\partial_t \pi(x,t) + \idiv (\pi(x,t) v(x))  & = 0
&& \text{in}\,\,\Omega \times (0,1],
\\
\pi(x,0) &=  \pi_0(x)
&& \text{in}\,\,\Omega.
\end{aligned}
\end{align}
\end{subequations}

The state equation \cref{e:varopt:continuity} models a mass-preserving transport map for $\pi_0(x)$ subjected to $v$. The first term of the objective functional is a squared $L^2$-distance that measures the proximity of $\pi$ at time $t=1$ (terminal state) and the density $\pi_1$. The regularization model stipulates $H^2$-regularity on $v$, i.e., we set $s=2$ in~\cref{e:regop}. The optimality conditions are a biharmonic PDE.
The admissible set $\fs{V}_{\scr{ad}}$ for the control variable is the one introduced in~\cref{e:admissible}. The admissible set for the state variable is
\begin{equation}\label{e:admissible-ot}
\fs{P}_{\scr{ad}} \defeq \Big\{ \pi \in C(\bar{\Omega}\times[0,1]) \,:\, \pi \ge 0 \,\,\text{and}\,\, \int_{\Omega} \pi(x,t)\,\mathrm{d}x = \int_{\Omega} \pi_0(x)\,\mathrm{d}x \,\,\text{for all}\,\, t \in [0,1] \Big\};
\end{equation}

\noindent that is, $\fs{P}_{\scr{ad}}$ consists of the nonnegative functions of constant mass. We use $C(\bar{\Omega}\times[0,1])$ instead of $\fs{M}_{\scr{ad}}$, since~\cref{e:varopt:continuity} is in conservative form, so that its solution involves the determinant of the Jacobian of the associated deformation map. For periodic boundary conditions, the continuity equation~\cref{e:varopt:continuity} preserves both properties, so $\fs{P}_{\scr{ad}}$ is invariant under the transport.

We generate two probability densities $\pi_0$ and $\pi_1$ of equal mass. The spatial mesh is of size $256 \times 256$. The number of time steps is set to $n_t = 16$. The tolerance for the optimizer is set to $\epsilon_{\scr{rel}} = 10^{-3}$. No pre-smoothing is applied to the data. The regularization parameter $\alpha$ is set to $\alpha = 10^{-3}$.

\begin{figure}
\centering
\includegraphics[width=0.6\textwidth]{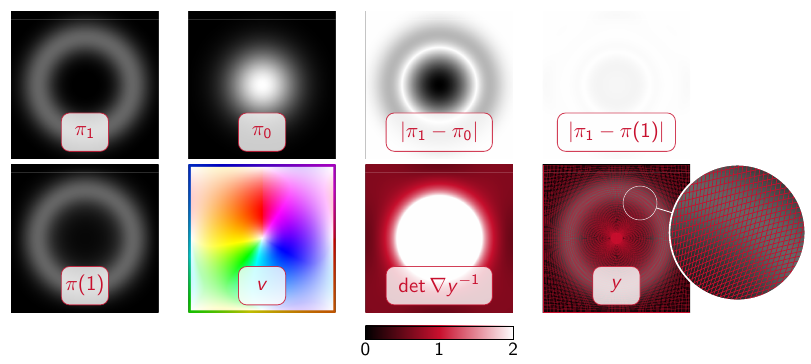}
\caption{We show exemplary results for mass preserving transport.  The results were obtained with \gangmres. Top row (from left to right): ($i$) the target density $\pi_1$; $(ii)$ the template density $\pi_0$ (probability density to be transported), $(iii)$ the residual differences between $\pi_0$ and $\pi_1$ (white: small difference; black: large difference); and $(iv)$ the residual differences between the terminal state $\pi$ at $t=1$ and $\pi_1$ after solving for the optimal $v$. Bottom row (from left to right): $(i)$ final state $\pi$ at $t=1$; $(ii)$ optimal control variable $v$ (color indicates orientation); $(iii)$ determinant of the inverse deformation gradient $\det\igrad y^{-1}$ (the values are all positive at the grid points, numerical evidence of local orientation preservation and of no detected folding, which does not by itself establish that the computed deformation map $y$ is a diffeomorphism; the colormap has fixed limits $[0,2]$ (black: $0$; red: $1$, i.e., volume preserving; white: $2$) and clips values outside this range; since the sampled determinant is strictly positive in all results reported in this work, clipping can only occur at the upper end); and $(iv)$ computed deformation map $y$.}
\label{f:otgauss-256x256-h2s-exemplary-ngmres}
\end{figure}

We report convergence results for the \gangmres\ scheme in \Cref{t:otgauss-ngmres-new}. We consider various hyperparameter choices $w$ and $p = (\sigma,\tau)$. We show an exemplary result in \Cref{f:otgauss-256x256-h2s-exemplary-ngmres}.

\begin{table}
\caption{Convergence results for \gangmres\ for mapping a probability distribution $\pi_0$ to a distribution $\pi_1$ (optimal mass transport). Reported quantities and layout are as in \Cref{t:nirep-ngmres-sg}; runtimes with $\ast$ indicate that the scheme did not reach the tolerance. We use $n_{\scr{iter}} = 200$ and $\epsilon_{\scr{rel}} = 10^{-3}$. Runtimes are the median of five repetitions.}
\label{t:otgauss-ngmres-new}
\tabadjust\setcounter{run}{0}
\begin{tabular}{FrcrrrrRRRRRE}\toprule
\multicolumn{7}{c}{} & \multicolumn{6}{c}{\cellcolor{gray!20}time (in seconds)} \\
run & $w$ & $(\sigma,\tau)$ & \#iter & \#pdes & dist & grad & pdes & $q_k$ & $g$ & sg & ls & tts \\
\midrule
\rowcolor{p1color!20}
\runid &  10 &   (5,5) & 14 & 88 & \snum{2.495505e-03} & \snum{5.241528e-04} & \fnum{7.254811e+00} (\fnum{9.268168e-01}) & \fnum{4.773483e+00} & \fnum{2.008294e+00} & \fnum{8.269018e-01} & \fnum{2.106208e-02} & \fnum{7.827665e+00}\\
\runid &  15 &         & 14 & 89 & \snum{2.506658e-03} & \snum{8.969961e-04} & \fnum{7.308253e+00} (\fnum{9.252184e-01}) & \fnum{4.833343e+00} & \fnum{2.009484e+00} & \fnum{8.287617e-01} & \fnum{2.358938e-02} & \fnum{7.898950e+00}\\
\runid &  20 &         & 14 & 89 & \snum{2.506658e-03} & \snum{8.969961e-04} & \fnum{7.385016e+00} (\fnum{9.268824e-01}) & \fnum{4.882514e+00} & \fnum{2.035304e+00} & \fnum{8.430949e-01} & \fnum{2.384713e-02} & \fnum{7.967587e+00}\\
\runid &  25 &         & 14 & 89 & \snum{2.506658e-03} & \snum{8.969961e-04} & \fnum{7.262604e+00} (\fnum{9.256920e-01}) & \fnum{4.803082e+00} & \fnum{1.983392e+00} & \fnum{8.226670e-01} & \fnum{2.391108e-02} & \fnum{7.845594e+00}\\
\midrule
\runid &  10 &   (4,2) & 16 & 107 & \snum{2.497044e-03} & \snum{2.807262e-04} & \fnum{8.630600e+00} (\fnum{9.240797e-01}) & \fnum{5.606018e+00} & \fnum{2.299763e+00} & \fnum{1.169375e+00} & \fnum{3.270617e-02} & \fnum{9.339671e+00}\\
\runid &  15 &         & 19 & 128 & \snum{2.501882e-03} & \snum{6.454976e-04} & \fnum{1.008853e+01} (\fnum{9.196341e-01}) & \fnum{6.571820e+00} & \fnum{2.683555e+00} & \fnum{1.404341e+00} & \fnum{4.608250e-02} & \fnum{1.097016e+01}\\
\runid &  20 &         & 19 & 129 & \snum{2.504740e-03} & \snum{7.788904e-04} & \fnum{1.018446e+01} (\fnum{9.226944e-01}) & \fnum{6.529143e+00} & \fnum{2.680651e+00} & \fnum{1.474305e+00} & \fnum{4.550367e-02} & \fnum{1.103774e+01}\\
\runid &  25 &         & 19 & 129 & \snum{2.504740e-03} & \snum{7.788904e-04} & \fnum{1.020496e+01} (\fnum{9.211601e-01}) & \fnum{6.571868e+00} & \fnum{2.690968e+00} & \fnum{1.488898e+00} & \fnum{4.765613e-02} & \fnum{1.107838e+01}\\
\midrule
\runid &  10 &   (2,4) & 19 & 99 & \snum{2.494782e-03} & \snum{6.123988e-04} & \fnum{7.999793e+00} (\fnum{9.247893e-01}) & \fnum{5.188734e+00} & \fnum{2.713609e+00} & \fnum{5.767058e-01} & \fnum{1.853987e-02} & \fnum{8.650394e+00}\\
\runid &  15 &         & 19 & 99 & \snum{2.495188e-03} & \snum{4.090861e-04} & \fnum{8.199878e+00} (\fnum{9.222744e-01}) & \fnum{5.258860e+00} & \fnum{2.774132e+00} & \fnum{6.015210e-01} & \fnum{2.374492e-02} & \fnum{8.890931e+00}\\
\runid &  20 &         & 19 & 99 & \snum{2.472156e-03} & \snum{6.939374e-04} & \fnum{8.185683e+00} (\fnum{9.202443e-01}) & \fnum{5.181511e+00} & \fnum{2.809289e+00} & \fnum{5.988587e-01} & \fnum{2.526833e-02} & \fnum{8.895119e+00}\\
\runid &  25 &         & 19 & 99 & \snum{2.472156e-03} & \snum{6.939374e-04} & \fnum{7.646821e+00} (\fnum{9.191270e-01}) & \fnum{4.914781e+00} & \fnum{2.612752e+00} & \fnum{5.586981e-01} & \fnum{2.500404e-02} & \fnum{8.319656e+00}\\
\midrule
\runid &  10 &   (6,3) & 200 & 4597 & \snum{2.499871e-03} & \snum{1.004471e-03} & \fnum{3.748859e+02} (\fnum{9.801087e-01}) & \fnum{3.762588e+02} & \fnum{1.985070e+00} & \fnum{3.558060e+00} & \fnum{9.320533e-02} & $\ast$\fnum{3.824942e+02}\\
\runid &  15 &         & 13 & 91 & \snum{2.494819e-03} & \snum{9.527471e-04} & \fnum{7.640177e+00} (\fnum{9.270143e-01}) & \fnum{5.068240e+00} & \fnum{1.895229e+00} & \fnum{1.013969e+00} & \fnum{2.915871e-02} & \fnum{8.241704e+00}\\
\runid &  20 &         & 13 & 91 & \snum{2.494819e-03} & \snum{9.527471e-04} & \fnum{7.554037e+00} (\fnum{9.284956e-01}) & \fnum{4.996787e+00} & \fnum{1.891173e+00} & \fnum{1.011630e+00} & \fnum{2.424767e-02} & \fnum{8.135782e+00}\\
\runid &  25 &         & 13 & 91 & \snum{2.494819e-03} & \snum{9.527471e-04} & \fnum{7.455375e+00} (\fnum{9.274657e-01}) & \fnum{4.937799e+00} & \fnum{1.868075e+00} & \fnum{9.972985e-01} & \fnum{2.398575e-02} & \fnum{8.038438e+00}\\
\midrule
\runid &  10 &   (3,6) & 19 & 100 & \snum{2.500464e-03} & \snum{7.274186e-04} & \fnum{8.191401e+00} (\fnum{9.260699e-01}) & \fnum{5.255412e+00} & \fnum{2.758490e+00} & \fnum{5.901579e-01} & \fnum{1.774775e-02} & \fnum{8.845337e+00}\\
\runid &  15 &         & 19 & 100 & \snum{2.492856e-03} & \snum{7.074966e-04} & \fnum{8.096212e+00} (\fnum{9.219739e-01}) & \fnum{5.224820e+00} & \fnum{2.710865e+00} & \fnum{5.783293e-01} & \fnum{2.165442e-02} & \fnum{8.781390e+00}\\
\runid &  20 &         & 19 & 100 & \snum{2.475072e-03} & \snum{6.844089e-04} & \fnum{8.224907e+00} (\fnum{9.248727e-01}) & \fnum{5.281447e+00} & \fnum{2.763232e+00} & \fnum{5.908588e-01} & \fnum{2.289275e-02} & \fnum{8.893015e+00}\\
\runid &  25 &         & 19 & 100 & \snum{2.475072e-03} & \snum{6.844089e-04} & \fnum{8.098087e+00} (\fnum{9.243041e-01}) & \fnum{5.213690e+00} & \fnum{2.714270e+00} & \fnum{5.819087e-01} & \fnum{2.119950e-02} & \fnum{8.761280e+00}\\
\bottomrule
\end{tabular}
\end{table}

The most important observation is that the \gangmres\ scheme remains effective for the reformulation of our problem to account for a mass-preserving transport map. The solver remains effective for a range of hyperparameter values, providing excellent agreement between the transported density $\pi$ at $t=1$ and the target density $\pi_1$.
\section{Conclusions}\label{s:conclusions}

We have proposed a novel scheme to accelerate first-order optimization algorithms for PDE-constrained optimization problems governed by transport equations. We have compared it to the state-of-the-art~\cite{mang2015:inexact, mang2024:claire, mang2017:semilagrangian, mang2017:lagrangian, mang2016:distributed, mang2016:constrained, brunn2020:multinode, brunn2021:claire, brunn2021:fast, mang2019:claire} for intensity preserving transport maps, and we have tested it for mass-preserving transport maps (optimal transport) and incompressible flows of diffeomorphisms (Stokes flow). The most important observations are
\begin{enumerate*}[label=\it(\roman*)]
\item The proposed \gangmres\ scheme converges where the baseline \rpgd\ algorithm does not. For $\alpha = 10^{-3}$, \rpgd\ does not reach the tolerance within 200 iterations, whereas \gangmres\ does so in a median of 38 (\hands) and 93 (\nirep) iterations.
\item \gangmres\ does not, in general, outperform the (quasi-)Newton methods. The limited-memory BFGS method in the Hilbert space induced by the regularization model is faster by a factor of $1.63$ on \hands\ and $1.52$ on \nirep, by between $1.49$ and $1.80$ across the regularization parameters, and at the tighter tolerance of the mesh convergence study wherever it converges. Against \nk, the outcome depends on the preconditioner. \gangmres\ requires fewer PDE solves than \nk\ with the spectral and the two-level preconditioner (against the latter, between $1.0\times$ and $3.2\times$ fewer at all seven mesh points), and fewer than \nk\ with the zero-velocity preconditioner at six of the seven mesh points and, on \nirep, for $\alpha \ge 10^{-2}$; for $\alpha \le 10^{-3}$, including both image registration problems at our default $\alpha = 10^{-3}$, it requires more. In runtime, it is faster than the best \nk\ variant for $\alpha \ge 5\times 10^{-3}$ and slower for smaller values.
\item The proposed \gangmres\ scheme remains effective for a broad class of transport dominated PDE-constrained optimization problems, and it remains relatively insensitive to refinements in the discretization as long as the input data maintains the same smoothness level as the data presented on the coarsest mesh.
\item \gangmres\ is suited to moderately accurate solutions, as needed in the applications we target. On \nirep, tightening the tolerance increases its runtime less than that of \nk\ and \lbfgs, but below roughly $10^{-2}$ it no longer reaches the tolerance within $200$ iterations for any hyperparameter choice, whereas \nk\ and \lbfgs\ still converge. \nk\ remains preferable when high accuracy is needed on non-smooth data.
\item The proposed method is sensitive to vanishing regularization parameters. The largest eigenvalue of the preconditioned reduced space Hessian grows like $1/\alpha$ (see \Cref{s:appendix:theory}); the line search compensates by reducing the step size proportionally, so that residual components associated with small eigenvalues decrease only slowly. Addressing this sensitivity requires additional work.
\item If we select $w$ in $\{10,15,20,25\}$ \gangmres\ remains quite stable in terms of the time-to-solution with respect to changes in $p = (\sigma,\tau)$. Larger windows do not necessarily improve the runtime, and we recommend $\sigma\geq \tau$, which reduces the iteration count.
\item The \gangmres\ scheme significantly outperforms the \gaaa\ scheme, which reaches the tolerance for almost none of the hyperparameter choices on the two image registration problems.
\end{enumerate*}

In our future work, we plan to integrate the prototype implementation presented in this manuscript into our 3D graphics processing unit accelerated package. We also plan to explore how to address the sensitivity with respect to the regularization parameter $\alpha$, and to repeat the comparison with the Hilbert space quasi-Newton scheme in three dimensions, where the balance between the cost of storing secant pairs and that of applying the inverse Riesz map differs substantially from our two-dimensional prototype.

\section*{Acknowledgments}

This work was partly supported by the National Science Foundation (NSF) through the grants DMS-2145845 and DMS-2012825 (AM). Any opinions, findings, and conclusions or recommendations expressed herein are those of the authors and do not necessarily reflect the views of the NSF.

\appendix
\crefalias{section}{appendix}
\crefalias{subsection}{subappendix}

\section{Hardware Configuration and Additional Results}

\subsection{Hardware \& Software}\label{s:appendix:hardware}
All experiments were run in \texttt{MATLAB} \texttt{R2025b} under macOS Tahoe 26.5.2 on an Apple MacBook Pro (Mac16,8) with an Apple M4 Pro chip (12-core CPU with eight performance and four efficiency cores; 24 GB of unified memory).

\subsection{Proofs of the Convergence Results}\label{s:appendix:proofs}

We prove the results presented in \Cref{s:convergence}. Throughout, we mostly omit the subscripts of the norms: $\|\cdot\|$ denotes $\|\cdot\|_{\fs{V}}$ for elements of $\fs{V}$ (iterates and steps), $\|\cdot\|_{\fs{V}^\ast}$ for elements of $\fs{V}^\ast$ (gradients and their combinations), and the induced operator norm for operators, e.g., $\|\cdot\|_{\fs{V}\to\fs{V}^\ast}$ for $g'$ and $G$, and $\|\cdot\|_{\fs{V}\to\fs{V}}$ for $q'$. In particular, $e_k \defeq \|g(v^{(k)})\|_{\fs{V}^\ast}$.

\begin{proof}[Proof of \Cref{t:global}]
Let $s^{(k)} \defeq -(\alpha\tilde{\fs{L}})^{-1}g(v^{(k)})$, so that $q_k(v^{(k)})-v^{(k)} = \rho^{(k)}s^{(k)}$; dividing~\cref{e:qk-norm} by $\rho^{(k)}$ gives $\langle g(v^{(k)}), s^{(k)}\rangle_{L^2} = -e_k^2$ and $\|s^{(k)}\| = e_k$.   Since $g$ is the gradient of $\operatorname{obj}$ and Lipschitz continuous, $\operatorname{obj}(v+s) \le \operatorname{obj}(v) + \langle g(v), s\rangle_{L^2} + \frac{L_g}{2}\|s\|^2$ for $v$ and $v+s$ in the convex set. Let $e_k > 0$; otherwise, the algorithm stops at a zero of $g$. If the line search of Step~I rejects its first trial step, it has rejected $2\rho^{(k)}$, i.e., $\operatorname{obj}(v^{(k)}+2\rho^{(k)}s^{(k)}) \ge \operatorname{obj}(v^{(k)}) - 2c\rho^{(k)}e_k^2$; together with the previous inequality, this gives $-2\rho^{(k)} + 2L_g(\rho^{(k)})^2 \ge -2c\rho^{(k)}$, i.e., $\rho^{(k)} \ge (1-c)/L_g$. Otherwise, $\rho^{(k)}$ is the first trial step, which is $1$ for $k = 0$ and at least $\rho^{(k-1)}$ for $k \ge 1$ (see~\Cref{s:line-search}). By induction, $\rho^{(k)} \ge \rho_{\min} \defeq \min\{1,(1-c)/L_g\}$. Every iterate is either the \fp\ iterate $u^{(k)}$, which satisfies~\cref{e:armijo} at $v^{(k)}$ along $s^{(k)}$, or $u^{(k)} + \rho_{\scr{sg}}^{(k)}\hat{s}$ with $\langle g(u^{(k)}),\hat{s}\rangle_{L^2} < 0$ and $\rho_{\scr{sg}}^{(k)}$ satisfying~\cref{e:armijo} at $u^{(k)}$ along $\hat{s}$. In both cases,
\[
\operatorname{obj}(v^{(k+1)}) \le \operatorname{obj}(u^{(k)}) \le \operatorname{obj}(v^{(k)}) - c\rho^{(k)}e_k^2 \le \operatorname{obj}(v^{(k)}) - c\rho_{\min}e_k^2 .
\]
Summing over $k$ and rearranging gives $c\rho_{\min}\sum_k e_k^2 \le \operatorname{obj}(v^{(0)}) - \inf\operatorname{obj} < \infty$, so that $e_k \to 0$. The last statement of the theorem follows, since all norms on the finite-dimensional space $\fs{V}$ are equivalent. For the inexact gradient of the remark after \Cref{t:global}, we apply the descent inequality to $\nabla\operatorname{obj}$ instead of $g$: since $\langle\nabla\operatorname{obj}(v^{(k)}), s^{(k)}\rangle_{L^2} \le -e_k^2 + \|\nabla\operatorname{obj}(v^{(k)})-g(v^{(k)})\|\,\|s^{(k)}\| \le -(1-\epsilon_{\scr{g}})e_k^2$, a rejected trial step $2\rho^{(k)}$ gives $-2(1-\epsilon_{\scr{g}})\rho^{(k)} + 2L_g(\rho^{(k)})^2 \ge -2c\rho^{(k)}$, i.e., $\rho^{(k)} \ge (1-c-\epsilon_{\scr{g}})/L_g$. The rest of the proof is unchanged, since the Armijo conditions use the computed $g$, and $\|\nabla\operatorname{obj}(v^{(k)})\| \le (1+\epsilon_{\scr{g}})e_k \to 0$.
\end{proof}

\begin{proof}[Proof of \Cref{l:contraction}]
 Since $\langle\fs{B}v,w\rangle_{\fs{V}} = \langle\fs{H}v,w\rangle_{L^2}$, the operator $\fs{B}$ is self-adjoint and positive definite in $\fs{V}$; hence $q'(v^\ast) = \fs{I}-\rho\fs{B}$ is self-adjoint in $\fs{V}$ as well. The operator norm of $q'(v^\ast)$ on $\fs{V}$ therefore equals its spectral radius $\max\{|1-\rho\lambda_{\min}(\fs{B})|, |1-\rho\lambda_{\max}(\fs{B})|\}$, which is less than one if and only if $\rho\lambda_{\max}(\fs{B}) < 2$, since $0 < \lambda_{\min}(\fs{B}) \le \lambda_{\max}(\fs{B})$. (i) By continuity, $\|q'(v)\| \le \kappa < 1$ on a ball around $v^\ast$, on which $q$ is therefore a contraction. (ii) If $q$ were a contraction with constant $\kappa < 1$ in some norm $\|\cdot\|_\ast$ on a neighborhood of $v^\ast$, then $\|q'(v^\ast)d\|_\ast = \lim_{t\to 0^+}\|q(v^\ast+td)-q(v^\ast)\|_\ast/t \le \kappa\|d\|_\ast$ for all $d$,  so that the spectral radius of $q'(v^\ast)$ would not exceed $\kappa$. This contradicts the fact that $1-\rho\lambda_{\max}(\fs{B}) \le -1$ is an eigenvalue of $q'(v^\ast)$.
\end{proof}

\begin{proof}[Proof of \Cref{t:onestep}]
Write $a \defeq \zeta_{k+1}^{(k)}$, $b_j \defeq \zeta_j^{(k)}$, $x \defeq q_k(v^{(k)})$, and $z \defeq v^{(k+1)}$; below, $i$ and $j$ run over $\fs{W}_k$. For points $y$, $y'$ of the convex set, $g(y)-g(y') = G(y,y')(y-y')$ with $G(y,y') \defeq \int_0^1 g'(y'+t(y-y'))\,dt$, and~\cref{e:lipschitz} gives $\|G(z,y)-G(z,y')\|_{\fs{V}\to\fs{V}^\ast} \le \frac{\sigma_1}{2}\|y-y'\|$.  By~\cref{e:affine}, $z-x = \sum_i b_i(v^{(i)}-x)$ and $z-v^{(j)} = a(x-v^{(j)}) + \sum_i b_i(v^{(i)}-v^{(j)})$. Since $a+\sum_j b_j = 1$,
\begin{align*}
g(z)-\hat{g}^{(k+1)}
&= a\big(g(z)-g(x)\big) + \sum_j b_j\big(g(z)-g(v^{(j)})\big)\\
&= -a\sum_i b_i\,G(z,x)(x-v^{(i)}) + a\sum_j b_j\,G(z,v^{(j)})(x-v^{(j)}) + \sum_{i,j} b_ib_j\,G(z,v^{(j)})(v^{(i)}-v^{(j)})\\
&= a\sum_j b_j\big(G(z,v^{(j)})-G(z,x)\big)(x-v^{(j)}) + \sum_{i>j} b_ib_j\big(G(z,v^{(j)})-G(z,v^{(i)})\big)(v^{(i)}-v^{(j)}),
\end{align*}
where we have combined the terms $(i,j)$ and $(j,i)$ of the double sum. In the last line, each difference of two values of $G$ multiplies the difference of the corresponding points, so that the bound on $G$ above gives $\|g(z)-\hat{g}^{(k+1)}\|_{\fs{V}^\ast} \le \psi^{(k)}$,  and~\cref{e:ngmres-bound} follows from the triangle inequality and~\cref{e:theta-gamma}. For the accelerated candidate and the least-squares problem posed in~\cref{e:dualnorm}, $\hat{g}^{(k+1)}$ is the residual of the least-squares problem at its minimizer, for which $\beta^{(k)} = 0$ and $\beta^{(k)} = (-1,0,\ldots,0)$ are feasible; they give $g(x)$ and $g(v^{(k)})$, respectively, so that $\gamma^{(k+1)} \le 1$ and $\theta^{(k+1)} \le 1$. A \fp\ step gives $\hat{g}^{(k+1)} = g(x)$, i.e., $\gamma^{(k+1)} = 1$ and $\theta^{(k+1)} = \theta_{\scr{fp}}^{(k+1)}$. An iterate accepted in Step~III with $\rho_{\scr{sg}}^{(k)} \in (0,1)$ gives $\hat{g}^{(k+1)} = (1-\rho_{\scr{sg}}^{(k)})g(x) + \rho_{\scr{sg}}^{(k)}\hat{g}_{\scr{c}}$, where $\hat{g}_{\scr{c}}$ is the minimal residual of the least-squares problem in line~\ref{l:sg-lsq} of \Cref{a:gangmres-sg}, so that $\|\hat{g}_{\scr{c}}\|_{\fs{V}^\ast} \le \|g(x)\|_{\fs{V}^\ast}$ as above. Hence $\|\hat{g}^{(k+1)}\|_{\fs{V}^\ast} \le \|g(x)\|_{\fs{V}^\ast}$, i.e., $\gamma^{(k+1)} \le 1$ and $\theta^{(k+1)} \le \theta_{\scr{fp}}^{(k+1)}$.
\end{proof}

\noindent The proof of \Cref{t:local} rests on two lemmas, which we state under its assumptions.

\begin{lemma}\label[lemma]{l:errorbound}
Let $r_E \defeq \min\{r_0, 2/(\sigma_1c_E)\}$, with $r_E \defeq r_0$ if $\sigma_1 = 0$. Then $\|v-v^\ast\|_{\fs{V}} \le c_E\|g(v)\|_{\fs{V}^\ast}$ for all $v$ with $\|v-v^\ast\|_{\fs{V}} \le r_E$.
\end{lemma}

\begin{proof}
Let $d \defeq v-v^\ast$. By~\cref{e:lipschitz}, $\|g(v)-g'(v^\ast)d\| \le \frac{\sigma_1}{2}\|d\|^2$, and $\|g'(v^\ast)d\| \ge \|d\|/\|g'(v^\ast)^{-1}\| = 2\|d\|/c_E$. Hence $\|g(v)\| \ge (2/c_E - \sigma_1\|d\|/2)\|d\| \ge \|d\|/c_E$ for $\|d\| \le 2/(\sigma_1c_E)$.
\end{proof}

\begin{lemma}\label[lemma]{l:hot}
Let $v^{(j)}$, $j\in\fs{W}_k$, satisfy $\|v^{(j)}-v^\ast\|_{\fs{V}} \le r_E$, and let $M_k \defeq \max_{j\in\fs{W}_k}e_j$. Then $\psi^{(k)} \le C_h\bar\zeta^2M_k^2$ with $C_h \defeq \frac{\sigma_1}{4}(\bar\rho+2c_E)^2$.
\end{lemma}

\begin{proof}
By \Cref{l:errorbound}, $\|v^{(i)}-v^{(j)}\| \le c_E(e_i+e_j) \le 2c_EM_k$, and, by~\cref{e:qk-norm} and~\ref{as:rho}, $\|q_k(v^{(k)})-v^{(j)}\| \le \|q_k(v^{(k)})-v^{(k)}\| + \|v^{(k)}-v^{(j)}\| \le (\bar\rho+2c_E)M_k$. All distances in~\cref{e:hot} are thus bounded by $(\bar\rho+2c_E)M_k$. Their weights in~\cref{e:hot} sum to
\[
|\zeta_{k+1}^{(k)}|\sum_j|\zeta_j^{(k)}| + \sum_{i>j}|\zeta_i^{(k)}||\zeta_j^{(k)}|
\le |\zeta_{k+1}^{(k)}|\sum_j|\zeta_j^{(k)}| + \frac12\Big(\sum_j|\zeta_j^{(k)}|\Big)^2
\le \frac12\Big(|\zeta_{k+1}^{(k)}|+\sum_j|\zeta_j^{(k)}|\Big)^2 \le \frac12\bar\zeta^2
\]
by~\ref{as:zeta}. Hence $\psi^{(k)} \le \frac{\sigma_1}{2}\cdot\frac12\bar\zeta^2(\bar\rho+2c_E)^2M_k^2 = C_h\bar\zeta^2M_k^2$.
\end{proof}

\begin{proof}[Proof of \Cref{t:local}]
We assume $e_k > 0$, as~\ref{as:period} presupposes; otherwise, $v^{(k)} = v^\ast$ by \Cref{l:errorbound} and the induction below, and the iteration stops. Let $L_g \defeq \|g'(v^\ast)\| + \sigma_1r_0$, a Lipschitz constant of $g$ on $\fs{N}$, $\hat\kappa \defeq 1+L_g\bar\rho$, $K \defeq 2\hat\kappa^\omega$, and $s \defeq \lceil w/\omega\rceil$, where $\lceil\cdot\rceil$ and $\lfloor\cdot\rfloor$ denote rounding up and down to the nearest integer. Choose $\varepsilon > 0$ such that
\[
\eta \defeq C_h\bar\zeta^2K^2\hat\theta^{-2s}\varepsilon \le \min\Big\{\frac1\omega, \frac{\hat\theta-\bar\theta}{\omega\hat\kappa^\omega}\Big\}
\qquad\text{and}\qquad
\bar\zeta(\bar\rho+c_E)K\varepsilon \le r_E,
\]
and let $\delta \defeq \min\{r_E, \varepsilon/L_g\}$, so that $e_0 \le L_g\delta \le \varepsilon$. We show by induction over $n$ that
\begin{equation}\label{e:induction}
e_{n\omega} \le \hat\theta^ne_0,
\qquad
e_j \le K\hat\theta^{\lfloor j/\omega\rfloor}e_0
\;\;\text{and}\;\;
\|v^{(j)}-v^\ast\| \le r_E
\;\;\text{for all } j \le n\omega .
\end{equation}

\noindent For $n = 0$, this holds since $\delta \le r_E$. Assume~\cref{e:induction} for some $n$, and let $k = n\omega+l$, $0 \le l < \omega$, be such that the last two bounds in~\cref{e:induction} hold for all $j \le k$; for $l = 0$, this is the induction hypothesis. Each $j\in\fs{W}_k$ satisfies $j \ge k-w \ge (n-s)\omega$, so that $e_j \le K\hat\theta^{\max\{n-s,0\}}e_0$; hence $M_k \le K\varepsilon$ and $M_k^2 \le K^2\hat\theta^{n-2s}e_0^2$. By~\cref{e:qk-norm},~\ref{as:rho}, and \Cref{l:errorbound}, $\|q_k(v^{(k)})-v^\ast\| \le (\bar\rho+c_E)e_k \le r_E$, and, by~\cref{e:affine} and~\ref{as:zeta}, $\|v^{(k+1)}-v^\ast\| \le \bar\zeta(\bar\rho+c_E)M_k \le r_E$. All points that enter the step thus lie in the ball of radius $r_E$ around $v^\ast$, which is convex, so that \Cref{t:onestep} and \Cref{l:hot} apply and give $\psi^{(k)} \le C_h\bar\zeta^2M_k^2 \le \eta\hat\theta^ne_0$. Moreover, $\|g(q_k(v))\| \le \|g(v)\| + L_g\|q_k(v)-v\| \le \hat\kappa\|g(v)\|$ by~\cref{e:qk-norm}, so that $\theta^{(k+1)} \le \hat\kappa$ by~\cref{e:theta-fp}. Unrolling the one-step bounds~\cref{e:ngmres-bound-theta} of the steps $n\omega,\ldots,k$, each of which satisfies the bounds just derived, yields
\begin{equation}\label{e:unrolled}
e_{k+1} \le \Big(\prod_{m=n\omega}^{k}\theta^{(m+1)}\Big)e_{n\omega} + \sum_{m=n\omega}^{k}\Big(\prod_{m'=m+1}^{k}\theta^{(m'+1)}\Big)\psi^{(m)}
\le \hat\kappa^\omega(1+\omega\eta)\,\hat\theta^ne_0 \le K\hat\theta^ne_0,
\end{equation}

\noindent which completes the induction over $l$ for $k+1 < (n+1)\omega$. For $k+1 = (n+1)\omega$, \ref{as:period} bounds the first product in~\cref{e:unrolled} by $\bar\theta$, so that $e_{(n+1)\omega} \le (\bar\theta + \omega\hat\kappa^\omega\eta)\hat\theta^ne_0 \le \hat\theta^{n+1}e_0$, which completes the induction over $n$. Hence $e_k \le K\hat\theta^{\lfloor k/\omega\rfloor}e_0 \le (K/\hat\theta)\,\hat\theta^{k/\omega}e_0$ for all $k$, which is the first claim with $C \defeq K/\hat\theta$; the second follows from \Cref{l:errorbound}.
\end{proof}

\subsection{Verification of the Convergence Theory and the Role of the Hyperparameters}\label{s:appendix:theory}

We verify the statements of \Cref{s:convergence} numerically; \Cref{s:results:theory} summarizes the results. Unless stated otherwise, all numbers are for the \nirep\ dataset without the safeguard. The measurements of $\theta^{(k+1)}$, $\gamma^{(k+1)}$, and the factors in~\cref{e:theta-fp} use the least-squares problem posed in~\cref{e:dualnorm}; all other numbers are for the Euclidean variant of \Cref{s:results}.

\subsubsection{Contractivity}

For each regularization parameter, we take the solution computed in the regularization study by \gangmres\ with $w=20$, $p=(5,1)$ and the safeguard; these runs are reported in \Cref{t:nirep-ngmres-regul-new}. At this solution, we estimate the largest eigenvalue of $\fs{B} = (\alpha\tilde{\fs{L}})^{-1}\fs{H}$ with a Krylov eigenvalue solver, where $\fs{H}$ is the Gauss--Newton approximation of the reduced space Hessian, both on all fields and on mean-free fields, and we record how the energy of the top eigenvector is distributed over the Fourier modes. In \Cref{t:theory-contraction}, we compare the eigenvalue with the step size $\rho$ accepted by the line search of Step~I.

\begin{table}
\caption{We report the largest eigenvalue of $\fs{B} = (\alpha\tilde{\fs{L}})^{-1}\fs{H}$ at the computed solution, where $\fs{H}$ is the Gauss--Newton approximation of $g'(v^\ast)$, the step size $\rho$ accepted there by the line search of Step~I, and their product, for the regularization parameters $\alpha$ we consider (\nirep; \gangmres\ with $w=20$, $p=(5,1)$ and the safeguard switched on). The columns labeled mean-free report the same quantities for $\fs{B}$ restricted to mean-free fields; the last two columns report the fraction of the energy of the top eigenvector in the mean ($\xi = 0$) and in the modes with $|\xi| = 1$. The contraction condition~\cref{e:contraction} requires $\rho\lambda_{\max}(\fs{B}) < 2$.}\label{t:theory-contraction}
\tabadjust
\begin{tabular}{rrrrrrrr}\toprule
& & & & \multicolumn{2}{c}{mean-free} & \multicolumn{2}{c}{energy of top eigenvector} \\
\cmidrule(lr){5-6}\cmidrule(lr){7-8}
$\alpha$ & $\lambda_{\max}(\fs{B})$ & $\rho$ & $\rho\lambda_{\max}(\fs{B})$ & $\lambda_{\max}$ & $\rho\lambda_{\max}$ & $\xi = 0$ & $|\xi| = 1$ \\\midrule
$10^{-1}$            &   $39.7$ & $1.25\times 10^{-1}$ & $4.96$  &   $35.1$ & $4.39$  & $0.32$ & $0.67$ \\
$5\times 10^{-2}$    &   $79.2$ & $6.25\times 10^{-2}$ & $4.95$  &   $69.7$ & $4.36$  & $0.33$ & $0.66$ \\
$10^{-2}$            &  $397.8$ & $3.13\times 10^{-2}$ & $12.43$ &  $346.7$ & $10.84$ & $0.34$ & $0.65$ \\
$5\times 10^{-3}$    &  $795.9$ & $3.91\times 10^{-3}$ & $3.11$  &  $692.7$ & $2.71$  & $0.35$ & $0.64$ \\
$10^{-3}$            & $3985.7$ & $4.88\times 10^{-4}$ & $1.95$  & $3466.4$ & $1.69$  & $0.35$ & $0.64$ \\
$5\times 10^{-4}$    & $7991.7$ & $4.88\times 10^{-4}$ & $3.90$  & $6951.7$ & $3.39$  & $0.35$ & $0.64$ \\
\bottomrule
\end{tabular}
\end{table}

For the Gauss--Newton approximation, condition~\cref{e:contraction} is violated at five of the six regularization parameters; at $\alpha = 10^{-3}$, the accepted step happens to fall just below the threshold. The largest eigenvalue scales like $1/\alpha$ ($\alpha\lambda_{\max}(\fs{B})$ lies between $3.96$ and $4.00$), whereas the accepted step size, a power of two between $2^{-3}$ and $2^{-11}$, depends on the decrease of the objective along the search direction and decreases only roughly in proportion to $\alpha$ ($\rho/\alpha$ between $0.5$ and $3.1$). A step of order $\alpha$ therefore barely reduces the residual components along eigenvalues of $\fs{B}$ of order one, consistent with the deterioration of the scheme as $\alpha$ decreases. The top eigenvector holds $32$ to $35\%$ of its energy in the mean and $64$ to $67\%$ in the modes with $|\xi| = 1$; on mean-free fields, which the modification of the kernel of $\fs{L}$ (see \Cref{s:discretization}) does not affect, $\lambda_{\max}(\fs{B})$ is only $11$ to $13\%$ smaller, and the condition fails at the same five parameters. The iteration still makes progress because the line search enforces descent, which is a weaker requirement than contraction. These numbers concern the Gauss--Newton approximation, which neglects the second-order terms of the data mismatch. Since \Cref{l:contraction} concerns the exact Jacobian $g'(v^\ast)$ and these terms can change its eigenvalues, the numbers indicate, but do not establish, that the maps $q_k$ are not contractive; establishing it would require the exact Jacobian or a bound on the neglected terms.

\subsubsection{Bounds on NGMRES Steps}
We run the variant with the least-squares problem posed in~\cref{e:dualnorm} without the safeguard on \nirep\ ($\alpha = 10^{-3}$) for twelve hyperparameter settings: $p=(5,1)$ with $w\in\{1,5,10,15,20,25,50\}$, $w=20$ with $p\in\{(4,2),(2,4)\}$, and $p=(1,0)$ with $w\in\{5,20,25\}$. All twelve runs reach the tolerance. \Cref{t:theory-theta} reports the quantities of this section for each of them. On \ngmres\ steps, $\theta^{(k+1)}$ and $\gamma^{(k+1)}$, measured in the dual norm $\|\cdot\|_{\fs{V}^\ast}$, never exceed one, as they must by construction. Their maxima over the twelve runs are $0.9982$ and $0.9999$, both attained for $w=1$, $p=(5,1)$. If we evaluate them along the same trajectories in the Euclidean norm instead, at least one of them exceeds one in eleven of the twelve runs, by up to $9\%$; the exception is $w=20$, $p=(2,4)$. The median of $\theta^{(k+1)}$ lies at or below the observed rate of reduction of $\|g\|_{\fs{V}^\ast}$ in all twelve runs, but it does not predict that rate.

\begin{table}
\caption{Quantities of the convergence theory on \ngmres\ steps for the twelve runs of the variant with the least-squares problem posed in~\cref{e:dualnorm} (\nirep; $\alpha = 10^{-3}$; without the safeguard). We report the number of iterations, the median (med.) and the maximum of $\theta = \theta^{(k+1)}$, the maximum of $\gamma = \gamma^{(k+1)}$, the angle in~\cref{e:gamma-angle} that corresponds to the median of $\gamma^{(k+1)}$, the median of $\sum_j|\zeta_j^{(k)}|$, and the median observed rate of reduction of $\|g\|_{\fs{V}^\ast}$ per iteration, all measured in the dual norm $\|\cdot\|_{\fs{V}^\ast}$. The last two columns report the maxima of $\theta^{(k+1)}$ and $\gamma^{(k+1)}$ along the same trajectories when they are measured in the Euclidean norm. The entry $>0.9999$ is $0.999999$; it rounds to one at four digits but lies below one.}\label{t:theory-theta}
\tabadjust
\begin{tabular}{rcrrrrrrrrr}\toprule
& & & \multicolumn{6}{c}{dual norm} & \multicolumn{2}{c}{Euclidean norm} \\
\cmidrule(lr){4-9}\cmidrule(lr){10-11}
$w$ & $p$ & \#iter & med.\ $\theta$ & max $\theta$ & max $\gamma$ & angle & med.\ $\sum_j|\zeta_j|$ & med.\ rate & max $\theta$ & max $\gamma$ \\\midrule
 1 & $(5,1)$ & 191 & $0.9832$ & $0.9982$ & $>0.9999$ & $68.0^\circ$ & $  2.21$ & $0.9904$ & $1.0254$ & $1.0277$ \\
 5 & $(5,1)$ & 104 & $0.9498$ & $0.9838$ & $0.9961$ & $53.4^\circ$ & $  3.17$ & $0.9588$ & $1.0015$ & $1.0154$ \\
10 & $(5,1)$ &  95 & $0.9447$ & $0.9816$ & $0.9917$ & $52.9^\circ$ & $  4.70$ & $0.9528$ & $1.0194$ & $1.0261$ \\
15 & $(5,1)$ &  81 & $0.9042$ & $0.9855$ & $0.9938$ & $51.6^\circ$ & $  5.67$ & $0.9463$ & $1.0194$ & $1.0248$ \\
20 & $(5,1)$ &  79 & $0.7207$ & $0.9705$ & $0.9827$ & $38.6^\circ$ & $ 10.06$ & $0.9355$ & $1.0194$ & $1.0248$ \\
25 & $(5,1)$ &  70 & $0.6407$ & $0.9691$ & $0.9804$ & $31.7^\circ$ & $ 15.37$ & $0.9057$ & $1.0194$ & $1.0248$ \\
50 & $(5,1)$ &  63 & $0.2738$ & $0.8323$ & $0.8845$ & $18.4^\circ$ & $ 54.89$ & $0.8871$ & $1.0462$ & $1.0462$ \\
\midrule
20 & $(4,2)$ &  75 & $0.6300$ & $0.9783$ & $0.9855$ & $41.1^\circ$ & $ 18.15$ & $0.9327$ & $1.0132$ & $1.0206$ \\
20 & $(2,4)$ &  91 & $0.4930$ & $0.8730$ & $0.9179$ & $26.5^\circ$ & $169.07$ & $0.9807$ & $0.9828$ & $0.9870$ \\
\midrule
 5 & $(1,0)$ & 103 & $0.9581$ & $0.9813$ & $0.8845$ & $8.9^\circ$ & $  2.70$ & $0.9581$ & $0.9956$ & $1.0541$ \\
20 & $(1,0)$ &  89 & $0.9069$ & $0.9793$ & $0.8997$ & $11.1^\circ$ & $  8.19$ & $0.9525$ & $1.0756$ & $1.0902$ \\
25 & $(1,0)$ &  75 & $0.6258$ & $0.9834$ & $0.8880$ & $20.8^\circ$ & $ 16.70$ & $0.9439$ & $1.0756$ & $1.0902$ \\
\bottomrule
\end{tabular}
\end{table}

For the alternating sequences, the observed rate includes the \fp\ steps, which increase the residual, as \Cref{t:theory-cycle} shows. For $p=(1,0)$, which has no \fp\ steps, the gap grows with the coefficients: it vanishes for $w=5$ (median $\sum_j|\zeta_j^{(k)}| = 2.7$) and is $4.6$ percentage points for $w=20$ ($8.2$) and $32$ for $w=25$ ($16.7$), consistent with the higher order terms of~\cref{e:ngmres-bound}. The bound concerns $\|g\|_{\fs{V}^\ast}$ and does not transfer verbatim to our stopping criterion, which uses the $\ell^\infty$-norm of the reduced gradient. The setting $w=25$, $p=(1,0)$, which stagnates in the Euclidean variant in the sweep without the safeguard of \Cref{s:results:safeguard}, converges in this variant.

\subsubsection{Alternating Scheme}

For the four runs with $w=20$ of \Cref{t:theory-theta}, \Cref{t:theory-cycle} reports the medians of the factors in~\cref{e:theta-fp}, measured in the dual norm, on \ngmres\ steps ($\theta_{\scr{ng}}$; first case) and on \fp\ steps ($\theta_{\scr{fp}}$; second case). Since the linear factor of a period is the product of its $\sigma+\tau$ factors, we also report the average factor per period, the largest product of the factors over a completed period, which is the quantity that~\ref{as:period} bounds, and the \fp-only factor, as defined in the caption.

\begin{table}
\caption{Factors of~\cref{e:theta-fp}, measured in the dual norm: median of $\theta_{\scr{ng}}$ over the \ngmres\ steps and of $\theta_{\scr{fp}}$ over the \fp\ steps of each run, the average factor per period $\big(\prod_{k=0}^{K-1}\theta^{(k+1)}\big)^{(\sigma+\tau)/K}$ with $K$ the number of iterations (per period), the largest product of the factors over a completed period (max per period), and the geometric mean of the $\theta_{\scr{fp}}$ alone raised to the power $\sigma+\tau$ (\fp\ only). The runs are the four with $w=20$ of \Cref{t:theory-theta}; the sequence $p=(1,0)$ has no \fp\ steps (---).}\label{t:theory-cycle}
\tabadjust
\begin{tabular}{lrrrrrr}\toprule
configuration & \#iter & median $\theta_{\scr{ng}}$ & median $\theta_{\scr{fp}}$ & avg.\ per period & max per period & \fp\ only \\\midrule
$w=20$, $p=(5,1)$ &  79 & $0.7207$ & $1.6995$ & $0.1581$ & $2.5841$ & $40.4070$ \\
$w=20$, $p=(4,2)$ &  75 & $0.6300$ & $1.1993$ & $0.1291$ & $1.0592$ & $3.9852$ \\
$w=20$, $p=(2,4)$ &  91 & $0.4930$ & $0.9980$ & $0.2422$ & $0.8758$ & $1.1729$ \\
$w=20$, $p=(1,0)$ &  89 & $0.9069$ & ---      & $0.6059$ & $0.9793$ & ---      \\
\bottomrule
\end{tabular}
\end{table}

As required, $\theta_{\scr{ng}}$ never exceeds one. The median of $\theta_{\scr{fp}}$ exceeds one for $p=(5,1)$ and $p=(4,2)$ and is close to one for $p=(2,4)$, consistent with the noncontractivity suggested above. The average factor per period is far below the \fp-only factor, as expected from the local comparison in \Cref{s:convergence}; this does not by itself show that \gangmres\ is faster than the \fp\ iteration. The largest product, however, exceeds one for $p=(5,1)$ and $p=(4,2)$ ($2.58$ and $1.06$), so that these two runs converge although they do not satisfy~\ref{as:period}.

\subsubsection{Choice of Norm and Window Size}
We also ran the dual-norm variant on all four hyperparameter sweeps of \Cref{s:results:conv} and studied the effect of the window size; for brevity, we only summarize the results. With the dual norm, up to four fewer of the $86$ settings per sweep reach the tolerance, and the settings that converge with either norm need $0.97$ to $1.04$ times as many PDE solves in the median. A larger window decreases the angle in~\cref{e:gamma-angle} but increases $\sum_j|\zeta_j^{(k)}|$, and the iteration count follows neither trend. Without the safeguard, the \fp\ steps appear to stabilize the iteration; on \nirep, two of the eight settings with $\tau = 0$ miss the tolerance, against three of the $78$ with $\tau > 0$.

\end{document}